\documentclass[a4paper]{article} 

\usepackage{latexsym}
\usepackage{amsfonts}
\usepackage{amssymb}
\usepackage{amsmath}
\usepackage{amscd}
\usepackage{eucal}
\usepackage[all, knot]{xy}
\xyoption{arc}
\xyoption{rotate}
\xyoption{dvips}
\usepackage{amsthm}
\usepackage{pict2e}

\usepackage{pstricks}
\usepackage{pst-node}
\usepackage{verbatim} 
\usepackage{graphics}
\usepackage{graphicx}
\usepackage{color}   
\usepackage{framed}

\usepackage{bbm}

\usepackage{marvosym}  
\usepackage{stmaryrd}
\usepackage{textcomp}
\usepackage{wasysym}

\renewcommand{\arraystretch}{1.3} 

\psset{unit=1mm,linewidth=0.5pt,
arrowlength=1.1, arrowinset=0.7,arrowsize=4.5pt,
doublesep=1pt,
labelsep=2pt,nodesep=2pt}
  
\def\1{\ensuremath{\mathbbm{1}}}%
\def\2{\ensuremath{\mathbbm{2}}}%

\newcommand{\E}{\ensuremath{\bb{E}}}

\newcommand{\cB}{{\cl{B}}}

\newcommand{\cE}{{\cl{E}}}

\newcommand{\cJ}{{\cl{J}}}

\renewcommand{\bf}{\bfseries}

\newcommand{\iso}{\cong} 
\newcommand{\catequiv}{\simeq} 
\newcommand{\lra}{\ensuremath{\longrightarrow}} 

\renewcommand{\:}{\colon}

\newcommand{\cat}[1]{\ensuremath{\textrm{\bfseries {\upshape {#1}}}}}

\newcommand{\ob}{\mbox{\upshape{ob}\hspace{1.3pt}}}

\newcommand{\Set}{{\cat{Set}}}
\newcommand{\Cat}{{\cat{Cat}}}

\newcommand{\Bicon}{{\cat{$\cB$-Icon}}}

\newcommand{\cl}[1]{\ensuremath{\mathcal {#1}}}
\newcommand{\bb}[1]{\ensuremath{\mathbb {#1}}}
\newcommand{\ed}{\end{document}}
\newcommand{\map}[1]{\ensuremath{\stackrel{{#1}}{\lra}}}

\newcommand{\demph}[1]{{\bfseries #1}}

\newcommand{\hh}[1]{\hspace*{#1}}

\newcommand{\tra}{\psset{unit=0.1cm,nodesep=0pt} \pspicture(8,0)
\pcline{->}(1,1.1)(7,1.1) \endpspicture}

\newcommand{\tmapsto}{\psset{unit=0.1cm,nodesep=0pt} \pspicture(8,0) 
\pcline{|->}(1,1.2)(7,1.2) \endpspicture}

\newcommand{\tramap}[1]{\psset{unit=0.1cm,nodesep=0pt,labelsep=2pt} \pspicture(8,4)
\pcline{->}(1,1.1)(7,1.1)\naput{\ensuremath{\scriptstyle{#1}}} \endpspicture}

\newcommand{\tmap}{\tramap}

\newcommand{\ltramap}[1]{\psset{unit=0.1cm,nodesep=0pt} \pspicture(15,4)
\pcline{->}(1.5,1.1)(13.5,1.1)\naput{\ensuremath{\scriptstyle{#1}}} \endpspicture}

\newcommand{\ltmap}{\ltramap}

\newcommand{\trta}{\psset{unit=0.1cm,nodesep=0pt} \pspicture(8,0)
\pcline[doubleline=true,arrowinset=0.6,arrowlength=0.8]{->}(1,1.1)(7,1.1) \endpspicture}

\newcommand{\Tra}{\trta}

\newcommand{\noi}{\noindent}
\newcommand{\dend}{\ensuremath{\hfill \maltese}\end{mydefinition}}

\newcommand{\igc}[1]{\begin{tabular}[c]{c}\includegraphics[width={#1}]}
\newcommand{\igt}[1]{\begin{tabular}[t]{c}\\[-24pt] \includegraphics[width={#1}]}
\newcommand{\ei}{\end{tabular}}

\newcommand{\setunit}[1]{\setlength{\unitlength}{#1}}

\newcommand{\degenarrow}{
\begin{picture}(0,0)
\qbezier[20](0,0)(-3,-2.25)(-6,-4.5)
\put(-6,-4.5){\vector(-12,-9){0.1}}
\end{picture}}

\theoremstyle{plain}

\newtheorem{theorem}{Theorem}[section]

\newtheorem{proposition}[theorem]{Proposition}

\newtheorem{conjecture}[theorem]{Conjecture}

\theoremstyle{definition}

\newtheorem{definition}[theorem]{Definition}

\newtheorem{example}[theorem]{Example}

\newtheorem{examples}[theorem]{Examples}
\newtheorem{nonexample}[theorem]{Non-example}
\newtheorem{remark}[theorem]{Remark}
\newtheorem{remarks}[theorem]{Remarks}
\newtheorem{exercise}[theorem]{Exercise}
\newtheorem{note}[theorem]{Note}
\newtheorem{question}[theorem]{Question}
\newtheorem{questions}[theorem]{Questions}
\newtheorem{algorithm}[theorem]{Algorithm}
\newtheorem{method}[theorem]{Method}

\newenvironment{mydefinition}{\begin{definition}} {\end{definition}}

\newenvironment{prf}{\vspace{1ex}\begin{sloppypar}{\noindent\upshape
{\bfseries Proof.}}} {{\hspace*{\fill}
$\Box$}\end{sloppypar}\vspace{2ex}}

{ \end{sf}\end{framed}\end{minipage}
\end{center}}

\newcommand{\bs}{\boldsymbol}

\renewcommand{\dot}{\centerdot}

\begin{document}

\title{Iterated icons}

\author{Eugenia Cheng and Nick Gurski \\  Department of Mathematics, University of Sheffield \\E-mail: e.cheng@sheffield.ac.uk, nick.gurski@sheffield.ac.uk}


\maketitle

\begin{abstract}
We study the totality of categories weakly enriched in a monoidal bicategory using a notion of enriched icon as 2-cells.  We show that when the monoidal bicategory in question is symmetric then this process can be iterated.  We show that starting from the symmetric monoidal bicategory \cat{Cat} and performing the construction twice yields a convenient symmetric monoidal bicategory of partially strict tricategories.  We show that restricting to the doubly degenerate ones immediately gives the correct bicategory of ``2-tuply monoidal categories'' missing from our earlier studies of the Periodic Table.  We propose a generalisation to all $k$-tuply monoidal $n$-categories.
\end{abstract}

\section*{Introduction}

In this paper we use an iterated icon construction to study the degenerate $n$-categories that appear in the Periodic Table. We began our analysis of the Periodic Table of $n$-categories in \cite{cg1} (degenerate categories and bicategories) and \cite{cg3} (degenerate tricategories).  A $k$-degenerate $n$-category is one in which the lowest $k$ dimensions are degenerate, that is, contain only one cell.  In this case the lowest non-trivial dimension is $k$, and we can perform a ``dimension shift'' to forget the degenerate dimensions and regard the $k$-cells as the 0-cells of an $(n-k)$-category.  These cells come equipped with $k$ monoidal structures coming from the $k$ types of composition they had in the old structure; thus a degenerate category is regarded as a monoid, and a degenerate bicategory as a monoidal category.  In general we have the following shift:
\begin{center}
\begin{tabular}{llcl}
 & $k$-degenerate $m$-category & $\tra$ & $k$-monoidal $(m-k)$-category \\
or & $k$-degenerate $(n+k)$-category & $\tra$ & $k$-monoidal $n$-category.
\end{tabular}
\end{center}

The first type of indexing is more natural if one is starting with a fixed $m$ and varying the amount of degeneracy as we did in \cite{cg1,cg3}.  The second form is more appropriate if one is fixing the \emph{codimension} and varying the number of monoidal structures, which is what we will end up doing in the present work.

The periodic table of $n$-categories as proposed by Baez and Dolan \cite{bd3} gives predictions of what sort of structures should arise as $k$-monoidal $n$-categories.  The idea is that as the $k$ monoidal structures came from different types of composition in the original $(n+k)$-category, the interchange laws from that structure give interactions between the different monoidal structures which, as a result, amount to weak ``commutations'' such as braidings and so on.  The first few entries are given in the following table.

\scalebox{0.8}{
\setunit{2mm}  
\begin{picture}(80,92)(12,-6)

\put(0,36.5){
\renewcommand{\tabcolsep}{1.3em}
\renewcommand{\arraystretch}{1}

\begin{tabular}{|lllll|} \hline &&&&\\
set & category & 2-category & 3-category & \makebox(0,0)[r]{$\cdots$}  \\[40pt]
monoid & monoidal category & monoidal 2-category & monoidal 3-category & \makebox(0,0)[r]{$\cdots$} \\
\makebox(2.5,0.5)[l]{$\equiv$}category with & \makebox(2.5,0.5)[l]{$\equiv$}2-category with & \makebox(2.5,0.5)[l]{$\equiv$}3-category with & \makebox(2.5,0.5)[l]{$\equiv$}4-category with &\\ 
\makebox(2.5,0.5)[l]{}only one object & \makebox(2.5,0.5)[l]{}only one object & \makebox(2.5,0.5)[l]{}only one object & \makebox(2.5,0.5)[l]{}only one object &\\[35pt]
{\bfseries commutative}& braided monoidal & braided monoidal & braided monoidal & \makebox(0,0)[r]{$\cdots$}  \\
{\bfseries monoid}& category & 2-category & 3-category &\\
\makebox(2.5,0.5)[l]{$\equiv$}2-category with & \makebox(2.5,0.5)[l]{$\equiv$}3-category with & \makebox(2.5,0.5)[l]{$\equiv$}4-category with & \makebox(2.5,0.5)[l]{$\equiv$}5-category with &\\ 
\makebox(2.5,0.5)[l]{}only one object & \makebox(2.5,0.5)[l]{}only one object & \makebox(2.5,0.5)[l]{}only one object & \makebox(2.5,0.5)[l]{}only one object &\\
\makebox(2.5,0.5)[l]{}only one 1-cell & \makebox(2.5,0.5)[l]{}only one 1-cell & \makebox(2.5,0.5)[l]{}only one 1-cell & \makebox(2.5,0.5)[l]{}only one 1-cell &\\[30pt]
\makebox(8,0)[r]{\large $\prime\prime$}& {{\bfseries symmetric monoidal}} & sylleptic monoidal & {sylleptic monoidal} & \makebox(0,0)[r]{$\cdots$}  \\
& {\bfseries category} & 2-category & 3-category &\\
\makebox(2.5,0.5)[l]{$\equiv$}3-category with & \makebox(2.5,0.5)[l]{$\equiv$}4-category with & \makebox(2.5,0.5)[l]{$\equiv$}5-category with & \makebox(2.5,0.5)[l]{$\equiv$}6-category with &\\ 
\makebox(2.5,0.5)[l]{}only one object & \makebox(2.5,0.5)[l]{}only one object & \makebox(2.5,0.5)[l]{}only one object & \makebox(2.5,0.5)[l]{}only one object &\\
\makebox(2.5,0.5)[l]{}only one 1-cell & \makebox(2.5,0.5)[l]{}only one 1-cell & \makebox(2.5,0.5)[l]{}only one 1-cell & \makebox(2.5,0.5)[l]{}only one 1-cell &\\
\makebox(2.5,0.5)[l]{}only one 2-cell & \makebox(2.5,0.5)[l]{}only one 2-cell & \makebox(2.5,0.5)[l]{}only one 2-cell & \makebox(2.5,0.5)[l]{}only one 2-cell &\\ [20pt]
\makebox(8,0)[r]{\large $\prime\prime$}& \makebox(8,0)[r]{\large $\prime\prime$} & \makebox(16,0)[l]{{\bfseries symmetric monoidal}} & \makebox(8,0)[r]{\Large ?}& \makebox(0,0)[r]{$\cdots$}  \\
& & {\bfseries 2-category} & &\\
\makebox(2.5,0.5)[l]{$\equiv$}4-category with & \makebox(2.5,0.5)[l]{$\equiv$}5-category with & \makebox(2.5,0.5)[l]{$\equiv$}6-category with & \makebox(2.5,0.5)[l]{$\equiv$}7-category with &\\ 
\makebox(2.5,0.5)[l]{}only one object & \makebox(2.5,0.5)[l]{}only one object & \makebox(2.5,0.5)[l]{}only one object & \makebox(2.5,0.5)[l]{}only one object &\\
\makebox(2.5,0.5)[l]{}only one 1-cell & \makebox(2.5,0.5)[l]{}only one 1-cell & \makebox(2.5,0.5)[l]{}only one 1-cell & \makebox(2.5,0.5)[l]{}only one 1-cell &\\
\makebox(2.5,0.5)[l]{}only one 2-cell & \makebox(2.5,0.5)[l]{}only one 2-cell & \makebox(2.5,0.5)[l]{}only one 2-cell & \makebox(2.5,0.5)[l]{}only one 2-cell &\\ 
\makebox(2.5,0.5)[l]{}only one 3-cell & \makebox(2.5,0.5)[l]{}only one 3-cell & \makebox(2.5,0.5)[l]{}only one 3-cell & \makebox(2.5,0.5)[l]{}only one 3-cell &\\[30pt] 
\makebox(8,0)[r]{\large $\prime\prime$}& \makebox(8,0)[r]{\large $\prime\prime$}& \makebox(8,0)[r]{\large $\prime\prime$}& {\bfseries symmetric monoidal}& \makebox(0,0)[r]{$\cdots$} \\
&&& {\bfseries 3-category} &\\
\makebox(7.5,0)[r]{$\vdots$} & \makebox(7.5,0)[r]{$\vdots$}
& \makebox(7.5,0)[r]{$\vdots$} & \makebox(7.5,0)[r]{$\vdots$} & \\
&&&& \\
&&&& \\
 
\hline
\end{tabular}}

\put(22,76){\degenarrow}
\put(22,61){\degenarrow}
\put(22,46){\degenarrow}
\put(22,30){\degenarrow}
\put(22,12){\degenarrow}

\put(43.5,76){\degenarrow}
\put(43.5,61){\degenarrow}
\put(43.5,46){\degenarrow}
\put(43.5,30){\degenarrow}
\put(43.5,12){\degenarrow}

\put(64,76){\degenarrow}
\put(64,61){\degenarrow}
\put(64,46){\degenarrow}
\put(64,30){\degenarrow}
\put(64,12){\degenarrow}


\end{picture}}

\renewcommand{\arraystretch}{1.3}

\vspace{3em}

The process of fixing an overall dimension $m$ and increasing the amount of degeneracy amounts to moving diagonally down-and-left in the table.  The process of fixing codimension $n$ and adding monoidal structures consists of moving down a column.  The hypothesis is then that there is some kind of equivalence between the following three types of structure:

\begin{itemize}
\item $k$-degenerate $n$-categories
\item $k$-tuply monoidal $(n-k)$-categories
\item the $(n,k)$th entry of the Periodic Table.
\end{itemize}

One of the lessons learnt from \cite{cg1,cg3} is that we must take care about the precise formulation of these equivalences.  For example, we wish to compare the totality of commutative monoids with the totality of doubly degenerate bicategories, but the former naturally form a category whereas the latter form a tricategory.  Now, bicategories do form a category but this is a sort of ``accident'' and does not generalise into higher dimensions; $n$-categories and weak maps do not form a category for any higher $n$ (for algebraic definitions).  Moreover, the category of degenerate bicategories is in any case not equivalent to the category of commutative monoids, as shown in \cite{cg1}.  

The general problem is that we are trying to compare a totality of $k$-degenerate $(n+k)$-categories, which is naturally an $(n+k+1)$-dimensional structure, with a totality of $k$-monoidal $n$-categories, which is naturally an $(n+1)$-dimensional structure.  

In \cite{cg1} we investigated ways of 
\begin{itemize}
 \item decreasing the dimensions of the former, by forgetting higher-dimensional maps, or
\item increasing the dimensions of the latter, by adding identities as higher cells.
\end{itemize}
We did this in a somewhat ad hoc way on a case-by-case basis, with mixed results.

We are now able to propose a much more systematic approach based on Lack's theory of icons \cite{lac2,lp1}.  The aim of Lack's work was to construct a ``convenient 2-category of bicategories''.  It is well-known that bicategories and homomorphisms form a category, and that bicategories, homomorphisms, transformations and modifications form a tricategory, but simply forgetting the modifications does not give us a 2-category as the transformations do not compose coherently enough.

Lack's solutions is to restrict to ``Identity Component Oplax Natural transformations'' (ICONs), and re-arrange their composition a little.  This gives a 2-dimensional totality of bicategories that is not only coherent (i.e. it is a 2-category) but also convenient---it contains enough 2-dimensional information to express a strong form of the coherence theorem for bicategories.

The key, for our purposes, is that it gives the correct totality for comparing degenerate bicategories with monoidal categories which, after all, naturally form a bicategory.  That is, we consider $\cat{Icon}$, the 2-category of bicategories, weak functors and icons.  The full sub-2-category on the degenerate bicategories is then equivalent to the 2-category of monoidal categories.  This is done by Lack \cite{lac2}.

We now generalise this construction in order to iterate it.  In \cite{gg1} various convenient lower-dimensional totalities of tricategories are constructed (that is, lower than the natural four dimensions) but none of these gives the right framework for comparing doubly degenerate tricategories with braided monoidal categories, for example.  

Our approach is to consider bicategories as categories ``weakly enriched'' in \cat{Cat}, and then generalise this to replace \Cat\ with any symmetric monoidal bicategory \cB.    We then make a bicategory with
\begin{itemize}
 \item 0-cells: weak $\cB$-categories,
\item 1-cells: weak $\cB$-functors, and
\item 2-cells: $\cB$-icons.
\end{itemize}
We name this 2-category $\cat{$\cB$-Icon}$ after its 2-cells.  Our key theorem \ref{smb} shows that if $\cB$ is symmetric monoidal then so is $\cat{$\cB$-Icon}$, so that we can iterate the construction.

We then have the following inductive process: starting with $\cB_0 = \Cat$ and putting $\cB_{k+1} = \cat{$\cB_k$-Icon}$, we finally take the full sub-2-category on suitably degenerate objects.  This gives the correct totalities of monoidal categories as shown in the following schematic diagram.

\[\begin{array}{lclcl}
\cB_0 &=& \cat{Cat}  \\
\cB_1 &=& \cat{$\cB_0$-Icon} & \ltmap{\mbox{\small\sf{degen}}} & \cat{MonCat} \\[4pt]
\cB_2 &=& \cat{$\cB_1$-Icon} & \ltmap{\mbox{\small\sf{2-degen}}} & \cat{BrMonCat} \\[4pt]
\cB_3 &=& \cat{$\cB_2$-Icon} & \ltmap{\mbox{\small\sf{3-degen}}} & \cat{SymMonCat} \\[4pt]
\cB_4 &=& \cat{$\cB_3$-Icon} & \ltmap{\mbox{\small\sf{4-degen}}} & \cat{SymMonCat} \\[4pt]
&&& \vdots  
  \end{array}\]

Note that the basic definition of $\cat{\cB-Icon}$ only requires that \cB\ be monoidal.  However, in order to iterate the icon construction
\[\cl{B} \tmapsto \cat{\cB-Icon},\] 
we need $\cat{\cB-Icon}$ to inherit the relevant structure of $\cB$.  We have the following structures.
\[\begin{array}{ccc}
   \cB &\tmapsto& \cat{\cB-Icon} \\[4pt]
\mbox{monoidal} & & \mbox{---} \\
\mbox{braided monoidal} && \mbox{monoidal} \\
\mbox{symmetric monoidal} && \mbox{symmetric monoidal}
  \end{array}\]
Thus we see that for the iteration to proceed, we need symmetric monoidal bicategories.  This is the main technical content of this work.  The work was first presented \cite{cg4} but suffered some delay due to the technical complications in the construction of the symmetric monoidal bicategory \Bicon.  These were resolved in \cite{go1}.  In the meantime, Shulman \cite{shul2} gave an elegant abstract account of \cB-icons via 2-monads, but does not focus on the totality.  His focus is on the 0-cells, and he does not study iteration, so he does not need a \emph{symmetric} monoidal structure on \cB, or on \cat{\cB-Icon}.

%
%
%
%
%

This paper is structured as follows.  In Section~\ref{basic} we recall the basic theory of icons as in \cite{lac2}.  In Section~\ref{monbicatbackground} we recall the definitions of various kinds of monoidal bicategory (braided, sylleptic, symmetric), and briefly remind the reader of the relevant coherence theorems that we will use repeatedly in later proofs.  Then in Section~\ref{genicons} we generalise the definitions of Section~\ref{basic} to the enriched setting and prove that if \cB\ is a symmetric monoidal bicategory then so is \cat{$\cl{B}$-Icon} so we can iterate the icon construction. We study this iteration in Section~\ref{hd}.  We show that this gives the correct totalities of degenerate structures to correspond with the second column of the Periodic Table, and make some conjectures about how to generalise this for the higher-dimensional columns.  

\section{Icons}\label{basic}

In this section we recall Lack's definition of icon \cite{lac2, lp1}.  The idea is to form a convenient 2-category of bicategories.  Bicategories together with functors, transformations and modifications form a tricategory, and while truncating this to the level of functors does yield a category, truncating to the level of transformations does not yield a coherent 2-dimensional structure.  This is because transformations between bicategories do not in general compose in a strictly associative or unital way, nor is interchange strict.

Lack's idea is to introduce a new type of transformation called \emph{icon} which can be thought of as `Identity Component Oplax Natural transformation''.  However, composition is not merely that of transformations---composing identity 1-cell components $I_a$ would give a component $I_a \circ I_a$, not itself necessarily an identity.

Instead, Lack's insight is to drop all reference to the 1-cell components of these ``transformations'' and only retain the 2-cell components as data.  This results in the following definition.

\begin{mydefinition} 

Let $X$ and $Y$ be bicategories and let $F,G \: X \tra Y$ be weak functors such that $Fa=Ga$ for all objects $a \in X$.  An \demph{icon}
\[
\psset{unit=0.1cm,labelsep=2pt,nodesep=2pt}
\pspicture(20,20)

\rput(0,10){\rnode{a1}{$X$}}  
\rput(20,10){\rnode{a2}{$Y$}}  

\ncarc[arcangle=45]{->}{a1}{a2}\naput{{\scriptsize $F$}}
\ncarc[arcangle=-45]{->}{a1}{a2}\nbput{{\scriptsize $G$}}

\pcline[linewidth=0.6pt,doubleline=true,arrowinset=0.6,arrowlength=0.8,arrowsize=0.5pt 2.1]{->}(10,13)(10,7)  \naput{{\scriptsize $\alpha$}}


\endpspicture
\]
is given by, for all pairs of 0-cells $a,b \in X$ a 2-cell

\[
\psset{unit=0.1cm,labelsep=2pt,nodesep=2pt}
\pspicture(0,-4)(30,24)

\rput(0,10){\rnode{a1}{$X(a,b)$}}  
\rput(30,13){\rnode{a2}{$Y(Fa,Fb)$}}  
\rput(31,10){\rotatebox{90}{$=$}}
\rput(30,7){\rnode{a3}{$Y(Ga,Gb)$}}  

\ncarc[arcangle=45]{->}{a1}{a2}\naput[npos=0.6]{{\scriptsize $F$}}
\ncarc[arcangle=-45]{->}{a1}{a3}\nbput[npos=0.6]{{\scriptsize $G$}}

{
\rput[c](17,8.5){\psset{unit=1mm,doubleline=true,arrowinset=0.6,arrowlength=0.5,arrowsize=0.5pt 2.1,nodesep=0pt,labelsep=2pt}
\pcline{->}(0,3)(0,0) \nbput{{\scriptsize $\alpha_{ab}$}}}}


\endpspicture
\]
satisfying the following axioms.  Note that to save space we will omit the name of the bicategory when writing hom-categories, thus $X(a,b)$ becomes $(a,b)$, and $Y(Fa,Fb)$ becomes $(Fa,Fb)$, and so on.

\begin{itemize}
 \item Composition:

\[
\psset{unit=0.1cm,labelsep=2pt,nodesep=2pt}
\pspicture(0,-15)(80,40)

%

\rput(0,30){\rnode{a1}{$(b,c)(a,b)$}}  
\rput(30,33){\rnode{a2}{\makebox[4em][l]{$(Fb,Fc)(Fa,Fb)$}}}  
\rput(38,30){\rotatebox{90}{$=$}}
\rput(30,27){\rnode{a3}{\makebox[4em][l]{$(Gb,Gc)(Ga,Gb)$}}}  

\rput(0,3){\rnode{b1}{$(a,c)$}}  
\rput(30,3){\rnode{b2}{\makebox[3em][l]{$(Ga,Gc)$}}}

\ncarc[arcangle=45]{->}{a1}{a2}\naput[npos=0.56]{{\scriptsize $FF$}}
\ncarc[arcangle=-45]{->}{a1}{a3}\nbput[npos=0.56]{{\scriptsize $GG$}}

\ncline{->}{a1}{b1} \nbput{{\scriptsize $m$}}
\ncline{->}{a3}{b2} \naput{{\scriptsize $m'$}}

\ncarc[arcangle=-45]{->}{b1}{b2}\nbput[npos=0.5]{{\scriptsize $G$}}

{
\rput[c](17,28.5){\psset{unit=1mm,doubleline=true,arrowinset=0.6,arrowlength=0.5,arrowsize=0.5pt 2.1,nodesep=0pt,labelsep=2pt}
\pcline{->}(0,3)(0,0) \nbput{{\scriptsize $\alpha\alpha$}}}}

{
\rput[c](17,8.5){\psset{unit=1mm,doubleline=true,arrowinset=0.6,arrowlength=0.5,arrowsize=0.5pt 2.1,nodesep=0pt,labelsep=2pt}
\pcline{->}(0,3)(0,0) \nbput{{\scriptsize $\phi^G$}}}}


\rput(60,20){$=$}

\endpspicture
\psset{unit=0.1cm,labelsep=2pt,nodesep=2pt}
\pspicture(0,-15)(50,40)

%
%

\rput(0,30){\rnode{a1}{$(b,c)(a,b)$}}  
\rput(30,30){\rnode{a2}{\makebox[4em][l]{$(Fb,Fc)(Fa,Fb)$}}}  

\rput(0,3){\rnode{b1}{$(a,c)$}}  
\rput(30,6){\rnode{b2}{\makebox[3em][l]{$(Fa,Fc)$}}}  

\rput(32,3){\rotatebox{90}{$=$}}
\rput(30,0){\rnode{b3}{\makebox[3em][l]{$(Ga,Gc)$}}}  

\ncarc[arcangle=45]{->}{a1}{a2}\naput[npos=0.56]{{\scriptsize $FF$}}

\ncline{->}{a1}{b1} \nbput{{\scriptsize $m$}}
\ncline{->}{a2}{b2} \naput{{\scriptsize $m'$}}

\ncarc[arcangle=45]{->}{b1}{b2}\naput[npos=0.56]{{\scriptsize $F$}}
\ncarc[arcangle=-45]{->}{b1}{b3}\nbput[npos=0.56]{{\scriptsize $G$}}

{
\rput[c](16,21){\psset{unit=1mm,doubleline=true,arrowinset=0.6,arrowlength=0.5,arrowsize=0.5pt 2.1,nodesep=0pt,labelsep=2pt}
\pcline{->}(0,3)(0,0) \nbput{{\scriptsize $\phi^F$}}}}

{
\rput[c](17,0.5){\psset{unit=1mm,doubleline=true,arrowinset=0.6,arrowlength=0.5,arrowsize=0.5pt 2.1,nodesep=0pt,labelsep=2pt}
\pcline{->}(0,3)(0,0) \nbput{{\scriptsize $\alpha$}}}}


\endpspicture
\]

\item Unit:

\[
\psset{unit=0.15cm,labelsep=2pt,nodesep=2pt}
\pspicture(45,20)


\rput(2,10){\rnode{a1}{$\1$}}
\rput(25,18){\rnode{a2}{$(a,a)$}}
\rput(25,2){\rnode{a3}{$(Fa,Fa)=(Ga,Ga)$}}

\ncline{->}{a1}{a2} \naput{{\scriptsize $I$}}
\ncline{->}{a1}{a3} \nbput{{\scriptsize $I'$}}

\ncarc[arcangle=45]{->}{a2}{a3}\naput[npos=0.5]{{\scriptsize $G$}}
\ncarc[arcangle=-45]{->}{a2}{a3}\nbput[npos=0.5]{{\scriptsize $F$}}

{
\rput[c](14,9){\psset{unit=1mm,doubleline=true,arrowinset=0.6,arrowlength=0.5,arrowsize=0.5pt 2.1,nodesep=0pt,labelsep=-2pt}
\pcline{->}(0,0)(3,3) \naput{{\scriptsize $\phi^F$}}}}

{
\rput[c](24,10){\psset{unit=1mm,doubleline=true,arrowinset=0.6,arrowlength=0.5,arrowsize=0.5pt 2.1,nodesep=0pt,labelsep=2pt}
\pcline{->}(0,0)(4,0) \nbput{{\scriptsize $\alpha$}}}}

\rput(40,10){$=$}

\endpspicture
\psset{unit=0.15cm,labelsep=2pt,nodesep=2pt}
\pspicture(30,20)


\rput(2,10){\rnode{a1}{$\1$}}
\rput(25,18){\rnode{a2}{$(a,a)$}}
\rput(25,2){\rnode{a3}{$(Fa,Fa)=(Ga,Ga)$}}

\ncline{->}{a1}{a2} \naput{{\scriptsize $I$}}
\ncline{->}{a1}{a3} \nbput{{\scriptsize $I'$}}

\ncarc[arcangle=45]{->}{a2}{a3}\naput[npos=0.5]{{\scriptsize $G$}}

{
\rput[c](17,9){\psset{unit=1mm,doubleline=true,arrowinset=0.6,arrowlength=0.5,arrowsize=0.5pt 2.1,nodesep=0pt,labelsep=-2pt}
\pcline{->}(0,0)(3,3) \naput{{\scriptsize $\phi^G$}}}}


\endpspicture\]

\end{itemize}

\end{mydefinition}

Lack notes the following result.

\begin{theorem}
Bicategories, weak functors and icons form a 2-category.
\end{theorem}

\section{Background on monoidal bicategories}\label{monbicatbackground}

In this section we provide the full definition of symmetric monoidal bicategory, which is necessary for the constructions in Section~\ref{genicons}.  We build this up gradually.  We first recall the definition of monoidal bicategory, which can also be found in \cite{cg3} and \cite{gur3}.  We then define the notion of a braiding, as in \cite{gur4}, followed by syllepsis and symmetry as in \cite{go1} and \cite{sch1}.  While these definitions are not new, we feel it is useful to gather them all together in one place, which does not appear to have been done before.  

\begin{definition} \label{tricatdef} A monoidal bicategory $M$ consists of the following data subject to the following axioms: \\
DATA:
\begin{itemize}
\item an underlying bicategory, also denoted $M$;
\item a functor $\otimes:M \times M \rightarrow M$;
\item a functor $I: 1 \rightarrow M$, where 1 denotes the unit bicategory;
\item an adjoint equivalence $\bs{a}$
\[
\xymatrix{
M \times M \times M \ar[rr]^{\otimes \times 1} \ar[d]_{ 1 \times \otimes} & & M \times M \ar[d]^{\otimes} \ar@{}[dll]|{\Downarrow \bs{a}}\\
M \times M \ar[rr]_{\otimes} & & M }
\]
in $\mathbf{Bicat}(M^{3}, M)$;
\item adjoint equivalences $\bs{l}$ and $\bs{r}$
\[
\xy
{\ar^-{I \times 1} (0,0)*+{M}; (15, 15)*+{M \times M} };
{\ar_{1} (0,0)*+{M}; (30,0)*+{M} };
{\ar^{\otimes} (15, 15)*+{M \times M}; (30,0)*+{M} };
(14.5,5)*{\Downarrow \bs{l}};
{\ar_{1} (50,0)*+{M}; (80,0)*+{M} };
{\ar^{1 \times I} (50,0)*+{M}; (65,15)*+{M \times M} };
{\ar^{\otimes} (65,15)*+{M \times M}; (80,0)*+{M} };
(64.5,5)*{\Downarrow \bs{r}}
\endxy
\]
in $\mathbf{Bicat}(M, M)$;
\item an invertible 2-cell $\pi$ (i.e., an invertible modification)
\[
\def\objectstyle{\scriptstyle}
\def\labelstyle{\scriptstyle}
\xy0;/r.20pc/:
{\ar_{1 \times \otimes} (-65,0)*+{M^{3}}; (-50, -16)*+{M^{2}} };
{\ar_{\otimes} (-50,-16)*+{M^{2}}; (-22, -16)*+{M} };
{\ar_{1 \times 1 \times \otimes} (-50, 16)*+{M^{4}}; (-65,0)*+{M^{3}} };
{\ar|{1 \times \otimes \times 1} (-50, 16)*+{M^{4}}; (-35,0)*+{M^{3}} };
{\ar^{1 \times \otimes} (-35,0)*+{M^{3}}; (-50, -16)*+{M^{2}} };
{\ar^{\otimes \times 1 \times 1} (-50, 16)*+{M^{4}}; (-22,16)*+{M^{3}} };
{\ar^{\otimes \times 1} (-35,0)*+{M^{3}}; (-7,0)*+{M^{2}} };
{\ar^{\otimes \times 1}(-22,16)*+{M^{3}}; (-7,0)*+{M^{2}} };
{\ar^{\otimes} (-7,0)*+{M^{2}};  (-22, -16)*+{M} };
{\ar@3{->}^{\pi} (-3.5,0)*{}; (3.5,0)*{} };
{\ar_{1 \times \otimes} (7,0)*+{M^{3}}; (22,-16)*+{M^{2}} };
{\ar_{\otimes} (22,-16)*+{M^{2}}; (50,-16)*+{M} };
{\ar^{\otimes \times 1} (7,0)*+{M^{3}}; (35,0)*+{M^{2}} };
{\ar_{1 \times 1 \times \otimes} (22,16)*+{M^{4}}; (7,0)*+{M^{3}} };
{\ar^{\otimes \times 1 \times 1} (22,16)*+{M^{4}}; (50,16)*+{M^{3}} };
{\ar|{1 \times \otimes} (50,16)*+{M^{3}}; (35,0)*+{M^{2}} };
{\ar_{\otimes} (35,0)*+{M^{2}}; (50,-16)*+{M} };
{\ar^{\otimes \times 1} (50,16)*+{M^{3}}; (65,0)*+{M^{2}} };
{\ar^{\otimes} (65,0)*+{M^{2}}; (50,-16)*+{M} };
(49,-1)*{\scriptstyle \Leftarrow}; (49,2)*{\scriptstyle a};
(24,-7)*{\scriptstyle \Downarrow a}; (24,8)*{\scriptstyle =};
(-24,8)*{\scriptstyle \Downarrow a \times 1};
(-24,-7)*{\scriptstyle \Downarrow a};
(-49,2)*{\scriptstyle 1 \times a}; (-49,-1)*{\scriptstyle \Leftarrow}
\endxy
\]
in the bicategory $\mathbf{Bicat}(M^{4}, M)$;
\item invertible modifications
\[
\xy0;/r.20pc/:
{\ar@/^1.5pc/^{1} (0,0)*+{M^{2}}; (50,-10)*+{M^{2}} };
{\ar|{1 \times I \times 1} (0,0)*+{M^{2}}; (20,-10)*+{M^{3}} };
{\ar@/_1.5pc/_{1} (0,0)*+{M^{2}}; (20,-30)*+{M^{2}} };
{\ar_{\otimes \times 1} (20,-10)*+{M^{3}}; (50,-10)*+{M^{2}} };
{\ar^{1 \times \otimes} (20,-10)*+{M^{3}}; (20,-30)*+{M^{2}} };
{\ar_{\otimes} (20,-30)*+{M^{2}}; (50,-30)*+{M} };
{\ar^{\otimes} (50,-10)*+{M^{2}}; (50,-30)*+{M} };
{\ar@/^1.5pc/^{1} (65,0)*+{M^{2}}; (115,-10)*+{M^{2}} };
{\ar@/_1.5pc/_{1} (65,0)*+{M^{2}}; (85,-30)*+{M^{2}} };
{\ar_{\otimes} (85,-30)*+{M^{2}}; (115,-30)*+{M} };
{\ar^{\otimes} (115,-10)*+{M^{2}}; (115,-30)*+{M} };
(25,-3)*{\scriptstyle \Downarrow r^{\dot} \times 1}; (12,-17)*{\scriptstyle \Leftarrow 1 \times l}; (36,-20)*{\scriptstyle \Downarrow a};
{\ar@{=>}_{1} (100,-10)*{}; (80,-20)*{} };
{\ar@3{->}^{\mu} (55,-15)*{}; (65,-15)*{} }
\endxy
\]
\[
\def\objectstyle{\scriptstyle}
\def\labelstyle{\scriptstyle}
\xy0;/r.20pc/:
{\ar^{I \times 1 \times 1} (0,0)*+{M^{2}}; (20,8)*+{M^{3}} };
{\ar^{\otimes \times 1} (20,8)*+{M^{3}}; (40,0)*+{M^{2}} };
{\ar_{\otimes} (0,0)*+{M^{2}}; (0,-16)*+{M} };
{\ar_{1} (0,-16)*+{M}; (40,-16)*+{M} };
{\ar^{\otimes} (40,0)*+{M^{2}}; (40,-16)*+{M} };
{\ar_{1} (0,0)*+{M^{2}}; (40,0)*+{M^{2}} };
{\ar_{\otimes} (60,0)*+{M^{2}}; (60,-16)*+{M} };
{\ar_{1} (60,-16)*+{M}; (100,-16)*+{M} };
{\ar^{I \times 1 \times 1} (60,0)*+{M^{2}}; (80,8)*+{M^{3}} };
{\ar^{\otimes \times 1} (80,8)*+{M^{3}}; (100,0)*+{M^{2}} };
{\ar^{\otimes} (100,0)*+{M^{2}}; (100,-16)*+{M} };
{\ar^{I \times 1} (60,-16)*+{M}; (80,-4)*+{M^{2}} };
{\ar^{\otimes} (80,-4)*+{M^{2}}; (100,-16)*+{M} };
{\ar_{1 \times \otimes} (80,8)*+{M^{3}}; (80,-4)*+{M^{2}} };
(20,-7)*{\scriptstyle =}; (18,3.3)*{\scriptstyle \Downarrow l \times 1};
(70,-2)*{\scriptstyle =}; (90,-2)*{\scriptstyle \Downarrow a}; (80,-10)*{\scriptstyle \Downarrow l};
{\ar@3{->}^{\lambda} (45,-3)*{}; (55,-3)*{} }
\endxy
\]
\[
\def\objectstyle{\scriptstyle}
\def\labelstyle{\scriptstyle}
\xy0;/r.20pc/:
{\ar^{1} (0,0)*+{M}; (40,0)*+{M} };
{\ar^{\otimes} (0,-16)*+{M^{2}}; (0,0)*+{M} };
{\ar_{1 \times 1 \times I} (0,-16)*+{M^{2}}; (20,-24)*+{M^{3}} };
{\ar_{1 \times \otimes} (20,-24)*+{M^{3}}; (40,-16)*+{M^{2}} };
{\ar_{\otimes} (40,-16)*+{M^{2}}; (40,0)*+{M} };
{\ar^{1} (0,-16)*+{M^{2}}; (40,-16)*+{M^{2}} };
{\ar^{\otimes} (60,-16)*+{M^{2}}; (60,0)*+{M} };
{\ar^{1} (60,0)*+{M}; (100,0)*+{M} };
{\ar_{1 \times 1 \times I} (60,-16)*+{M^{2}}; (80,-24)*+{M^{3}} };
{\ar_{1 \times \otimes} (80,-24)*+{M^{3}}; (100,-16)*+{M^{2}} };
{\ar_{\otimes} (100,-16)*+{M^{2}}; (100,0)*+{M} };
{\ar_{1 \times I} (60,0)*+{M}; (80,-12)*+{M^{2}} };
{\ar_{\otimes} (80,-12)*+{M^{2}}; (100,0)*+{M} };
{\ar^{\otimes \times 1} (80,-24)*+{M^{3}}; (80,-12)*+{M^{2}} };
(20,-19)*{\scriptstyle \Downarrow 1 \times r^{\dot}}; (20,-7)*{\scriptstyle =};
(80,-6)*{\scriptstyle \Downarrow r^{\dot}}; (70,-14)*{\scriptstyle =}; (90,-14)*{\scriptstyle \Rightarrow a};
{\ar@3{->}^{\rho} (45,-10)*{}; (55,-10)*{} }
\endxy
\]
\end{itemize}
AXIOMS:
\begin{itemize}
\item The following equation of 2-cells holds for all objects $a,b,c,d,e$ in $M$, where we have used parentheses instead of $\otimes$ for compactness and the unmarked isomorphisms are naturality isomorphisms for $a$.
\[
\def\objectstyle{\scriptstyle}
\def\labelstyle{\scriptstyle}
\xy0;/r.18pc/:
{\ar^{(a1)1} (0,0)*+{(((ab)c)d)e}; (0,20)*+{((a(bc))d)e} };
{\ar^{a1} (0,20)*+{((a(bc))d)e}; (20,30)*+{(a((bc)d))e} };
{\ar^{(1a)1} (20,30)*+{(a((bc)d))e}; (50,35)*+{(a(b(cd)))e} };
{\ar^{a} (50,35)*+{(a(b(cd)))e}; (80,30)*+{a((b(cd))e)} };
{\ar^{1a} (80,30)*+{a((b(cd))e)}; (100,20)*+{a(b((cd)e))} };
{\ar^{1(1a)} (100,20)*+{a(b((cd)e))}; (100,0)*+{a(b(c(de)))} };
{\ar_{a} (0,0)*+{(((ab)c)d)e}; (30,-15)*+{((ab)c)(de)} };
{\ar_{a} (30,-15)*+{((ab)c)(de)}; (70,-15)*+{(ab)(c(de))} };
{\ar_{a} (70,-15)*+{(ab)(c(de))}; (100,0)*+{a(b(c(de)))} };
{\ar^{a} (20,30)*+{(a((bc)d))e}; (50,20)*+{a(((bc)d)e)} };
{\ar^{1(a1)} (50,20)*+{a(((bc)d)e)}; (80,30)*+{a((b(cd))e)} };
{\ar_{1a} (50,20)*+{a(((bc)d)e)}; (70,5)*+{a((bc)(de))} };
{\ar^{1a} (70,5)*+{a((bc)(de))}; (100,0)*+{a(b(c(de)))} };
{\ar_{a} (0,20)*+{((a(bc))d)e}; (30,5)*+{(a(bc))(de)} };
{\ar_{a} (30,5)*+{(a(bc))(de)}; (70,5)*+{a((bc)(de))} };
{\ar_{a(11)} (30,-15)*+{((ab)c)(de)}; (30,5)*+{(a(bc))(de)} };
{\ar^{(a1)1} (0,-70)*+{(((ab)c)d)e}; (0,-50)*+{((a(bc))d)e} };
{\ar^{a1} (0,-50)*+{((a(bc))d)e}; (20,-40)*+{(a((bc)d))e} };
{\ar^{(1a)1} (20,-40)*+{(a((bc)d))e}; (50,-35)*+{(a(b(cd)))e} };
{\ar^{a} (50,-35)*+{(a(b(cd)))e}; (80,-40)*+{a((b(cd))e)} };
{\ar^{1a} (80,-40)*+{a((b(cd))e)}; (100,-50)*+{a(b((cd)e))} };
{\ar^{1(1a)} (100,-50)*+{a(b((cd)e))}; (100,-70)*+{a(b(c(de)))} };
{\ar_{a} (0,-70)*+{(((ab)c)d)e}; (30,-85)*+{((ab)c)(de)} };
{\ar_{a} (30,-85)*+{((ab)c)(de)}; (70,-85)*+{(ab)(c(de))} };
{\ar_{a} (70,-85)*+{(ab)(c(de))}; (100,-70)*+{a(b(c(de)))} };
{\ar_{a1} (0,-70)*+{(((ab)c)d)e}; (35,-55)*+{((ab)(cd))e} };
{\ar_{a} (35,-55)*+{((ab)(cd))e}; (65,-65)*+{(ab)((cd)e)} };
{\ar_{(11)a} (65,-65)*+{(ab)((cd)e)}; (70,-85)*+{(ab)(c(de))} };
{\ar^{a} (65,-65)*+{(ab)((cd)e)}; (100,-50)*+{a(b((cd)e))} };
{\ar_{a1} (35,-55)*+{((ab)(cd))e}; (50,-35)*+{(a(b(cd)))e} };
{\ar@{=} (50,-20)*{}; (50,-30)*{} };
(50,27)*{\cong}; (13,4)*{\cong};
(30,16)*{\Downarrow \pi}; (78,16)*{\Downarrow 1 \pi}; (60,-3)*{\Downarrow \pi};
(82,-68)*{\cong}; (20,-50)*{\Downarrow \pi 1}; (70,-50)*{\Downarrow \pi}; (35,-68)*{\Downarrow \pi}
\endxy
\]
\item The following equation of 2-cells holds for all objects $a,b,c$ in $M$, where the unmarked isomorphisms are either naturality isomorphisms for $a$ or unique coherence isomorphisms from the hom-bicategory.
\[
\def\objectstyle{\scriptstyle}
\def\labelstyle{\scriptstyle}
\xy0;/r.16pc/:
{\ar^{(r^{\dot}1)1} (0,0)*+{(ab)c}; (30,15)*+{((aI)b)c} };
{\ar^{a1} (30,15)*+{((aI)b)c}; (60,30)*+{(a(Ib))c} };
{\ar^{(1l)1} (60,30)*+{(a(Ib))c}; (90,15)*+{(ab)c} };
{\ar^{a} (90,15)*+{(ab)c}; (90,-15)*+{a(bc)} };
{\ar_{a} (0,0)*+{(ab)c}; (30,-15)*+{a(bc)} };
{\ar@/_1.5pc/_{1} (30,-15)*+{a(bc)}; (90,-15)*+{a(bc)} };
{\ar_{a} (30,15)*+{((aI)b)c}; (30,0)*+{(aI)(bc)} };
{\ar_{r^{\dot}(11)} (30,-15)*+{a(bc)}; (30,0)*+{(aI)(bc)} };
{\ar_{a} (60,30)*+{(a(Ib))c}; (60,10)*+{a((Ib)c)} };
{\ar_{1a} (60,10)*+{a((Ib)c)}; (60,-10)*+{a(I(bc))} };
{\ar^{1(l1)} (60,10)*+{a((Ib)c)}; (90,-15)*+{a(bc)} };
{\ar_{1l} (60,-10)*+{a(I(bc))}; (90,-15)*+{a(bc)} };
{\ar^{a} (30,0)*+{(aI)(bc)}; (60,-10)*+{a(I(bc))} };
{\ar^{(r^{\dot}1)1} (0,-70)*+{(ab)c}; (30,-55)*+{((aI)b)c} };
{\ar^{a1} (30,-55)*+{((aI)b)c}; (60,-40)*+{(a(Ib))c} };
{\ar^{(1l)1} (60,-40)*+{(a(Ib))c}; (90,-55)*+{(ab)c} };
{\ar^{a} (90,-55)*+{(ab)c}; (90,-75)*+{a(bc)} };
{\ar_{a} (0,-70)*+{(ab)c}; (30,-75)*+{a(bc)} };
{\ar_{1} (30,-75)*+{a(bc)}; (90,-75)*+{a(bc)} };
{\ar_{11} (0,-70)*+{(ab)c}; (90,-55)*+{(ab)c} };
{\ar@{=} (45,-25)*{}; (45,-40)*{} };
(15,0)*{\cong}; (45,10)*{\Downarrow \pi}; (75,10)*{\cong}; (70,-5)*{\Downarrow 1 \lambda}; (50,-15)*{\Downarrow \mu};
(50,-53)*{\Downarrow \mu 1}; (55,-65)*{\cong}
\endxy
\]
\item The following equation of 2-cells holds for all objects $a,b,c$ in $M$.
\[
\def\objectstyle{\scriptstyle}
\def\labelstyle{\scriptstyle}
\xy0;/r.16pc/:
{\ar^{a} (0,-15)*+{(ab)c}; (0,15)*+{a(bc)} };
{\ar^{1(r^{\dot}1)} (0,15)*+{a(bc)}; (30,30)*+{a((bI)c)} };
{\ar^{1a} (30,30)*+{a((bI)c)}; (60,15)*+{a(b(Ic))} };
{\ar^{1(1l)} (60,15)*+{a(b(Ic))}; (90,0)*+{a(bc)} };
{\ar@/_1.5pc/_{1} (0,-15)*+{(ab)c}; (60,-15)*+{(ab)c} };
{\ar_{a} (60,-15)*+{(ab)c}; (90,0)*+{a(bc)} };
{\ar^{(1r^{\dot})1} (0,-15)*+{(ab)c}; (30,10)*+{(a(bI))c} };
{\ar_{r^{\dot} 1} (0,-15)*+{(ab)c}; (30,-10)*+{((ab)I)c} };
{\ar_{a1} (30,-10)*+{((ab)I)c}; (30,10)*+{(a(bI))c} };
{\ar^{a} (30,-10)*+{((ab)I)c}; (60,0)*+{(ab)(Ic)} };
{\ar_{(11)l} (60,0)*+{(ab)(Ic)}; (60,-15)*+{(ab)c} };
{\ar^{a} (60,0)*+{(ab)(Ic)}; (60,15)*+{a(b(Ic))} };
{\ar_{a} (30,10)*+{(a(bI))c}; (30,30)*+{a((bI)c)} };
{\ar^{a} (0,-75)*+{(ab)c}; (0,-55)*+{a(bc)} };
{\ar^{1(r^{\dot}1)} (0,-55)*+{a(bc)}; (30,-40)*+{a((bI)c)} };
{\ar^{1a} (30,-40)*+{a((bI)c)}; (60,-55)*+{a(b(Ic))} };
{\ar^{1(1l)} (60,-55)*+{a(b(Ic))}; (90,-70)*+{a(bc)} };
{\ar_{1} (0,-75)*+{(ab)c}; (60,-75)*+{(ab)c} };
{\ar_{a} (60,-75)*+{(ab)c}; (90,-70)*+{a(bc)} };
{\ar_{11} (0,-55)*+{a(bc)}; (90,-70)*+{a(bc)} };
{\ar@{=} (45,-25)*{}; (45,-40)*{} };
(15,12)*{\cong}; (19,-6)*{\Downarrow \rho 1}; (45,10)*{\Downarrow \pi}; (75,0)*{\cong}; (40,-15)*{\Downarrow \mu};
(40,-55)*{\Downarrow 1 \mu}; (40,-70)*{\cong}
\endxy
\]
\end{itemize}
\end{definition}

\begin{definition} Let $B = (B, \otimes, I, \mathbf{a}, \mathbf{l}, \mathbf{r}, \pi, \mu, \rho, \lambda)$ be a monoidal bicategory.  Then a \textit{braiding} for $B$ consists of
\begin{itemize}
\item an adjoint equivalence $\bs{R}: \otimes \Rightarrow \otimes \circ \tau$ in $\mathbf{Bicat}(B \times B, B)$, where we define $\tau:B \times B \rightarrow B \times B$ to interchange the coordinates;
\item an invertible modification $R_{(-|-,-)}$ as displayed below;
\[
\xy
{\ar^{R1} (0,0)*+{(AB)C}; (15,10)*+{(BA)C} };
{\ar^{a} (15,10)*+{(BA)C}; (45,10)*+{B(AC)} };
{\ar^{1R} (45,10)*+{B(AC)}; (60,0)*+{B(CA)} };
{\ar_{a} (0,0)*+{(AB)C}; (15,-10)*+{A(BC)} };
{\ar_{R} (15,-10)*+{A(BC)}; (45,-10)*+{(BC)A} };
{\ar_{a} (45,-10)*+{(BC)A}; (60,0)*+{B(CA)} };
{\ar@{=>}^{R_{(A|B,C)}} (29, 7)*{}; (29,-7)*{} }
\endxy
\]
\item and an invertible modification $R_{(-,-|-)}$ as displayed below;
\[
\xy
{\ar^{1R} (0,0)*+{A(BC)}; (15,10)*+{A(CB)} };
{\ar^{a^{\centerdot}} (15,10)*+{A(CB)}; (45,10)*+{(AC)B} };
{\ar^{R1} (45,10)*+{(AC)B}; (60,0)*+{(CA)B} };
{\ar_{a^{\centerdot}} (0,0)*+{A(BC)}; (15,-10)*+{(AB)C} };
{\ar_{R} (15,-10)*+{(AB)C}; (45,-10)*+{C(AB)} };
{\ar_{a^{\centerdot}} (45,-10)*+{C(AB)}; (60,0)*+{(CA)B} };
{\ar@{=>}^{R_{(A,B|C)}} (29, 7)*{}; (29,-7)*{} }
\endxy
\]
\end{itemize}
all subject to the following four axioms.
\[
\xy
{\ar^{a} (0,0)*+{\scriptstyle (AB)(CD)}; (8,8)*+{\scriptstyle A(B(CD))} };
{\ar^{1(1R)} (8,8)*+{\scriptstyle A(B(CD))}; (16,16)*+{\scriptstyle A(B(DC))} };
{\ar^{1a^{\centerdot}} (16,16)*+{\scriptstyle A(B(DC))}; (24,24)*+{\scriptstyle A((BD)C)} };
{\ar^{1(R1)} (24,24)*+{\scriptstyle A((BD)C)}; (50,24)*+{\scriptstyle A((DB)C)} };
{\ar^{1a} (50,24)*+{\scriptstyle A((DB)C)}; (76,24)*+{\scriptstyle A(D(BC))} };
{\ar^{a^{\centerdot}} (76,24)*+{\scriptstyle A(D(BC))}; (84,16)*+{\scriptstyle (AD)(BC)} };
{\ar^{a^{\centerdot}} (84,16)*+{\scriptstyle (AD)(BC)}; (92,8)*+{\scriptstyle ((AD)B)C} };
{\ar^{(R1)1} (92,8)*+{\scriptstyle ((AD)B)C}; (100,0)*+{\scriptstyle ((DA)B)C} };
{\ar_{a^{\centerdot}} (0,0)*+{\scriptstyle (AB)(CD)}; (20,-14)*+{\scriptstyle ((AB)C)D} };
{\ar_{R} (20,-14)*+{\scriptstyle ((AB)C)D}; (50,-14)*+{\scriptstyle D((AB)C)} };
{\ar_{a^{\centerdot}} (50,-14)*+{\scriptstyle D((AB)C)}; (80,-14)*+{\scriptstyle (D(AB))C} };
{\ar_{a^{\centerdot}1} (80,-14)*+{\scriptstyle (D(AB))C}; (100,0)*+{\scriptstyle ((DA)B)C} };
{\ar^{1a^{\centerdot}} (8,8)*+{\scriptstyle A(B(CD))}; (42,16)*+{\scriptstyle A((BC)D)} };
{\ar^{1R} (42,16)*+{\scriptstyle A((BC)D)}; (76,24)*+{\scriptstyle A(D(BC))} };
{\ar_{a^{\centerdot}} (42,16)*+{\scriptstyle A((BC)D)}; (30,0)*+{\scriptstyle (A(BC))D} };
{\ar_{a^{\centerdot}1} (30,0)*+{\scriptstyle (A(BC))D}; (20,-14)*+{\scriptstyle ((AB)C)D} };
{\ar^{R} (30,0)*+{\scriptstyle (A(BC))D}; (50,0)*+{\scriptstyle D(A(BC))} };
{\ar^{a^{\centerdot}} (50,0)*+{\scriptstyle D(A(BC))}; (74,0)*+{\scriptstyle (DA)(BC)} };
{\ar_{a^{\centerdot}} (74,0)*+{\scriptstyle (DA)(BC)}; (100,0)*+{\scriptstyle ((DA)B)C} };
{\ar^{1a^{\centerdot}} (50,0)*+{\scriptstyle D(A(BC))}; (50,-14)*+{\scriptstyle D((AB)C)} };
{\ar_{R1} (84,16)*+{\scriptstyle (AD)(BC)}; (74,0)*+{\scriptstyle (DA)(BC)} };
(16,0)*{\Downarrow \pi}; (34,20)*{\Downarrow 1R_{(B,C|D)}}; (62,11)*{\Downarrow R_{(A,BC|D)}}; (36,-7)*{\cong}; (87,4)*{\cong}; (70,-6)*{\Downarrow \pi};
{\ar@{=} (50,-16)*{}; (50,-25)*{} };
{\ar^{a} (0,-51)*+{\scriptstyle (AB)(CD)}; (8,-43)*+{\scriptstyle A(B(CD))} };
{\ar^{1(1R)} (8,-43)*+{\scriptstyle A(B(CD))}; (16,-35)*+{\scriptstyle A(B(DC))} };
{\ar^{1a^{\centerdot}} (16,-35)*+{\scriptstyle A(B(DC))}; (24,-27)*+{\scriptstyle A((BD)C)} };
{\ar^{1(R1)} (24,-27)*+{\scriptstyle A((BD)C)}; (50,-27)*+{\scriptstyle A((DB)C)} };
{\ar^{1a} (50,-27)*+{\scriptstyle A((DB)C)}; (76,-27)*+{\scriptstyle A(D(BC))} };
{\ar^{a^{\centerdot}} (76,-27)*+{\scriptstyle A(D(BC))}; (84,-35)*+{\scriptstyle (AD)(BC)} };
{\ar^{a^{\centerdot}} (84,-35)*+{\scriptstyle (AD)(BC)}; (92,-43)*+{\scriptstyle ((AD)B)C} };
{\ar^{(R1)1} (92,-43)*+{\scriptstyle ((AD)B)C}; (100,-51)*+{\scriptstyle ((DA)B)C} };
{\ar_{a^{\centerdot}} (0,-51)*+{\scriptstyle (AB)(CD)}; (20,-67.5)*+{\scriptstyle ((AB)C)D} };
{\ar_{R} (20,-67.5)*+{\scriptstyle ((AB)C)D}; (50,-67.5)*+{\scriptstyle D((AB)C)} };
{\ar_{a^{\centerdot}} (50,-67.5)*+{\scriptstyle D((AB)C)}; (80,-67.5)*+{\scriptstyle (D(AB))C} };
{\ar_{a^{\centerdot}1} (80,-67.5)*+{\scriptstyle (D(AB))C}; (100,-51)*+{\scriptstyle ((DA)B)C} };
{\ar_{1R} (0,-51)*+{\scriptstyle (AB)(CD)}; (24,-51)*+{\scriptstyle (AB)(DC)} };
{\ar^{a} (24,-51)*+{\scriptstyle (AB)(DC)}; (16,-35)*+{\scriptstyle A(B(DC))} };
{\ar_{a^{\centerdot}} (24,-51)*+{\scriptstyle (AB)(DC)}; (48,-51)*+{\scriptstyle ((AB)D)C} };
{\ar_{a1} (48,-51)*+{\scriptstyle ((AB)D)C}; (48,-44)*+{\scriptstyle (A(BD))C} };
{\ar_{(1R)1} (48,-44)*+{\scriptstyle (A(BD))C}; (65,-35)*+{\scriptstyle (A(DB))C} };
{\ar^{a^{\centerdot}} (48,-44)*+{\scriptstyle (A(BD))C}; (24,-27)*+{\scriptstyle A((BD)C)} };
{\ar_{a^{\centerdot}} (65,-35)*+{\scriptstyle (A(DB))C}; (50,-27)*+{\scriptstyle A((DB)C)} };
{\ar_{a^{\centerdot}1} (65,-35)*+{\scriptstyle (A(DB))C}; (92,-43)*+{\scriptstyle ((AD)B)C} };
{\ar^{R1} (48,-51)*+{\scriptstyle ((AB)D)C}; (80,-67.5)*+{\scriptstyle (D(AB))C} };
(10,-47)*{\cong}; (31,-43)*{\Downarrow \pi}; (48,-35)*{\cong}; (70,-31)*{\Downarrow \pi}; (72,-50)*{\Downarrow R_{(A,B|D)}1}; (30,-59)*{\Downarrow R_{(AB,C|D)}}\endxy
\]

\[
\xy
{\ar^{R1} (0,0)*+{\scriptstyle (AB)(CD)}; (7,9)*+{\scriptstyle (BA)(CD)} };
{\ar^{a^{\centerdot}} (7,9)*+{\scriptstyle (BA)(CD)}; (14,18)*+{\scriptstyle ((BA)C)D} };
{\ar^{a1} (14,18)*+{\scriptstyle ((BA)C)D}; (38,18)*+{\scriptstyle (B(AC))D} };
{\ar^{(1R)1} (38,18)*+{\scriptstyle (B(AC))D}; (62,18)*+{\scriptstyle (B(CA))D} };
{\ar^{a^{\centerdot}1} (62,18)*+{\scriptstyle (B(CA))D}; (86,18)*+{\scriptstyle ((BC)A)D} };
{\ar^{a} (86,18)*+{\scriptstyle ((BC)A)D}; (93,9)*+{\scriptstyle (BC)(AD)} };
{\ar^{1R} (93,9)*+{\scriptstyle (BC)(AD)}; (100,0)*+{\scriptstyle (BC)(DA)} };
{\ar_{a} (0,0)*+{\scriptstyle (AB)(CD)}; (12,-14)*+{\scriptstyle A(B(CD))} };
{\ar_{R}  (12,-14)*+{\scriptstyle A(B(CD))}; (37,-14)*+{\scriptstyle (B(CD))A} };
{\ar_{a} (37,-14)*+{\scriptstyle (B(CD))A}; (63,-14)*+{\scriptstyle B((CD)A)} };
{\ar_{1a} (63,-14)*+{\scriptstyle B((CD)A)}; (88,-14)*+{\scriptstyle B(C(DA))} };
{\ar_{a^{\centerdot}} (88,-14)*+{\scriptstyle B(C(DA))}; (100,0)*+{\scriptstyle (BC)(DA)} };
{\ar_{a} (7,9)*+{\scriptstyle (BA)(CD)}; (25,0)*+{\scriptstyle B(A(CD))} };
{\ar_{1 a^{\centerdot}} (25,0)*+{\scriptstyle B(A(CD))}; (40,7)*+{\scriptstyle B((AC)D)} };
{\ar^{a^{\centerdot}} (40,7)*+{\scriptstyle B((AC)D)}; (38,18)*+{\scriptstyle (B(AC))D} };
{\ar_{1R} (25,0)*+{\scriptstyle B(A(CD))}; (63,-14)*+{\scriptstyle B((CD)A)} };
{\ar_{1(R1)} (40,7)*+{\scriptstyle B((AC)D)}; (62,6)*+{\scriptstyle B((CA)D)} };
{\ar^{a^{\centerdot}} (62,6)*+{\scriptstyle B((CA)D)}; (62,18)*+{\scriptstyle (B(CA))D} };
{\ar_{1a} (62,6)*+{\scriptstyle B((CA)D)}; (75,-4)*+{\scriptstyle B(C(AD))} };
{\ar_{1(1R)} (75,-4)*+{\scriptstyle B(C(AD))}; (88,-14)*+{\scriptstyle B(C(DA))} };
{\ar^{a^{\centerdot}} (75,-4)*+{\scriptstyle B(C(AD))}; (93,9)*+{\scriptstyle (BC)(AD)} };
(26,10)*{\Downarrow \pi}; (50,11)*{\cong}; (74,10)*{\Downarrow \pi}; (20,-7)*{\Downarrow R_{(A|B,CD)}}; (50,-3)*{\Downarrow 1 R_{(A|C,D)}}; (88,-5)*{\cong};
{\ar@{=} (50,-17)*{}; (50,-25)*{} };
{\ar^{R1} (0,-47)*+{\scriptstyle (AB)(CD)}; (7,-38)*+{\scriptstyle (BA)(CD)} };
{\ar^{a^{\centerdot}} (7,-38)*+{\scriptstyle (BA)(CD)}; (14,-29)*+{\scriptstyle ((BA)C)D} };
{\ar^{a1} (14,-29)*+{\scriptstyle ((BA)C)D}; (38,-29)*+{\scriptstyle (B(AC))D} };
{\ar^{(1R)1} (38,-29)*+{\scriptstyle (B(AC))D}; (62,-29)*+{\scriptstyle (B(CA))D} };
{\ar^{a^{\centerdot}1} (62,-29)*+{\scriptstyle (B(CA))D}; (86,-29)*+{\scriptstyle ((BC)A)D} };
{\ar^{a} (86,-29)*+{\scriptstyle ((BC)A)D}; (93,-38)*+{\scriptstyle (BC)(AD)} };
{\ar^{1R} (93,-38)*+{\scriptstyle (BC)(AD)}; (100,-47)*+{\scriptstyle (BC)(DA)} };
{\ar_{a} (0,-47)*+{\scriptstyle (AB)(CD)}; (12,-61)*+{\scriptstyle A(B(CD))} };
{\ar_{R}  (12,-61)*+{\scriptstyle A(B(CD))}; (37,-61)*+{\scriptstyle (B(CD))A} };
{\ar_{a} (37,-61)*+{\scriptstyle (B(CD))A}; (63,-61)*+{\scriptstyle B((CD)A)} };
{\ar_{1a} (63,-61)*+{\scriptstyle B((CD)A)}; (88,-61)*+{\scriptstyle B(C(DA))} };
{\ar_{a^{\centerdot}} (88,-61)*+{\scriptstyle B(C(DA))}; (100,-47)*+{\scriptstyle (BC)(DA)} };
{\ar_{a^{\centerdot}} (0,-47)*+{\scriptstyle (AB)(CD)}; (29,-41)*+{\scriptstyle ((AB)C)D} };
{\ar_{(R1)1} (29,-41)*+{\scriptstyle ((AB)C)D}; (14,-29)*+{\scriptstyle ((BA)C)D} };
{\ar_{a1} (29,-41)*+{\scriptstyle ((AB)C)D}; (57,-35)*+{\scriptstyle (A(BC))D} };
{\ar_{R1} (57,-35)*+{\scriptstyle (A(BC))D}; (86,-29)*+{\scriptstyle ((BC)A)D} };
{\ar^{a} (57,-35)*+{\scriptstyle (A(BC))D}; (59,-43)*+{\scriptstyle A((BC)D)} };
{\ar^{R} (59,-43)*+{\scriptstyle A((BC)D)}; (61,-52)*+{\scriptstyle ((BC)D)A} };
{\ar^{1a^{\centerdot}} (12,-61)*+{\scriptstyle A(B(CD))}; (59,-43)*+{\scriptstyle A((BC)D)} };
{\ar_{a^{\centerdot}1} (37,-61)*+{\scriptstyle (B(CD))A}; (61,-52)*+{\scriptstyle ((BC)D)A} };
{\ar_{a} (61,-52)*+{\scriptstyle ((BC)D)A}; (100,-47)*+{\scriptstyle (BC)(DA)} };
(16,-40)*{\cong}; (40,-34)*{\Downarrow R_{(A|B,C)}1}; (28,-48)*{\Downarrow \pi}; (45,-52)*{\cong}; (77,-43)*{\Downarrow R_{(A|BC,D)}}; (76,-56)*{\Downarrow \pi}
\endxy
\]

\[
\xy
{\ar^{\scriptstyle (1R)1} (0,0)*+{\scriptscriptstyle (A(BC))D}; (0,12.6)*+{\scriptscriptstyle (A(CB))D} };
{\ar^{\scriptstyle a^{\centerdot}1} (0,12.6)*+{\scriptscriptstyle (A(CB))D}; (8,25.2)*+{\scriptscriptstyle ((AC)B)D} };
{\ar^{\scriptstyle a} (8,25.2)*+{\scriptscriptstyle ((AC)B)D}; (29,25.2)*+{\scriptscriptstyle (AC)(BD)} };
{\ar^{\scriptstyle R1} (29,25.2)*+{\scriptscriptstyle (AC)(BD)}; (50,25.2)*+{\scriptscriptstyle (CA)(BD)} };
{\ar^{\scriptstyle 1R} (50,25.2)*+{\scriptscriptstyle (CA)(BD)}; (71,25.2)*+{\scriptscriptstyle (CA)(DB)} };
{\ar^{\scriptstyle a} (71,25.2)*+{\scriptscriptstyle (CA)(DB)}; (92,25.2)*+{\scriptscriptstyle C(A(DB))} };
{\ar^{\scriptstyle 1a^{\centerdot}} (92,25.2)*+{\scriptscriptstyle C(A(DB))}; (100,12.6)*+{\scriptscriptstyle C((AD)B)} };
{\ar^{\scriptstyle 1(R1)} (100,12.6)*+{\scriptscriptstyle C((AD)B)}; (100,0)*+{\scriptscriptstyle C((DA)B)} };
{\ar_{\scriptstyle a^{\centerdot}1} (0,0)*+{\scriptscriptstyle (A(BC))D}; (17,-15.75)*+{\scriptscriptstyle ((AB)C)D} };
{\ar_{\scriptstyle a} (17,-15.75)*+{\scriptscriptstyle ((AB)C)D}; (39,-15.75)*+{\scriptscriptstyle (AB)(CD)} };
{\ar_{\scriptstyle R} (39,-15.75)*+{\scriptscriptstyle (AB)(CD)}; (61,-15.75)*+{\scriptscriptstyle (CD)(AB)} };
{\ar_{\scriptstyle a} (61,-15.75)*+{\scriptscriptstyle (CD)(AB)}; (83,-15.75)*+{\scriptscriptstyle C(D(AB))} };
{\ar_{\scriptstyle 1 a^{\centerdot}} (83,-15.75)*+{\scriptscriptstyle C(D(AB))}; (100,0)*+{\scriptscriptstyle C((DA)B)} };
{\ar^{\scriptstyle R1} (17,-15.75)*+{\scriptscriptstyle ((AB)C)D}; (23,-2.1)*+{\scriptscriptstyle (C(AB))D} };
{\ar^{\scriptstyle a^{\centerdot}1} (23,-2.1)*+{\scriptscriptstyle (C(AB))D}; (29,11.55)*+{\scriptscriptstyle ((CA)B)D} };
{\ar_{\scriptstyle (R1)1} (8,25.2)*+{\scriptscriptstyle ((AC)B)D}; (29,11.55)*+{\scriptscriptstyle ((CA)B)D} };
{\ar_{\scriptstyle a} (29,11.55)*+{\scriptscriptstyle ((CA)B)D}; (50,25.2)*+{\scriptscriptstyle (CA)(BD)} };
{\ar_{\scriptstyle a} (50,25.2)*+{\scriptscriptstyle (CA)(BD)}; (71,11.55)*+{\scriptscriptstyle C(A(BD))} };
{\ar_{\scriptstyle 1(1R)} (71,11.55)*+{\scriptscriptstyle C(A(BD))}; (92,25.2)*+{\scriptscriptstyle C(A(DB))} };
{\ar^{\scriptstyle 1a^{\centerdot}} (71,11.55)*+{\scriptscriptstyle C(A(BD))}; (60,-1.05)*+{\scriptscriptstyle C((AB)D)} };
{\ar^{\scriptstyle a} (23,-2.1)*+{\scriptscriptstyle (C(AB))D}; (60,-1.05)*+{\scriptscriptstyle C((AB)D)} };
{\ar^{\scriptstyle 1R} (60,-1.05)*+{\scriptscriptstyle C((AB)D)}; (83,-15.75)*+{\scriptscriptstyle C(D(AB))} };
(12,8)*{\scriptstyle \Downarrow R_{(A,B|C)}1}; (29,18)*{\scriptstyle \cong}; (50,7)*{\scriptstyle \Downarrow \pi}; (71,18)*{\scriptstyle \cong}; (83,6)*{\scriptstyle \Downarrow 1 R_{(A,B|D)}}; (50,-8)*{\scriptstyle \Downarrow R_{(AB|C,D)}};
{\ar@{=} (50,-20.5)*{}; (50,-27.5)*{} };
{\ar^{\scriptstyle (1R)1} (0,-64)*+{\scriptscriptstyle (A(BC))D}; (0,-48)*+{\scriptscriptstyle (A(CB))D} };
{\ar^{\scriptstyle a^{\centerdot}1} (0,-48)*+{\scriptscriptstyle (A(CB))D}; (8,-30)*+{\scriptscriptstyle ((AC)B)D} };
{\ar^{\scriptstyle a} (8,-30)*+{\scriptscriptstyle ((AC)B)D}; (29,-30)*+{\scriptscriptstyle (AC)(BD)} };
{\ar^{\scriptstyle R1} (29,-30)*+{\scriptscriptstyle (AC)(BD)}; (50,-30)*+{\scriptscriptstyle (CA)(BD)} };
{\ar^{\scriptstyle 1R} (50,-30)*+{\scriptscriptstyle (CA)(BD)}; (71,-30)*+{\scriptscriptstyle (CA)(DB)} };
{\ar^{\scriptstyle a} (71,-30)*+{\scriptscriptstyle (CA)(DB)}; (92,-30)*+{\scriptscriptstyle C(A(DB))} };
{\ar^{\scriptstyle 1a^{\centerdot}} (92,-30)*+{\scriptscriptstyle C(A(DB))}; (100,-48)*+{\scriptscriptstyle C((AD)B)} };
{\ar^{\scriptstyle 1(R1)} (100,-48)*+{\scriptscriptstyle C((AD)B)}; (100,-64)*+{\scriptscriptstyle C((DA)B)} };
{\ar_{\scriptstyle a^{\centerdot}1} (0,-64)*+{\scriptscriptstyle (A(BC))D}; (17,-82)*+{\scriptscriptstyle ((AB)C)D} };
{\ar_{\scriptstyle a} (17,-82)*+{\scriptscriptstyle ((AB)C)D}; (39,-82)*+{\scriptscriptstyle (AB)(CD)} };
{\ar_{\scriptstyle R} (39,-82)*+{\scriptscriptstyle (AB)(CD)}; (61,-82)*+{\scriptscriptstyle (CD)(AB)} };
{\ar_{\scriptstyle a} (61,-82)*+{\scriptscriptstyle (CD)(AB)}; (83,-82)*+{\scriptscriptstyle C(D(AB))} };
{\ar_{\scriptstyle 1 a^{\centerdot}} (83,-82)*+{\scriptscriptstyle C(D(AB))}; (100,-64)*+{\scriptscriptstyle C((DA)B)} };
{\ar^{\scriptstyle a} (0,-64)*+{\scriptscriptstyle (A(BC))D}; (15,-57)*+{\scriptscriptstyle A((BC)D)} };
{\ar^{\scriptstyle 1(R1)} (15,-57)*+{\scriptscriptstyle A((BC)D)}; (15,-45)*+{\scriptscriptstyle A((CB)D)} };
{\ar_{\scriptstyle a} (0,-48)*+{\scriptscriptstyle (A(CB))D}; (15,-45)*+{\scriptscriptstyle A((CB)D)} };
{\ar_{\scriptstyle 1a} (15,-45)*+{\scriptscriptstyle A((CB)D)}; (26,-38)*+{\scriptscriptstyle A(C(BD))} };
{\ar_{\scriptstyle a^{\centerdot}} (26,-38)*+{\scriptscriptstyle A(C(BD))}; (29,-30)*+{\scriptscriptstyle (AC)(BD)} };
{\ar_{\scriptstyle 1R} (29,-30)*+{\scriptscriptstyle (AC)(BD)}; (50,-42)*+{\scriptscriptstyle (AC)(DB)} };
{\ar_{\scriptstyle R1} (50,-42)*+{\scriptscriptstyle (AC)(DB)}; (71,-30)*+{\scriptscriptstyle (CA)(DB)} };
{\ar_{\scriptstyle a^{\centerdot}} (71,-30)*+{\scriptscriptstyle (CA)(DB)}; (68,-49)*+{\scriptscriptstyle ((CA)D)B} };
{\ar_{\scriptstyle a1} (68,-49)*+{\scriptscriptstyle ((CA)D)B}; (84,-45)*+{\scriptscriptstyle (C(AD))B} };
{\ar_{\scriptstyle a} (84,-45)*+{\scriptscriptstyle (C(AD))B}; (100,-48)*+{\scriptscriptstyle C((AD)B)} };
{\ar^{\scriptstyle (1R)1} (84,-45)*+{\scriptscriptstyle (C(AD))B}; (84,-60)*+{\scriptscriptstyle (C(DA))B} };
{\ar_{\scriptstyle a} (84,-60)*+{\scriptscriptstyle (C(DA))B}; (100,-64)*+{\scriptscriptstyle C((DA)B)} };
{\ar_{\scriptstyle a^{\centerdot}} (61,-82)*+{\scriptscriptstyle (CD)(AB)}; (75,-71)*+{\scriptscriptstyle ((CD)A)B} };
{\ar_{\scriptstyle a1} (75,-71)*+{\scriptscriptstyle ((CD)A)B};  (84,-60)*+{\scriptscriptstyle (C(DA))B} };
{\ar_{\scriptstyle 1(1R)} (26,-38)*+{\scriptscriptstyle A(C(BD))}; (38,-50)*+{\scriptscriptstyle A(C(DB))} };
{\ar_{\scriptstyle a^{\centerdot}} (38,-50)*+{\scriptscriptstyle A(C(DB))}; (50,-42)*+{\scriptscriptstyle (AC)(DB)} };
{\ar_{\scriptstyle a^{\centerdot}} (50,-42)*+{\scriptscriptstyle (AC)(DB)}; (60,-59)*+{\scriptscriptstyle ((AC)D)B} };
{\ar_{\scriptstyle (R1)1} (60,-59)*+{\scriptscriptstyle ((AC)D)B}; (68,-49)*+{\scriptscriptstyle ((CA)D)B} };
{\ar^{\scriptstyle 1a} (15,-57)*+{\scriptscriptstyle A((BC)D)}; (24.5,-70)*+{\scriptscriptstyle A(B(CD))} };
{\ar^{\scriptstyle a^{\centerdot}} (24.5,-70)*+{\scriptscriptstyle A(B(CD))}; (39,-82)*+{\scriptscriptstyle (AB)(CD)} };
{\ar_{\scriptstyle 1R} (24.5,-70)*+{\scriptscriptstyle A(B(CD))}; (32,-59)*+{\scriptscriptstyle A((CD)B)} };
{\ar_{\scriptstyle 1a} (32,-59)*+{\scriptscriptstyle A((CD)B)}; (38,-50)*+{\scriptscriptstyle A(C(DB))} };
{\ar_{\scriptstyle a^{\centerdot}} (32,-59)*+{\scriptscriptstyle A((CD)B)}; (50,-71)*+{\scriptscriptstyle (A(CD))B} };
{\ar_{\scriptstyle R1} (50,-71)*+{\scriptscriptstyle (A(CD))B}; (75,-71)*+{\scriptscriptstyle ((CD)A)B} };
{\ar_{\scriptstyle a^{\centerdot}1} (50,-71)*+{\scriptscriptstyle (A(CD))B}; (60,-59)*+{\scriptscriptstyle ((AC)D)B} };
(7,-54)*{\scriptstyle \cong};(11,-38)*{\scriptstyle \Downarrow \pi}; (50,-35)*{\scriptstyle \cong}; (35,-42)*{\scriptstyle \cong}; (60,-45)*{\scriptstyle \cong}; (84,-38)*{\scriptstyle \Downarrow \pi}; (91,-55)*{\scriptstyle \cong}; (83,-76)*{\scriptstyle \Downarrow \pi}; (13,-68)*{\scriptstyle \Downarrow \pi}; (47,-60)*{\scriptstyle \Downarrow \pi}; (25,-53)*{\scriptstyle \Downarrow 1 R_{(B|C,D)}}; (50,-76)*{\scriptstyle \Downarrow R_{(A,B|CD)}}; (69,-64)*{\scriptstyle \Downarrow R_{(A|C,D)}1}
\endxy
\]

\[
\xy
{\ar^{R1} (0,0)*+{(AB)C}; (10,15)*+{(BA)C} };
{\ar^{a} (10,15)*+{(BA)C}; (30,15)*+{B(AC)} };
{\ar^{1R} (30,15)*+{B(AC)}; (50,15)*+{B(CA)} };
{\ar^{a^{\centerdot}} (50,15)*+{B(CA)}; (70,15)*+{(BC)A} };
{\ar^{R1} (70,15)*+{(BC)A}; (90,15)*+{(CB)A} };
{\ar^{a} (90,15)*+{(CB)A}; (100,0)*+{C(BA)} };
{\ar_{a} (0,0)*+{(AB)C}; (10,-15)*+{A(BC)} };
{\ar_{1R} (10,-15)*+{A(BC)}; (30,-15)*+{A(CB)} };
{\ar_{a^{\centerdot}} (30,-15)*+{A(CB)}; (50,-15)*+{(AC)B} };
{\ar_{R1} (50,-15)*+{(AC)B}; (70,-15)*+{(CA)B} };
{\ar_{a} (70,-15)*+{(CA)B}; (90,-15)*+{C(AB)} };
{\ar_{1R} (90,-15)*+{C(AB)}; (100,0)*+{C(BA)} };
{\ar^{R} (10,-15)*+{A(BC)}; (30,0)*+{(BC)A} };
{\ar^{a}  (30,0)*+{(BC)A}; (50,15)*+{B(CA)} };
{\ar_{R1} (30,0)*+{(BC)A}; (90,15)*+{(CB)A} };
{\ar_{R} (30,-15)*+{A(CB)}; (90,15)*+{(CB)A} };
(18,3.75)*{\scriptstyle \Downarrow R_{(A|B,C)}}; (58,10)*{\cong}; (40,-3.75)*{\cong}; (75,-3.75)*{\scriptstyle \Downarrow R_{(A|C,B)}^{-1}};
{\ar@{=} (50,-19.5)*{}; (50,-31.5)*{} };
{\ar^{R1} (0,-52.5)*+{(AB)C}; (10,-37.5)*+{(BA)C} };
{\ar^{a} (10,-37.5)*+{(BA)C}; (30,-37.5)*+{B(AC)} };
{\ar^{1R} (30,-37.5)*+{B(AC)}; (50,-37.5)*+{B(CA)} };
{\ar^{a^{\centerdot}} (50,-37.5)*+{B(CA)}; (70,-37.5)*+{(BC)A} };
{\ar^{R1} (70,-37.5)*+{(BC)A}; (90,-37.5)*+{(CB)A} };
{\ar^{a} (90,-37.5)*+{(CB)A}; (100,-52.5)*+{C(BA)} };
{\ar_{a} (0,-52.5)*+{(AB)C}; (10,-67.5)*+{A(BC)} };
{\ar_{1R} (10,-67.5)*+{A(BC)}; (30,-67.5)*+{A(CB)} };
{\ar_{a^{\centerdot}} (30,-67.5)*+{A(CB)}; (50,-67.5)*+{(AC)B} };
{\ar_{R1} (50,-67.5)*+{(AC)B}; (70,-67.5)*+{(CA)B} };
{\ar_{a} (70,-67.5)*+{(CA)B}; (90,-67.5)*+{C(AB)} };
{\ar_{1R} (90,-67.5)*+{C(AB)}; (100,-52.5)*+{C(BA)} };
{\ar_{R} (0,-52.5)*+{(AB)C}; (90,-67.5)*+{C(AB)} };
{\ar^{a^{\centerdot}} (30,-37.5)*+{B(AC)}; (65,-52.5)*+{(BA)C} };
{\ar_{R} (65,-52.5)*+{(BA)C}; (100,-52.5)*+{C(BA)} };
{\ar^{R1} (0,-52.5)*+{(AB)C}; (65,-52.5)*+{(BA)C} };
(20,-45)*{\cong}; (17,-61.25)*{\scriptstyle \Downarrow R_{(A,B|C)}^{-1}}; (70,-45)*{\scriptstyle \Downarrow R_{(B,A|C)}}; (70,-57.5)*{\cong}
\endxy
\]
A \textit{braided monoidal bicategory} is a monoidal bicategory equipped with a braiding.
\end{definition}

While these axioms might look quite daunting, they are in fact just algebraic expressions of the notion that a pair of homotopies that each start at the braid $\gamma$ and end at the braid $\gamma'$ are in fact themselves homotopic.  For instance, the first axiom concerns the case of a braid with four strands in which the first three strands are braided past the final one.  This can be done either all at once (this is the 1-cell target of the pasting diagram) or it can be done step-by-step in which strand 3 is braided past strand 4, then strand 2 is braided past strand 4, and finally strand 1 is braided past strand 4 (this is the 1-cell source).  The two different composite 2-cells which are claimed to be equal in this axiom are just two different ways to transform the step-by-step method into the all-at-once method using the algebra available in a braided monoidal bicategory.  The other three axioms also have similar interpretations.  In fact, one could take the presentation of these axioms seriously, and view them as certain three-dimensional polytopes in which the two-dimensional faces are precisely the 2-cells in each equation above.  Doing so produces polytopes discovered by Bar-Natan in \cite{barn1}.

\begin{definition}
Let $B$ be a braided monoidal bicategory.  A \textit{syllepsis} for $B$ consists of an invertible modification with components $v_{xy}:R_{yx} \circ R_{xy} \Rightarrow 1$ such that the following two axioms hold; here we use the convention that if $\alpha$ is a 2-cell, then $\widehat{\alpha}$ denotes its mate.
\[
\def\objectstyle{\scriptstyle}
\def\labelstyle{\scriptstyle}
\xy0;/r.20pc/:
{\ar_{R1} (0,0)*+{(ab)c}; (15,10)*+{(ba)c} };
{\ar@/^1.5pc/^{R^{\dot}1} (0,0)*+{(ab)c}; (15,10)*+{(ba)c} };
{\ar^{a} (15,10)*+{(ba)c}; (35,10)*+{b(ac)} };
{\ar_{1R} (35,10)*+{b(ac)}; (50,0)*+{b(ca)} };
{\ar@/^1.5pc/^{1R^{\dot}} (35,10)*+{b(ac)}; (50,0)*+{b(ca)} };
{\ar_{a} (0,0)*+{(ab)c}; (15,-10)*+{a(bc)} };
{\ar_{R} (15,-10)*+{a(bc)}; (35,-10)*+{(bc)a} };
{\ar_{a} (35,-10)*+{(bc)a}; (50,0)*+{b(ca)} };
{\ar@{=>}^{R_{(a|b,c)}} (25, 5)*{}; (25,-5)*{} };
(60,0)*{=}; (6,8.5)*{\scriptscriptstyle \Downarrow \widehat{v}_{ab}1};(44,8.5)*{\scriptscriptstyle \Downarrow 1\widehat{v}_{ac}};
{\ar^{R^{\dot}1} (70,0)*+{(ab)c}; (85,10)*+{(ba)c} };
{\ar^{a} (85,10)*+{(ba)c}; (105,10)*+{b(ac)} };
{\ar^{1R^{\dot}} (105,10)*+{b(ac)}; (120,0)*+{b(ca)} };
{\ar_{a} (70,0)*+{(ab)c}; (85,-10)*+{a(bc)} };
{\ar^{R^{\dot}} (85,-10)*+{a(bc)}; (105,-10)*+{(bc)a} };
{\ar@/_1.5pc/_{R} (85,-10)*+{a(bc)}; (105,-10)*+{(bc)a} };
{\ar_{a} (105,-10)*+{(bc)a}; (120,0)*+{b(ca)} };
{\ar@{=>}^{\widehat{R}_{(a,b|c)}} (95, 5)*{}; (95,-5)*{} };
(95,-14)*{\scriptscriptstyle \Downarrow \widehat{v}_{a,bc}}
\endxy
\]
\[
\def\objectstyle{\scriptstyle}
\def\labelstyle{\scriptstyle}
\xy0;/r.20pc/:
{\ar_{1R} (0,0)*+{a(bc)}; (15,10)*+{a(cb)} };
{\ar@/^1.5pc/^{1R^{\dot}} (0,0)*+{a(bc)}; (15,10)*+{a(cb)} };
{\ar^{a^{\centerdot}} (15,10)*+{a(cb)}; (35,10)*+{(ac)b} };
{\ar_{R1} (35,10)*+{(ac)b}; (50,0)*+{(ca)b} };
{\ar@/^1.5pc/^{R^{\dot}1} (35,10)*+{(ac)b}; (50,0)*+{(ca)b} };
{\ar_{a^{\centerdot}} (0,0)*+{a(bc)}; (15,-10)*+{(ab)c} };
{\ar_{R} (15,-10)*+{(ab)c}; (35,-10)*+{c(ab)} };
{\ar_{a^{\centerdot}} (35,-10)*+{c(ab)}; (50,0)*+{(ca)b} };
{\ar@{=>}^{R_{(a,b|c)}} (25, 7)*{}; (25,-7)*{} };
(60,0)*{=}; (6,8.5)*{\scriptscriptstyle \Downarrow 1 \widehat{v}_{bc}};(44,8.5)*{\scriptscriptstyle \Downarrow \widehat{v}_{ac}1};
{\ar^{1R} (70,0)*+{a(bc)}; (85,10)*+{a(cb)} };
{\ar^{a^{\centerdot}} (85,10)*+{a(cb)}; (105,10)*+{(ac)b} };
{\ar^{R1} (105,10)*+{(ac)b}; (120,0)*+{(ca)b} };
{\ar_{a^{\centerdot}} (70,0)*+{a(bc)}; (85,-10)*+{(ab)c} };
{\ar^{R^{\dot}} (85,-10)*+{(ab)c}; (105,-10)*+{c(ab)} };
{\ar@/_1.5pc/_{R} (85,-10)*+{(ab)c}; (105,-10)*+{c(ab)} };
{\ar_{a^{\centerdot}} (105,-10)*+{c(ab)}; (120,0)*+{(ca)b} };
{\ar@{=>}^{\widehat{R}_{(a|b,c)}} (95, 5)*{}; (95,-5)*{} };
(95,-14)*{\scriptscriptstyle \Downarrow \widehat{v}_{ab,c}}
\endxy
\]

A sylleptic monoidal bicategory is a braided monoidal bicategory equipped with a syllepsis.
\end{definition}

\begin{definition}
A symmetric monoidal bicategory is a sylleptic monoidal bicategory in which the following axiom holds for every pair of objects $a,b$.
\[
\xy
{\ar^{R} (0,0)*+{ab}; (10,15)*+{ba} };
{\ar^{R} (10,15)*+{ba}; (30,15)*+{ab} };
{\ar^{R} (30,15)*+{ab}; (40,0)*+{ba} };
{\ar_{R} (0,0)*+{ab}; (40,0)*+{ba} };
{\ar_{1} (0,0)*+{ab}; (30,15)*+{ab} };
(13,10.5)*{\Downarrow v_{ab}}; (25,5)*{\cong};
{\ar^{R} (50,0)*+{ab}; (60,15)*+{ba} };
{\ar^{R} (60,15)*+{ba}; (80,15)*+{ab} };
{\ar^{R} (80,15)*+{ab}; (90,0)*+{ba} };
{\ar_{R} (50,0)*+{ab}; (90,0)*+{ba} };
{\ar_{1} (60,15)*+{ba}; (90,0)*+{ba} };
(65,5)*{\cong}; (77,10.5)*{\Downarrow v_{ba}};
(45,7.5)*{=}
\endxy
\]

\end{definition}

\begin{remark}
The proofs in the next section rely heavily on various coherence theorems that we will not state here; see \cite{gur3, gur4, go1}. Essentially, these say that ``every diagram of 2-cell constraints commutes''.  In practice, this means we can simplify calculations involving 2-cells in the following ways. 
\begin{itemize}
 \item We can disregard lower-dimensional coherence cells in the sources and targets of 2-cell diagrams, as there is a unique coherence isomorphism between any two possible interpretations.
\item We can manipulate coherence 2-cells in the diagram, more-or-less at will.
\end{itemize}

\end{remark}

\section{Generalised icons}\label{genicons}

In this section we are going to generalise the definition of icon in the following sense.  Given a monoidal bicategory $\cl{B}$ we will define a bicategory \cat{$\cl{B}$-Icon} with
\begin{itemize}
 \item 0-cells: categories weakly enriched in \cl{B}
\item 1-cells: weak functors
\item 2-cells: icons enriched in \cl{B}.
\end{itemize}
The case $\cl{B}=\cat{Cat}$ gives the original definition Lack in \cite{lac2}; Shulman gives an abstract account of the case where \cl{B} is a tensor distributive monoidal 2-category in \cite{shul2}.  Furthermore we will show that if \cl{B} is symmetric monoidal then \cat{\cl{B}-Icon} is also symmetric monoidal, so the construction can be iterated.  

We begin by defining the appropriate bicategory of graphs.  For the rest of this section we fix a monoidal bicategory \cl{B}.  Motivating examples include \cat{Cat} and \cat{Icon}. 

\begin{mydefinition} A \demph{\cl{B}-graph} $X$ is given by 
\begin{itemize}
 \item a set $X_0$ of objects, and
\item for all objects $x,y \in X_0$ a hom-object $X(x,y) \in \cl{B}$.
\end{itemize}
A \demph{morphism of \cl{B}-graphs} $F:X \lra Y$ is given by
\begin{itemize}
 \item a function $F: X_0 \lra Y_0$, and
\item for all objects $x,y \in X_0$ a morphism $F=F_{x,y}:X(x,y) \lra Y(Fx, Fy) \in \cl{B}$.
\end{itemize}
A \demph{2-morphism of \cl{B}-graphs} $\alpha:F \Longrightarrow G$ 
\[
\xy
(-10,0)*+{X}="1";
(10,0)*+{Y}="2";
{\ar@/^1.7pc/^{F} "1";"2"};
{\ar@/_1.7pc/_{G} "1";"2"};
{\ar@{=>}^<<<{\alpha} (0,3)*{};(0,-3)*{}} ;
\endxy
\]
exists only when $Fx=Gx$ for all $x \in X_0$.  When this holds, $\alpha$ is given by, for every pair $x,y \in X_0$ a 2-cell
\[
\xy
(-20,0)*+{X(x,y)}="1";
(20,0)*+{\hspace*{6em}Y(Fx,Fy)=Y(Gx,Gy)}="2";
{\ar@/^1.9pc/^{F} "1";"2"};
{\ar@/_1.9pc/_{G} "1";"2"};
{\ar@{=>}^<<<{\alpha_{x,y}} (0,3)*{};(0,-3)*{}} ;
\endxy
\]

It is straightforward to check that \cl{B}-graphs, their morphisms and their 2-morphisms form a bicategory, which we write as \cat{\cl{B}-Gph}.  Note that horizontal and vertical composition of 2-morphisms comes from the horizontal and vertical composition of 2-cells in \cl{B}.

\end{mydefinition}

\begin{mydefinition}\label{tractable}\label{weakbcat}

Let \cl{B} be a braided monoidal bicategory. A \demph{category weakly enriched in \cl{B}} or \demph{weak \cl{B}-category} $X$ is given by 

\begin{itemize}
 \item an underlying $\cl{B}$-graph $X$,
\item composition: for all $a,b,c \in X_0$ a 1-cell in \cl{B}
\[m_{abc} : X(b,c) \otimes X(a,b) \lra X(a,c),\]
\item identities: for all $a \in X_0$ a 1-cell in \cl{B}
\[I_a: \1 \lra X(a,a),\]

\item associativity constraints: for all $a,b,c,d \in X_0$ an invertible 2-cell in \cl{B}

 \[
 \xy
    (-30,28)*+{ \big(X(c,d) \otimes  X(b,c)\big) \otimes  X(a,b)}="tl";
    (-30,15)*+{ X(c,d) \otimes  \big(X(b,c) \otimes  X(a,b)\big)}="ml";
        {\ar_{\catequiv} "tl";"ml"};
(-30,-10)*+{X(c,d) \otimes  X(a,c)}="tr";
    (30,28)*+{X(b,d)\otimes  X(a,b)}="bl";
    (30,-10)*+{X(a,d)}="br";
        {\ar_-{1 \otimes m_{a,b,c}} "ml";"tr"};
        {\ar^-{m_{b,c,d} \otimes 1} "tl";"bl"};
        {\ar^-{m_{a,b,d}} "bl";"br"};
        {\ar_-{m_{a,c,d}} "tr";"br"};
    {\ar@{=>}^{\;a_{a,b,c,d}} (8,12);(-6,2)};
 \endxy
 \]
and

 \item unit constraints: for all $a$, $b$ invertible 2-cells in \cl{B}
  \[\xy
  (20,20)*+{X(a,b) \otimes X(a,a)}="tl";
  (-20,20)*+{X(a,b) \otimes \1 }="tr";
  (20,0)*+{ X(a,b)}="bl";
    {\ar^{m_{a,a,b}} "tl";"bl"};
    {\ar^{1\otimes I_a} "tr";"tl"};
    {\ar_{\catequiv} "tr";"bl"};
 {\ar@{=>}_{r_{a,b}} (15,15);(9,9)};
  \endxy
 \qquad \quad
  \xy
  (-20,20)*+{X(b,b) \otimes X(a,b)}="tl";
  (20,20)*+{ \1 \otimes X(a,b)}="tr";
  (-20,0)*+{ X(a,b)}="bl";
    {\ar_{m_{a,b,b}} "tl";"bl"};
    {\ar_{I_b \otimes 1} "tr";"tl"};
    {\ar^{\catequiv} "tr";"bl"};
{\ar@{=>}^{\ell_{a,b}} (-15,15);(-9,9)};
  \endxy.\]

\end{itemize}

This data satisfies the following axioms.  Here we write $X_{ab}$ for $X(a,b)$ and omit all $\otimes$ symbols.  Note also that the 2-cell $\mbox{\bf{a}}$ is the associativity constraint for the weak $\cB$-category $X$, whereas the 1-cell $a$ is the associativity constraint for the monoidal structure of $\cB$.



\newpage

\[\psset{unit=0.16cm,labelsep=0pt,nodesep=3pt}
\pspicture(0,10)(80,60)


\rput(40,13){\rnode{a4}{$X_{de}(X_{cd}(X_{bc}X_{ab}))$}} 

\rput(18,23){\rnode{b2}{$(X_{de}X_{cd})(X_{bc}X_{ab})$}} 
\rput(62,23){\rnode{b6}{$X_{de}(X_{cd}X_{ac})$}} 

\rput(40,30){\rnode{c4}{$X_{de}((X_{cd}X_{bc})X_{ab})$}} 

\rput(4,40){\rnode{d1}{$((X_{de}X_{cd})X_{bc})X_{ab}$}} 
\rput(30,40){\rnode{d3}{$(X_{de}(X_{cd}X_{bc}))X_{ab}$}} 
\rput(50,40){\rnode{d5}{$X_{de}(X_{bd}X_{ab})$}} 
\rput(76,40){\rnode{d7}{$X_{de}X_{ad}$}} 

\rput(40,50){\rnode{e4}{$(X_{de}X_{bd})X_{ab}$}} 

\rput(18,57){\rnode{f2}{$(X_{ce}X_{bc})X_{ab}$}} 
\rput(62,57){\rnode{f6}{$X_{ac}$}} 

\rput(40,67){\rnode{g4}{$X_{be}X_{ab}$}}

\psset{nodesep=3pt,labelsep=2pt,arrows=->}

\ncline{d1}{f2}\naput{{\scriptsize $(m1)1$}}  
\ncline{f2}{g4}\naput{{\scriptsize $m1$}} 
\ncline{g4}{f6}\naput{{\scriptsize $m$}} 
\ncline{d7}{f6}\nbput{{\scriptsize $m$}} 

\ncline{d1}{b2}\nbput{{\scriptsize $a$}}  
\ncline{b2}{a4}\nbput{{\scriptsize $a$}} 
\ncline{a4}{b6}\nbput{{\scriptsize $1(1m)$}} 
\ncline{b6}{d7}\nbput{{\scriptsize $1m$}} 

\ncline{d3}{e4}\naput[labelsep=0pt]{{\scriptsize $(1m)1$}}  
\ncline{e4}{d5}\naput{{\scriptsize $a$}} 
\ncline{d3}{c4}\nbput{{\scriptsize $a$}} 
\ncline{c4}{d5}\nbput[labelsep=0pt]{{\scriptsize $1(m1)$}} 

\ncline{d1}{d3}\naput{{\scriptsize $a1$}}  
\ncline{e4}{g4}\naput{{\scriptsize $m1$}} 
\ncline{d5}{d7}\naput{{\scriptsize $1m$}} 
\ncline{c4}{a4}\naput{{\scriptsize $1a$}} 

\rput[c](40,43){{\scriptsize $\iso$}}

{
\rput[c](25,48){\psset{unit=1mm,doubleline=true,arrowinset=0.6,arrowlength=0.5,arrowsize=0.5pt 2.1,nodesep=0pt,labelsep=2pt}
\pcline{->}(0,4)(0,0) \naput{{\scriptsize {\bfseries a}$1$}}}}

{
\rput[c](55,30){\psset{unit=1mm,doubleline=true,arrowinset=0.6,arrowlength=0.5,arrowsize=0.5pt 2.1,nodesep=0pt,labelsep=2pt}
\pcline{->}(0,4)(0,0) \naput{{\scriptsize $1${\bfseries a}}}}}

{
\rput[c](25,30){\psset{unit=1mm,doubleline=true,arrowinset=0.6,arrowlength=0.5,arrowsize=0.5pt 2.1,nodesep=0pt,labelsep=2pt}
\pcline{->}(0,4)(0,0) \naput{{\scriptsize {$\Pi$}}}}}

{
\rput[c](55,48){\psset{unit=1mm,doubleline=true,arrowinset=0.6,arrowlength=0.5,arrowsize=0.5pt 2.1,nodesep=0pt,labelsep=2pt}
\pcline{->}(0,4)(0,0) \naput{{\scriptsize {\bfseries a}}}}}

\endpspicture\]

$=$

\[\psset{unit=0.16cm,labelsep=0pt,nodesep=3pt}
\pspicture(80,70)


\rput(40,13){\rnode{a4}{$X_{de}(X_{cd}(X_{bc}X_ab))$}} 

\rput(18,23){\rnode{b2}{$(X_{de}X_{cd})(X_{bc}X_{ab})$}} 
\rput(62,23){\rnode{b6}{$X_{de}(X_{cd}X_{ac})$}} 

\rput(50,33){\rnode{c5}{$(X_{de}X_{cd})X_{ac}$}} 

\rput(4,40){\rnode{d1}{$((X_{de}X_{cd})X_{bc})X_{ab}$}} 
\rput(76,40){\rnode{d7}{$X_{de}X_{ad}$}} 

\rput(30,47){\rnode{e3}{$X_{ce}(X_{bc}X_{ab})$}} 
\rput(50,47){\rnode{e5}{$X_{ce}X_{ac}$}}

\rput(18,57){\rnode{f2}{$(X_{ce}X_{bc})X_{ab}$}} 
\rput(62,57){\rnode{f6}{$X_{ac}$}} 

\rput(40,67){\rnode{g4}{$X_{be}X_{ab}$}}

\psset{nodesep=3pt,labelsep=2pt,arrows=->}

\ncline{d1}{f2}\naput{{\scriptsize $(m1)1$}}  
\ncline{f2}{g4}\naput{{\scriptsize $m1$}} 
\ncline{g4}{f6}\naput{{\scriptsize $m$}} 
\ncline{d7}{f6}\nbput{{\scriptsize $m$}} 

\ncline{d1}{b2}\nbput{{\scriptsize $a$}}  
\ncline{b2}{a4}\nbput{{\scriptsize $a$}} 
\ncline{a4}{b6}\nbput{{\scriptsize $1(1m)$}} 
\ncline{b6}{d7}\nbput{{\scriptsize $1m$}} 

\ncline{b2}{e3}\naput[labelsep=1pt]{{\scriptsize $m(11)$}}  
\ncline{e3}{e5}\naput{{\scriptsize $1m$}} 
\ncline{b2}{c5}\nbput[labelsep=1pt]{{\scriptsize $(11)m$}} 
\ncline{c5}{e5}\nbput{{\scriptsize $1m$}} 

\ncline{f2}{e3}\naput{{\scriptsize $a$}}  
\ncline{e5}{f6}\naput{{\scriptsize $m$}} 
\ncline{c5}{b6}\naput{{\scriptsize $a$}} 

\rput[c](39,41){{\scriptsize $\iso$}}
\rput[c](39,39){{\scriptsize\textsf{functoriality of $\otimes$}}}

\rput[c](18,46){{\scriptsize $\iso$}}
\rput[c](18,44){{\scriptsize\textsf{naturality of $a$}}}

\rput[c](45,24){{\scriptsize $\iso$}}
\rput[c](45,22){{\scriptsize\textsf{naturality of $a$} }}

{
\rput[c](40,55){\psset{unit=1mm,doubleline=true,arrowinset=0.6,arrowlength=0.5,arrowsize=0.5pt 2.1,nodesep=0pt,labelsep=2pt}
\pcline{->}(0,4)(0,0) \naput{{\scriptsize {\bfseries a}}}}}



{
\rput[c](60,39){\psset{unit=1mm,doubleline=true,arrowinset=0.6,arrowlength=0.5,arrowsize=0.5pt 2.1,nodesep=0pt,labelsep=2pt}
\pcline{->}(0,3)(3,0) \naput{{\scriptsize {\bfseries a}}}}}

\endpspicture\]

Note that ignoring $\cB$ coherence this axiom becomes (more tractably):

\[\psset{unit=0.16cm,labelsep=2pt,nodesep=3pt}
\pspicture(0,10)(50,45)


\rput(20,12){\rnode{a2}{$(d,e)(a,d)$}}  
\rput(46,12){\rnode{a4}{$(a,e)$}}  
\rput(3,23){\rnode{b1}{$(d,e)(c,d)(a,c)$}}  
\rput(20,30){\rnode{c2}{$(d,e)(b,d)(a,b)$}}  
\rput(46,30){\rnode{c4}{$(b,e)(a,b)$}}  
\rput(3,40){\rnode{d1}{$(d,e)(c,d)(b,c)(a,b)$}}  
\rput(28,40){\rnode{d3}{$(c,e)(b,c)(a,b)$}}

\ncline{->}{d1}{d3} \naput{{\scriptsize $m11$}} 
\ncline{->}{d3}{c4} \naput{{\scriptsize $m1$}}
\ncline{->}{d1}{c2} \naput{{\scriptsize $1m1$}}
\ncline{->}{c2}{c4} \naput{{\scriptsize $m1$}} 

\ncline{->}{d1}{b1} \nbput{{\scriptsize $11m$}}
\ncline{->}{b1}{a2} \nbput{{\scriptsize $1m$}}
\ncline{->}{c2}{a2} \naput{{\scriptsize $1m$}} 

\ncline{->}{c4}{a4} \naput{{\scriptsize $m$}} 
\ncline{->}{a2}{a4} \nbput[labelsep=3pt]{{\scriptsize $m$}}


{
\rput[c](25,34){\psset{unit=1mm,doubleline=true,arrowinset=0.6,arrowlength=0.5,arrowsize=0.5pt 2.1,nodesep=0pt,labelsep=1pt}
\pcline{->}(3,3)(0,0) \naput{{\scriptsize {\bfseries a}1}}}}

{
\rput[c](11,25){\psset{unit=1mm,doubleline=true,arrowinset=0.6,arrowlength=0.5,arrowsize=0.5pt 2.1,nodesep=0pt,labelsep=2pt}
\pcline{->}(3,3)(0,0) \naput{{\scriptsize 1{\bfseries a}}}}}

{
\rput[c](33,20){\psset{unit=1mm,doubleline=true,arrowinset=0.6,arrowlength=0.5,arrowsize=0.5pt 2.1,nodesep=0pt,labelsep=2pt}
\pcline{->}(3,3)(0,0) \naput{{\scriptsize {\bfseries a}}}}}

\endpspicture\]

$=$

\[\psset{unit=0.16cm,labelsep=2pt,nodesep=3pt}
\pspicture(50,45)


\rput(20,12){\rnode{a2}{$(d,e)(a,d)$}}  
\rput(46,12){\rnode{a4}{$(a,e)$}}  
\rput(3,23){\rnode{b1}{$(d,e)(c,d)(a,c)$}}  

\rput(28,23){\rnode{b3}{$(c,e)(a,c)$}}  

\rput(46,30){\rnode{c4}{$(b,e)(a,b)$}}

\rput(3,40){\rnode{d1}{$(d,e)(c,d)(b,c)(a,b)$}}  
\rput(28,40){\rnode{d3}{$(c,e)(b,c)(a,b)$}}

\ncline{->}{d1}{d3} \naput{{\scriptsize $m11$}} 
\ncline{->}{d3}{b3} \naput{{\scriptsize $1m$}}
\ncline{->}{d1}{b1} \nbput{{\scriptsize $11m$}}
\ncline{->}{b1}{b3} \naput{{\scriptsize $m1$}}

\ncline{->}{d3}{c4} \naput{{\scriptsize $m1$}} 
\ncline{->}{c4}{a4} \naput{{\scriptsize $m$}}
\ncline{->}{b3}{a4} \naput{{\scriptsize $m$}} 

\ncline{->}{b1}{a2} \nbput{{\scriptsize $1m$}} 
\ncline{->}{a2}{a4} \nbput[labelsep=3pt]{{\scriptsize $m$}}


{
\rput[c](36,26){\psset{unit=1mm,doubleline=true,arrowinset=0.6,arrowlength=0.5,arrowsize=0.5pt 2.1,nodesep=0pt,labelsep=1pt}
\pcline{->}(3,3)(0,0) \naput{{\scriptsize {\bfseries a}}}}}

{
\rput[c](24,16){\psset{unit=1mm,doubleline=true,arrowinset=0.6,arrowlength=0.5,arrowsize=0.5pt 2.1,nodesep=0pt,labelsep=2pt}
\pcline{->}(3,3)(0,0) \naput{{\scriptsize {\bfseries a}}}}}

\rput[c](16,33){{\scriptsize $\iso$}}
\rput[c](16,31){{\scriptsize\textsf{functoriality of $\otimes$}}}

\endpspicture\]
where we now write $(a,b)$ for $X(a,b)$.  This is the notation that we will use for the rest of this work.  

The unit axiom is as follows.

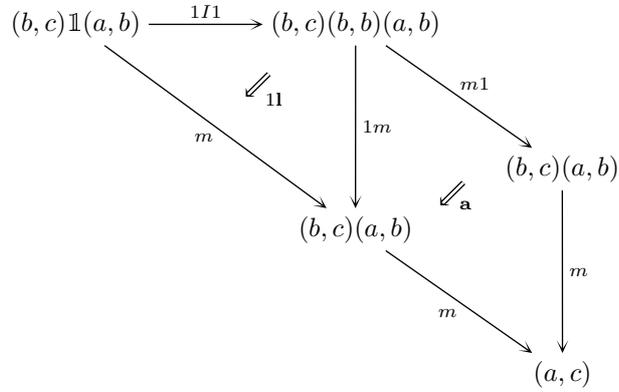
\[\psset{unit=0.16cm,labelsep=2pt,nodesep=3pt}
\pspicture(0,10)(50,45)


\rput(7,40){\rnode{a}{$(b,c)\1(a,b)$}}  
\rput(30,40){\rnode{b}{$(b,c)(b,b)(a,b)$}}  
\rput(47,28){\rnode{c}{$(b,c)(a,b)$}}  
\rput(30,23){\rnode{d}{$(b,c)(a,b)$}}  
\rput(47,11){\rnode{e}{$(a,c)$}}  

\ncline{->}{a}{b} \naput{{\scriptsize $1I1$}} 
\ncline{->}{b}{c} \naput{{\scriptsize $m1$}} 
\ncline{->}{c}{e} \naput{{\scriptsize $m$}} 
\ncline{->}{a}{d} \nbput{{\scriptsize $m$}} 
\ncline{->}{d}{e} \nbput{{\scriptsize $m$}} 

\ncline{->}{b}{d} \naput{{\scriptsize $1m$}} 

{
\rput[c](21,34){\psset{unit=1mm,doubleline=true,arrowinset=0.6,arrowlength=0.5,arrowsize=0.5pt 2.1,nodesep=0pt,labelsep=2pt}
\pcline{->}(3,3)(0,0) \naput{{\scriptsize 1{\bfseries l}}}}}

{
\rput[c](37,25){\psset{unit=1mm,doubleline=true,arrowinset=0.6,arrowlength=0.5,arrowsize=0.5pt 2.1,nodesep=0pt,labelsep=2pt}
\pcline{->}(3,3)(0,0) \naput{{\scriptsize {\bfseries a}}}}}

\endpspicture
\]

=

\[\psset{unit=0.16cm,labelsep=2pt,nodesep=3pt}
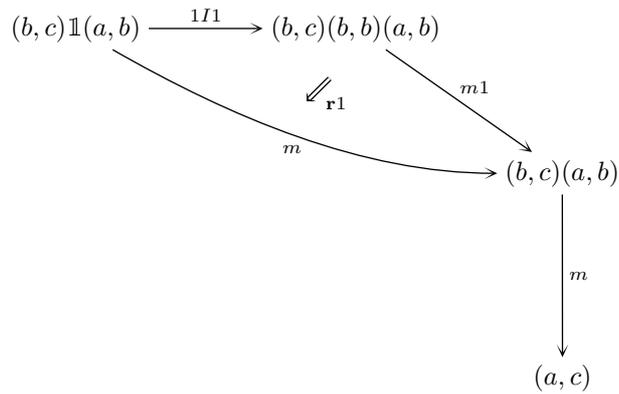
\pspicture(50,45)


\rput(7,40){\rnode{a}{$(b,c)\1(a,b)$}}  
\rput(30,40){\rnode{b}{$(b,c)(b,b)(a,b)$}}  
\rput(47,28){\rnode{c}{$(b,c)(a,b)$}}  
\rput(47,11){\rnode{e}{$(a,c)$}}  

\ncline{->}{a}{b} \naput{{\scriptsize $1I1$}} 
\ncline{->}{b}{c} \naput{{\scriptsize $m1$}} 
\ncline{->}{c}{e} \naput{{\scriptsize $m$}} 
\nccurve[angleA=-30,angleB=180]{->}{a}{c} \nbput{{\scriptsize $m$}} 

{
\rput[c](26,34){\psset{unit=1mm,doubleline=true,arrowinset=0.6,arrowlength=0.5,arrowsize=0.5pt 2.1,nodesep=0pt,labelsep=2pt}
\pcline{->}(3,3)(0,0) \naput{{\scriptsize {\bfseries r}1}}}}

\endpspicture
\]

\end{mydefinition}

\begin{mydefinition}
 A \demph{weak \cl{B}-functor} $F: X \tra X'$ is given by 

\begin{itemize}
 \item a function $F\: \ob X \tra \ob X'$,

\item 1-cells $F_{ab} \: X(a,b) \tra X'(Fa, Fb) \in \cB$, and

\item 2-cells

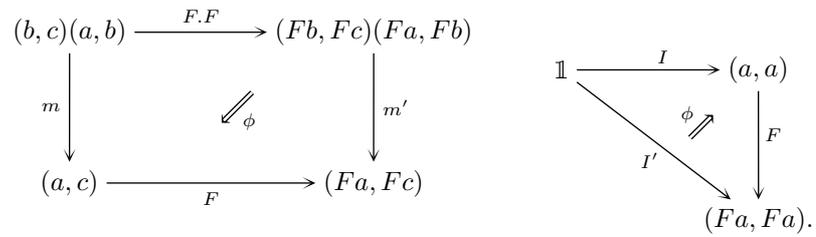
\[\psset{unit=0.1cm,labelsep=2pt,nodesep=3pt}
\pspicture(0,-5)(40,22)


\rput(0,20){\rnode{a1}{$(b,c)(a,b)$}} 
\rput(40,20){\rnode{a2}{$(Fb, Fc)(Fa,Fb)$}} 

\rput(0,0){\rnode{b1}{$(a,c)$}}   
\rput(40,0){\rnode{b2}{$(Fa,Fc)$}}  

\psset{nodesep=3pt,labelsep=3pt,arrows=->}
\ncline{a1}{a2}\naput{{\scriptsize $F.F$}} 
\ncline{b1}{b2}\nbput{{\scriptsize $F$}} 
\ncline{a1}{b1}\nbput{{\scriptsize $m$}} 
\ncline{a2}{b2}\naput{{\scriptsize $m'$}} 

{
\rput[c](20,8){\psset{unit=1mm,doubleline=true,arrowinset=0.6,arrowlength=0.5,arrowsize=0.5pt 2.1,nodesep=0pt,labelsep=1pt}
\pcline{->}(4,4)(0,0) \naput{{\scriptsize $\phi$}}}}


\endpspicture
\hh{7em}
\psset{unit=0.1cm,labelsep=2pt,nodesep=3pt,npos=0.4}
\pspicture(20,22)


\rput(0,20){\rnode{a1}{$\1$}}  
\rput(26,20){\rnode{a2}{$(a,a)$}}  
\rput(26,0){\rnode{a3}{$(Fa,Fa).$}}  

\ncline{->}{a1}{a2} \naput[npos=0.6]{{\scriptsize $I$}} 
\ncline{->}{a1}{a3} \nbput[npos=0.55]{{\scriptsize $I'$}} 
\ncline{->}{a2}{a3} \naput{{\scriptsize $F$}} 

{
\rput[c](17,11){\psset{unit=1mm,doubleline=true,arrowinset=0.6,arrowlength=0.5,arrowsize=0.5pt 2.1,nodesep=0pt,labelsep=1pt}
\pcline{->}(0,0)(3,3) \naput{{\scriptsize $\phi$}}}}

\endpspicture
\]

\end{itemize}

This data satisfies the following axioms:

\begin{itemize}
 \item associativity

\[\psset{unit=0.16cm,labelsep=2pt,nodesep=3pt}
\pspicture(0,10)(50,45)


\rput(20,12){\rnode{a2}{$(a,d)$}}  
\rput(46,12){\rnode{a4}{$(Fa,Fd)$}}  
\rput(0,23){\rnode{b1}{$(c,d)(a,c)$}}  
\rput(20,30){\rnode{c2}{$(b,d)(a,b)$}}  
\rput(46,30){\rnode{c4}{$(Fb,Fd)(Fa,Fb)$}}  
\rput(0,40){\rnode{d1}{$(c,d)(b,c)(a,b)$}}  
\rput(28,40){\rnode{d3}{$(Fc,Fd)(Fb,Fc)(Fa,Fb)$}}

\ncline{->}{d1}{d3} \naput{{\scriptsize $FFF$}} 
\ncline{->}{d3}{c4} \naput{{\scriptsize $m'1$}}
\ncline{->}{d1}{c2} \naput{{\scriptsize $m1$}}
\ncline{->}{c2}{c4} \naput{{\scriptsize $FF$}} 

\ncline{->}{d1}{b1} \nbput{{\scriptsize $1m$}}
\ncline{->}{b1}{a2} \nbput{{\scriptsize $m$}}
\ncline{->}{c2}{a2} \naput{{\scriptsize $m$}} 

\ncline{->}{c4}{a4} \naput{{\scriptsize $m'$}} 
\ncline{->}{a2}{a4} \nbput[labelsep=3pt]{{\scriptsize $F$}}


{
\rput[c](25,34){\psset{unit=1mm,doubleline=true,arrowinset=0.6,arrowlength=0.5,arrowsize=0.5pt 2.1,nodesep=0pt,labelsep=1pt}
\pcline{->}(3,3)(0,0) \naput{{\scriptsize $\phi. 1_F$}}}}

{
\rput[c](11,25){\psset{unit=1mm,doubleline=true,arrowinset=0.6,arrowlength=0.5,arrowsize=0.5pt 2.1,nodesep=0pt,labelsep=2pt}
\pcline{->}(3,3)(0,0) \naput{{\scriptsize {\bfseries a}}}}}

{
\rput[c](33,20){\psset{unit=1mm,doubleline=true,arrowinset=0.6,arrowlength=0.5,arrowsize=0.5pt 2.1,nodesep=0pt,labelsep=2pt}
\pcline{->}(3,3)(0,0) \naput{{\scriptsize $\phi$}}}}

\endpspicture\]

$=$

\[\psset{unit=0.16cm,labelsep=2pt,nodesep=3pt}
\pspicture(50,44)


\rput(20,12){\rnode{a2}{$(a,d)$}}  
\rput(46,12){\rnode{a4}{$(Fa,Fd)$}}  
\rput(0,23){\rnode{b1}{$(c,d)(a,c)$}}  

\rput(28,23){\rnode{b3}{$(Fc,Fd)(Fa,Fc)$}}  

\rput(46,30){\rnode{c4}{$(Fb,Fd)(Fa,Fb)$}}  
\rput(0,40){\rnode{d1}{$(c,d)(b,c)(a,b)$}}  
\rput(28,40){\rnode{d3}{$(Fc,Fd)(Fb,Fc)(Fa,Fb)$}}

\ncline{->}{d1}{d3} \naput{{\scriptsize $FFF$}} 
\ncline{->}{d3}{b3} \naput{{\scriptsize $1m'$}}
\ncline{->}{d1}{b1} \nbput{{\scriptsize $1m$}}
\ncline{->}{b1}{b3} \naput{{\scriptsize $FF$}}

\ncline{->}{d3}{c4} \naput{{\scriptsize $m'1$}} 
\ncline{->}{c4}{a4} \naput{{\scriptsize $m'$}}
\ncline{->}{b3}{a4} \naput{{\scriptsize $m'$}} 

\ncline{->}{b1}{a2} \nbput{{\scriptsize $m$}} 
\ncline{->}{a2}{a4} \nbput[labelsep=3pt]{{\scriptsize $F$}}


{
\rput[c](36,26){\psset{unit=1mm,doubleline=true,arrowinset=0.6,arrowlength=0.5,arrowsize=0.5pt 2.1,nodesep=0pt,labelsep=1pt}
\pcline{->}(3,3)(0,0) \naput{{\scriptsize $\mbox{{\bfseries a}}'$}}}}

{
\rput[c](24,16){\psset{unit=1mm,doubleline=true,arrowinset=0.6,arrowlength=0.5,arrowsize=0.5pt 2.1,nodesep=0pt,labelsep=2pt}
\pcline{->}(3,3)(0,0) \naput{{\scriptsize $\phi$}}}}

{
\rput[c](13,31){\psset{unit=1mm,doubleline=true,arrowinset=0.6,arrowlength=0.5,arrowsize=0.5pt 2.1,nodesep=0pt,labelsep=2pt}
\pcline{->}(3,3)(0,0) \naput{{\scriptsize $1_F.\phi$}}}}

\endpspicture\]

\newpage

\item right unit

\[\psset{unit=0.16cm,labelsep=2pt,nodesep=3pt}
\pspicture(0,8)(50,40)

%

\rput(20,12){\rnode{a2}{$(a,b)$}}  
\rput(46,12){\rnode{a4}{$(Fa,Fb)$}}  

\rput(20,30){\rnode{c2}{$(a,b)(a,a)$}}  
\rput(46,30){\rnode{c4}{$(Fa,Fb)(Fa,Fa)$}}  
\rput(0,40){\rnode{d1}{$(a,b)\1$}}  
\rput(28,40){\rnode{d3}{$(Fa,Fb)\1$}}

\ncline{->}{d1}{d3} \naput{{\scriptsize $F1$}} 
\ncline{->}{d3}{c4} \naput{{\scriptsize $1I'$}}
\ncline{->}{d1}{c2} \naput{{\scriptsize $1I$}}
\ncline{->}{c2}{c4} \naput{{\scriptsize $FF$}} 

\ncline{->}{d1}{a2} \nbput{{\scriptsize $l$}}
\ncline{->}{c2}{a2} \naput{{\scriptsize $m$}} 

\ncline{->}{c4}{a4} \naput{{\scriptsize $m'$}} 
\ncline{->}{a2}{a4} \nbput[labelsep=3pt]{{\scriptsize $F$}}


{
\rput[c](25,34){\psset{unit=1mm,doubleline=true,arrowinset=0.6,arrowlength=0.5,arrowsize=0.5pt 2.1,nodesep=0pt,labelsep=1pt}
\pcline{->}(3,3)(0,0) \naput{{\scriptsize $1_F.\phi$}}}}

{
\rput[c](12,26){\psset{unit=1mm,doubleline=true,arrowinset=0.6,arrowlength=0.5,arrowsize=0.5pt 2.1,nodesep=0pt,labelsep=2pt}
\pcline{->}(3,3)(0,0) \naput{{\scriptsize {\bfseries r}}}}}

{
\rput[c](33,20){\psset{unit=1mm,doubleline=true,arrowinset=0.6,arrowlength=0.5,arrowsize=0.5pt 2.1,nodesep=0pt,labelsep=2pt}
\pcline{->}(3,3)(0,0) \naput{{\scriptsize $\phi$}}}}

\endpspicture\]

$=$

\[\psset{unit=0.16cm,labelsep=2pt,nodesep=3pt}
\pspicture(50,44)

%

\rput(20,12){\rnode{a2}{$(a,b)$}}  
\rput(46,12){\rnode{a4}{$(Fa,Fb)$}}  

\rput(46,30){\rnode{c4}{$(Fa,Fb)(Fa,Fa)$}}  
\rput(0,40){\rnode{d1}{$(a,b)$}}  
\rput(28,40){\rnode{d3}{$(Fa,Fb)\1$}}  

\ncline{->}{d1}{d3} \naput{{\scriptsize $F1$}} 
\ncline{->}{d1}{a2} \nbput{{\scriptsize $r$}}

\ncline{->}{d3}{c4} \naput{{\scriptsize $1I'$}} 
\ncline{->}{c4}{a4} \naput{{\scriptsize $m'$}}
\ncline{->}{d3}{a4} \nbput{{\scriptsize $r$}} 

\ncline{->}{a2}{a4} \nbput[labelsep=3pt]{{\scriptsize $F$}}


{
\rput[c](39,25){\psset{unit=1mm,doubleline=true,arrowinset=0.6,arrowlength=0.5,arrowsize=0.5pt 2.1,nodesep=0pt,labelsep=1pt}
\pcline{->}(3,3)(0,0) \naput{{\scriptsize $\mbox{{\bfseries r}}'$}}}}

\rput[c](23,27){{\scriptsize $\iso$}}
\rput[c](23,25){{\scriptsize\textsf{naturality of $r$} }}

\endpspicture\]

\newpage

\item left unit

\[\psset{unit=0.16cm,labelsep=2pt,nodesep=3pt}
\pspicture(0,8)(50,40)

%

\rput(20,12){\rnode{a2}{$(a,b)$}}  
\rput(46,12){\rnode{a4}{$(Fa,Fb)$}}  

\rput(20,30){\rnode{c2}{$(b,b)(a,b)$}}  
\rput(46,30){\rnode{c4}{$(Fb,Fb)(Fa,Fb)$}}  
\rput(0,40){\rnode{d1}{$\1(a,b)$}}  
\rput(28,40){\rnode{d3}{$\1(Fa,Fb)$}}

\ncline{->}{d1}{d3} \naput{{\scriptsize $F1$}} 
\ncline{->}{d3}{c4} \naput{{\scriptsize $I'1$}}
\ncline{->}{d1}{c2} \naput{{\scriptsize $I1$}}
\ncline{->}{c2}{c4} \naput{{\scriptsize $FF$}} 

\ncline{->}{d1}{a2} \nbput{{\scriptsize $l$}}
\ncline{->}{c2}{a2} \naput{{\scriptsize $m$}} 

\ncline{->}{c4}{a4} \naput{{\scriptsize $m'$}} 
\ncline{->}{a2}{a4} \nbput[labelsep=3pt]{{\scriptsize $F$}}


{
\rput[c](25,34){\psset{unit=1mm,doubleline=true,arrowinset=0.6,arrowlength=0.5,arrowsize=0.5pt 2.1,nodesep=0pt,labelsep=1pt}
\pcline{->}(3,3)(0,0) \naput{{\scriptsize $\phi.1_F$}}}}

{
\rput[c](12,26){\psset{unit=1mm,doubleline=true,arrowinset=0.6,arrowlength=0.5,arrowsize=0.5pt 2.1,nodesep=0pt,labelsep=2pt}
\pcline{->}(3,3)(0,0) \naput{{\scriptsize {\bfseries l}}}}}

{
\rput[c](33,20){\psset{unit=1mm,doubleline=true,arrowinset=0.6,arrowlength=0.5,arrowsize=0.5pt 2.1,nodesep=0pt,labelsep=2pt}
\pcline{->}(3,3)(0,0) \naput{{\scriptsize $\phi$}}}}

\endpspicture\]

$=$

\[\psset{unit=0.16cm,labelsep=2pt,nodesep=3pt}
\pspicture(50,44)

%

\rput(20,12){\rnode{a2}{$(a,b)$}}  
\rput(46,12){\rnode{a4}{$(Fa,Fb)$}}  

\rput(46,30){\rnode{c4}{$(Fb,Fb)(Fa,Fb)$}}  
\rput(0,40){\rnode{d1}{$(a,b)$}}  
\rput(28,40){\rnode{d3}{$\1(Fa,Fb)$}}  

\ncline{->}{d1}{d3} \naput{{\scriptsize $1F$}} 
\ncline{->}{d1}{a2} \nbput{{\scriptsize $l$}}

\ncline{->}{d3}{c4} \naput{{\scriptsize $I'1$}} 
\ncline{->}{c4}{a4} \naput{{\scriptsize $m'$}}
\ncline{->}{d3}{a4} \nbput{{\scriptsize $l$}} 

\ncline{->}{a2}{a4} \nbput[labelsep=3pt]{{\scriptsize $F$}}


{
\rput[c](39,25){\psset{unit=1mm,doubleline=true,arrowinset=0.6,arrowlength=0.5,arrowsize=0.5pt 2.1,nodesep=0pt,labelsep=1pt}
\pcline{->}(3,3)(0,0) \naput{{\scriptsize $\mbox{{\bfseries l}}'$}}}}

\rput[c](23,27){{\scriptsize $\iso$}}
\rput[c](23,25){{\scriptsize\textsf{naturality of $l$} }}

\endpspicture\]

\end{itemize}

\end{mydefinition}


\begin{definition} \label{functorcomp}
 Composition of weak \cl{B}-functors is defined as follows.  Given composable weak $\cB$-functors as below
\[X \ltmap{(F, \phi^F)} Y \ltmap{(G, \phi^G)} Z\]
we define the composite $GF$.

\begin{itemize}
 \item The action on objects is given by the composite
\[\ob X \ltmap{F} \ob Y \ltmap{G} \ob Z.\]

\item The action on hom-categories is given by the composite
\[X(a,b) \ltmap{F_{ab}} Y(Fa,Fb) \ltmap{G_{Fa,Fb}} Z(GFa, GFb).\]

\item The composition constraint $\phi^{GF}$ is given by the following composite 2-cell.

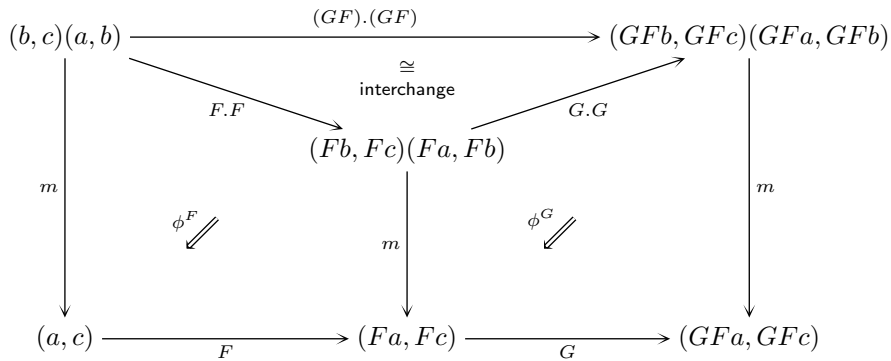
\[\psset{unit=0.1cm,labelsep=0pt,nodesep=3pt}
\pspicture(0,-10)(90,45)


\rput(0,40){\rnode{a1}{$(b,c)(a,b)$}} 
\rput(45,25){\rnode{a2}{$(Fb,Fc)(Fa,Fb)$}} 
\rput(90,40){\rnode{a3}{$(GFb,GFc)(GFa,GFb)$}} 

\rput(0,0){\rnode{b1}{$(a,c)$}}   
\rput(45,0){\rnode{b2}{$(Fa,Fc)$}}  
\rput(90,0){\rnode{b3}{$(GFa,GFc)$}}  

\psset{nodesep=3pt,labelsep=2pt,arrows=->}
\ncline{a1}{a3}\naput[labelsep=3pt]{{\scriptsize $(GF).(GF)$}} 
\ncline{a1}{a2}\nbput{{\scriptsize $F.F$}} %
\ncline{a2}{a3}\nbput{{\scriptsize $G.G$}} %

\ncline{b1}{b2}\nbput{{\scriptsize $F$}} 
\ncline{b2}{b3}\nbput{{\scriptsize $G$}} 

\ncline{a1}{b1}\nbput{{\scriptsize $m$}} 
\ncline{a2}{b2}\nbput[npos=0.5]{{\scriptsize $m$}} 
\ncline{a3}{b3}\naput{{\scriptsize $m$}} 

{\psset{doubleline=true,arrowinset=0.6,arrowlength=0.5,arrowsize=0.5pt 2.1,nodesep=0pt}
\rput[c](16,12){\pcline{->}(4,4)(0,0) \nbput[labelsep=-1pt]{{\scriptsize $\phi^F$}}}}

{\psset{doubleline=true,arrowinset=0.6,arrowlength=0.5,arrowsize=0.5pt 2.1,nodesep=0pt}
\rput[c](63,12){\pcline{->}(4,4)(0,0) \nbput[labelsep=0pt]{{\scriptsize $\phi^G$}}}}

\rput[c](45,36){{\scriptsize $\iso$}}
\rput[c](45,33){{\scriptsize\textsf{interchange}}}

\endpspicture\]

\item The unit constraint is given by the following composite 2-cell.

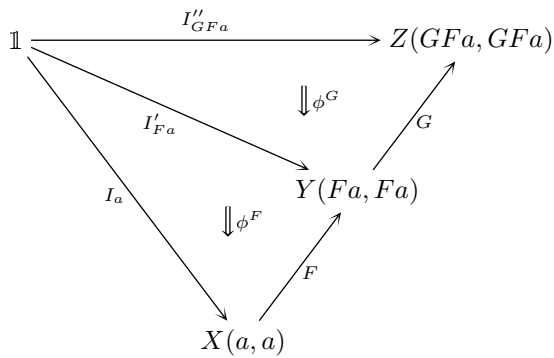
\[\psset{unit=0.1cm,labelsep=0pt,nodesep=3pt}
\pspicture(0,-10)(90,45)


\rput(0,40){\rnode{a1}{$\1$}} 

\rput(30,0){\rnode{b1}{$X(a,a)$}}   
\rput(45,20){\rnode{a2}{$Y(Fa,Fa)$}} 
\rput(60,40){\rnode{a3}{$Z(GFa,GFa)$}} 

\psset{nodesep=3pt,labelsep=0pt,arrows=->}
\ncline{a1}{a3}\naput[labelsep=3pt]{{\scriptsize $I''_{GFa}$}} 
\ncline{a1}{a2}\nbput{{\scriptsize $I'_{Fa}$}} %
\ncline{a1}{b1}\nbput{{\scriptsize $I_a$}} 

\ncline{a2}{a3}\nbput{{\scriptsize $G$}} %
\ncline{b1}{a2}\nbput{{\scriptsize $F$}} 

{\psset{doubleline=true,arrowinset=0.6,arrowlength=0.5,arrowsize=0.5pt 2.1,nodesep=0pt}
\rput[c](28,14){\pcline{->}(0,4)(0,0) \naput[labelsep=1pt]{{\scriptsize $\phi^F$}}}}

{\psset{doubleline=true,arrowinset=0.6,arrowlength=0.5,arrowsize=0.5pt 2.1,nodesep=0pt}
\rput[c](38,30){\pcline{->}(0,4)(0,0) \naput[labelsep=1pt]{{\scriptsize $\phi^G$}}}}

\endpspicture\]

\end{itemize}

We immediately check that this satisfies the axioms as follows.  

\begin{itemize}

\item Associativity:

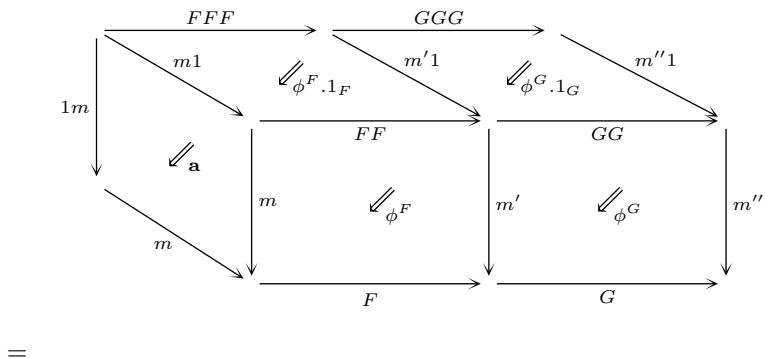
\[\psset{unit=0.12cm,labelsep=2pt,nodesep=3pt}
\pspicture(0,10)(80,45)



\rput(20,12){\rnode{a2}{$$}}  
\rput(46,12){\rnode{a4}{$$}}  
\rput(3,23){\rnode{b1}{$$}}  
\rput(20,30){\rnode{c2}{$$}}  
\rput(46,30){\rnode{c4}{$$}}  
\rput(3,40){\rnode{d1}{$$}}  
\rput(28,40){\rnode{d3}{$$}}  

\rput(53,40){\rnode{d5}{$$}}  
\rput(72,30){\rnode{c6}{$$}}  
\rput(72,12){\rnode{a6}{$$}}

\ncline{->}{d1}{d3} \naput{{\scriptsize $FFF$}} 
\ncline{->}{d3}{c4} \naput{{\scriptsize $m'1$}}
\ncline{->}{d1}{c2} \naput{{\scriptsize $m1$}}
\ncline{->}{c2}{c4} \nbput{{\scriptsize $FF$}} 

\ncline{->}{d1}{b1} \nbput{{\scriptsize $1m$}}
\ncline{->}{b1}{a2} \nbput{{\scriptsize $m$}}
\ncline{->}{c2}{a2} \naput{{\scriptsize $m$}} 

\ncline{->}{c4}{a4} \naput{{\scriptsize $m'$}} 
\ncline{->}{a2}{a4} \nbput[labelsep=3pt]{{\scriptsize $F$}}

\ncline{->}{d3}{d5} \naput{{\scriptsize $GGG$}} 
\ncline{->}{c4}{c6} \nbput{{\scriptsize $GG$}} 
\ncline{->}{a4}{a6} \nbput{{\scriptsize $G$}} 
\ncline{->}{d5}{c6} \naput{{\scriptsize $m''1$}} 
\ncline{->}{c6}{a6} \naput{{\scriptsize $m''$}}


{
\rput[c](23,34){\psset{unit=1mm,doubleline=true,arrowinset=0.6,arrowlength=0.5,arrowsize=0.5pt 2.1,nodesep=0pt,labelsep=-1pt}
\pcline{->}(3,3)(0,0) \naput{{\scriptsize $\phi^F. 1_F$}}}}

{
\rput[c](48,34){\psset{unit=1mm,doubleline=true,arrowinset=0.6,arrowlength=0.5,arrowsize=0.5pt 2.1,nodesep=0pt,labelsep=-1pt}
\pcline{->}(3,3)(0,0) \naput{{\scriptsize $\phi^G. 1_G$}}}}

{
\rput[c](11,25){\psset{unit=1mm,doubleline=true,arrowinset=0.6,arrowlength=0.5,arrowsize=0.5pt 2.1,nodesep=0pt,labelsep=2pt}
\pcline{->}(3,3)(0,0) \naput{{\scriptsize {\bfseries a}}}}}

{
\rput[c](33,20){\psset{unit=1mm,doubleline=true,arrowinset=0.6,arrowlength=0.5,arrowsize=0.5pt 2.1,nodesep=0pt,labelsep=0pt}
\pcline{->}(3,3)(0,0) \naput{{\scriptsize $\phi^F$}}}}

{
\rput[c](58,20){\psset{unit=1mm,doubleline=true,arrowinset=0.6,arrowlength=0.5,arrowsize=0.5pt 2.1,nodesep=0pt,labelsep=0pt}
\pcline{->}(3,3)(0,0) \naput{{\scriptsize $\phi^G$}}}}

\endpspicture\]

$=$

\[\psset{unit=0.12cm,labelsep=2pt,nodesep=3pt}
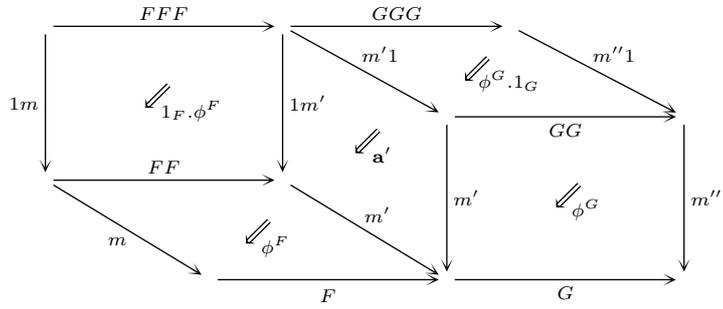
\pspicture(0,10)(80,44)


%
%
%

\rput(20,12){\rnode{a2}{$$}}  
\rput(46,12){\rnode{a4}{$$}}  
\rput(2,23){\rnode{b1}{$$}}  

\rput(28,23){\rnode{b3}{$$}}  

\rput(46,30){\rnode{c4}{$$}}  
\rput(2,40){\rnode{d1}{$$}}  
\rput(28,40){\rnode{d3}{$$}}

\ncline{->}{d1}{d3} \naput{{\scriptsize $FFF$}} 
\ncline{->}{d3}{b3} \naput{{\scriptsize $1m'$}}
\ncline{->}{d1}{b1} \nbput{{\scriptsize $1m$}}
\ncline{->}{b1}{b3} \naput{{\scriptsize $FF$}} 

\rput(53,40){\rnode{d5}{$$}}  
\rput(72,30){\rnode{c6}{$$}}  
\rput(72,12){\rnode{a6}{$$}}  

\ncline{->}{d3}{d5} \naput{{\scriptsize $GGG$}} 
\ncline{->}{c4}{c6} \nbput{{\scriptsize $GG$}} 
\ncline{->}{a4}{a6} \nbput{{\scriptsize $G$}} 
\ncline{->}{d5}{c6} \naput{{\scriptsize $m''1$}} 
\ncline{->}{c6}{a6} \naput{{\scriptsize $m''$}}

\ncline{->}{d3}{c4} \naput{{\scriptsize $m'1$}} 
\ncline{->}{c4}{a4} \naput{{\scriptsize $m'$}}
\ncline{->}{b3}{a4} \naput{{\scriptsize $m'$}} 

\ncline{->}{b1}{a2} \nbput{{\scriptsize $m$}} 
\ncline{->}{a2}{a4} \nbput[labelsep=3pt]{{\scriptsize $F$}}


{
\rput[c](36,26){\psset{unit=1mm,doubleline=true,arrowinset=0.6,arrowlength=0.5,arrowsize=0.5pt 2.1,nodesep=0pt,labelsep=1pt}
\pcline{->}(3,3)(0,0) \naput{{\scriptsize $\mbox{{\bfseries a}}'$}}}}

{
\rput[c](24,16){\psset{unit=1mm,doubleline=true,arrowinset=0.6,arrowlength=0.5,arrowsize=0.5pt 2.1,nodesep=0pt,labelsep=0pt}
\pcline{->}(3,3)(0,0) \naput{{\scriptsize $\phi^F$}}}}

{
\rput[c](13,31){\psset{unit=1mm,doubleline=true,arrowinset=0.6,arrowlength=0.5,arrowsize=0.5pt 2.1,nodesep=0pt,labelsep=0pt}
\pcline{->}(3,3)(0,0) \naput{{\scriptsize $1_F.\phi^F$}}}}

{
\rput[c](48,34){\psset{unit=1mm,doubleline=true,arrowinset=0.6,arrowlength=0.5,arrowsize=0.5pt 2.1,nodesep=0pt,labelsep=-1pt}
\pcline{->}(3,3)(0,0) \naput{{\scriptsize $\phi^G. 1_G$}}}}

{
\rput[c](58,20){\psset{unit=1mm,doubleline=true,arrowinset=0.6,arrowlength=0.5,arrowsize=0.5pt 2.1,nodesep=0pt,labelsep=0pt}
\pcline{->}(3,3)(0,0) \naput{{\scriptsize $\phi^G$}}}}

\endpspicture\]

$=$

\[\psset{unit=0.12cm,labelsep=2pt,nodesep=3pt}
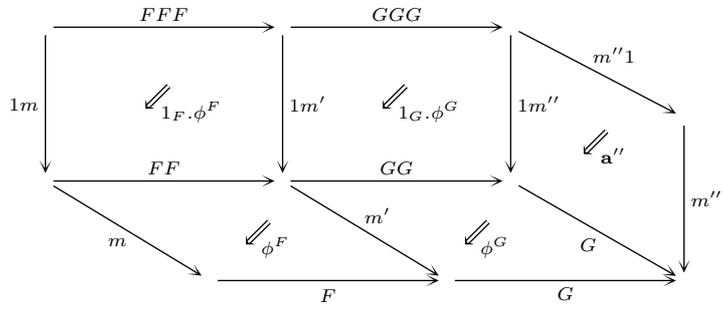
\pspicture(0,10)(80,44)


%
%
%

\rput(20,12){\rnode{a2}{$$}}  
\rput(46,12){\rnode{a4}{$$}}  
\rput(2,23){\rnode{b1}{$$}}  

\rput(28,23){\rnode{b3}{$$}}  

\rput(53,23){\rnode{b5}{$$}}  

\rput(2,40){\rnode{d1}{$$}}  
\rput(28,40){\rnode{d3}{$$}}

\ncline{->}{d1}{d3} \naput{{\scriptsize $FFF$}} 
\ncline{->}{d3}{b3} \naput{{\scriptsize $1m'$}}
\ncline{->}{d1}{b1} \nbput{{\scriptsize $1m$}}
\ncline{->}{b1}{b3} \naput{{\scriptsize $FF$}} 

\rput(53,40){\rnode{d5}{$$}}  
\rput(72,30){\rnode{c6}{$$}}  
\rput(72,12){\rnode{a6}{$$}}  

\ncline{->}{d3}{d5} \naput{{\scriptsize $GGG$}} 
\ncline{->}{a4}{a6} \nbput{{\scriptsize $G$}} 
\ncline{->}{d5}{c6} \naput{{\scriptsize $m''1$}} 
\ncline{->}{c6}{a6} \naput{{\scriptsize $m''$}} 

\ncline{->}{b3}{b5} \naput{{\scriptsize $GG$}} 
\ncline{->}{d5}{b5} \naput{{\scriptsize $1m''$}} 
\ncline{->}{b5}{a6} \nbput{{\scriptsize $G$}}

\ncline{->}{b3}{a4} \naput{{\scriptsize $m'$}} 

\ncline{->}{b1}{a2} \nbput{{\scriptsize $m$}} 
\ncline{->}{a2}{a4} \nbput[labelsep=3pt]{{\scriptsize $F$}}


{
\rput[c](61,26){\psset{unit=1mm,doubleline=true,arrowinset=0.6,arrowlength=0.5,arrowsize=0.5pt 2.1,nodesep=0pt,labelsep=1pt}
\pcline{->}(3,3)(0,0) \naput{{\scriptsize $\mbox{{\bfseries a}}''$}}}}

{
\rput[c](24,16){\psset{unit=1mm,doubleline=true,arrowinset=0.6,arrowlength=0.5,arrowsize=0.5pt 2.1,nodesep=0pt,labelsep=0pt}
\pcline{->}(3,3)(0,0) \naput{{\scriptsize $\phi^F$}}}}

{
\rput[c](13,31){\psset{unit=1mm,doubleline=true,arrowinset=0.6,arrowlength=0.5,arrowsize=0.5pt 2.1,nodesep=0pt,labelsep=0pt}
\pcline{->}(3,3)(0,0) \naput{{\scriptsize $1_F.\phi^F$}}}}

{
\rput[c](39,31){\psset{unit=1mm,doubleline=true,arrowinset=0.6,arrowlength=0.5,arrowsize=0.5pt 2.1,nodesep=0pt,labelsep=0pt}
\pcline{->}(3,3)(0,0) \naput{{\scriptsize $1_G.\phi^G$}}}}

{
\rput[c](48,16){\psset{unit=1mm,doubleline=true,arrowinset=0.6,arrowlength=0.5,arrowsize=0.5pt 2.1,nodesep=0pt,labelsep=0pt}
\pcline{->}(3,3)(0,0) \naput{{\scriptsize $\phi^G$}}}}

\endpspicture\]

\item Right unit (and similarly left unit):

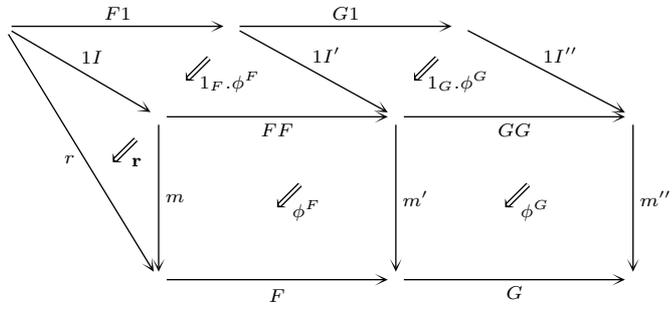
\[\psset{unit=0.12cm,labelsep=2pt,nodesep=3pt}
\pspicture(0,10)(80,45)

%


\rput(20,12){\rnode{a2}{$$}}  
\rput(46,12){\rnode{a4}{$$}}  
\rput(20,30){\rnode{c2}{$$}}  
\rput(46,30){\rnode{c4}{$$}}  
\rput(3,40){\rnode{d1}{$$}}  
\rput(28,40){\rnode{d3}{$$}}  

\rput(53,40){\rnode{d5}{$$}}  
\rput(72,30){\rnode{c6}{$$}}  
\rput(72,12){\rnode{a6}{$$}}

\ncline{->}{d1}{d3} \naput{{\scriptsize $F1$}} 
\ncline{->}{d3}{c4} \naput{{\scriptsize $1I'$}}
\ncline{->}{d1}{c2} \naput{{\scriptsize $1I$}}
\ncline{->}{c2}{c4} \nbput{{\scriptsize $FF$}} 

\ncline{->}{d1}{a2} \nbput{{\scriptsize $r$}}
\ncline{->}{c2}{a2} \naput{{\scriptsize $m$}} 

\ncline{->}{c4}{a4} \naput{{\scriptsize $m'$}} 
\ncline{->}{a2}{a4} \nbput[labelsep=3pt]{{\scriptsize $F$}}

\ncline{->}{d3}{d5} \naput{{\scriptsize $G1$}} 
\ncline{->}{c4}{c6} \nbput{{\scriptsize $GG$}} 
\ncline{->}{a4}{a6} \nbput{{\scriptsize $G$}} 
\ncline{->}{d5}{c6} \naput{{\scriptsize $1I''$}} 
\ncline{->}{c6}{a6} \naput{{\scriptsize $m''$}}


{
\rput[c](23,34){\psset{unit=1mm,doubleline=true,arrowinset=0.6,arrowlength=0.5,arrowsize=0.5pt 2.1,nodesep=0pt,labelsep=-1pt}
\pcline{->}(3,3)(0,0) \naput{{\scriptsize $1_F.\phi^F$}}}}

{
\rput[c](48,34){\psset{unit=1mm,doubleline=true,arrowinset=0.6,arrowlength=0.5,arrowsize=0.5pt 2.1,nodesep=0pt,labelsep=-1pt}
\pcline{->}(3,3)(0,0) \naput{{\scriptsize $1_G.\phi^G$}}}}

{
\rput[c](15,25){\psset{unit=1mm,doubleline=true,arrowinset=0.6,arrowlength=0.5,arrowsize=0.5pt 2.1,nodesep=0pt,labelsep=2pt}
\pcline{->}(3,3)(0,0) \naput{{\scriptsize {\bfseries r}}}}}

{
\rput[c](33,20){\psset{unit=1mm,doubleline=true,arrowinset=0.6,arrowlength=0.5,arrowsize=0.5pt 2.1,nodesep=0pt,labelsep=0pt}
\pcline{->}(3,3)(0,0) \naput{{\scriptsize $\phi^F$}}}}

{
\rput[c](58,20){\psset{unit=1mm,doubleline=true,arrowinset=0.6,arrowlength=0.5,arrowsize=0.5pt 2.1,nodesep=0pt,labelsep=0pt}
\pcline{->}(3,3)(0,0) \naput{{\scriptsize $\phi^G$}}}}

\endpspicture\]

$=$

\[\psset{unit=0.12cm,labelsep=2pt,nodesep=3pt}
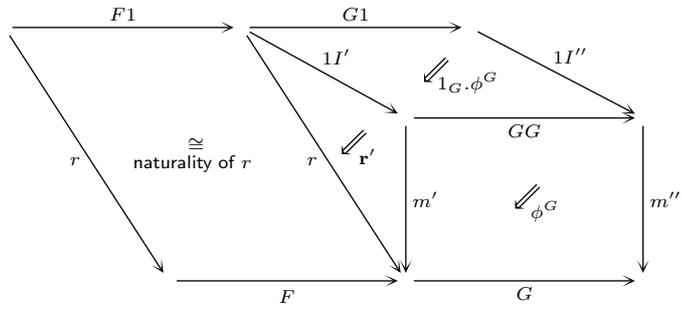
\pspicture(0,10)(80,44)

%

%
%
%

\rput(20,12){\rnode{a2}{$$}}  
\rput(46,12){\rnode{a4}{$$}}  


\rput(46,30){\rnode{c4}{$$}}  
\rput(2,40){\rnode{d1}{$$}}  
\rput(28,40){\rnode{d3}{$$}}

\ncline{->}{d1}{d3} \naput{{\scriptsize $F1$}} 
\ncline{->}{d1}{a2} \nbput{{\scriptsize $r$}}

\rput(53,40){\rnode{d5}{$$}}  
\rput(72,30){\rnode{c6}{$$}}  
\rput(72,12){\rnode{a6}{$$}}  

\ncline{->}{d3}{d5} \naput{{\scriptsize $G1$}} 
\ncline{->}{c4}{c6} \nbput{{\scriptsize $GG$}} 
\ncline{->}{a4}{a6} \nbput{{\scriptsize $G$}} 
\ncline{->}{d5}{c6} \naput{{\scriptsize $1I''$}} 
\ncline{->}{c6}{a6} \naput{{\scriptsize $m''$}}

\ncline{->}{d3}{c4} \naput{{\scriptsize $1I'$}} 
\ncline{->}{c4}{a4} \naput{{\scriptsize $m'$}}
\ncline{->}{d3}{a4} \nbput{{\scriptsize $r$}} 

\ncline{->}{a2}{a4} \nbput[labelsep=3pt]{{\scriptsize $F$}}


{
\rput[c](39,26){\psset{unit=1mm,doubleline=true,arrowinset=0.6,arrowlength=0.5,arrowsize=0.5pt 2.1,nodesep=0pt,labelsep=1pt}
\pcline{->}(3,3)(0,0) \naput{{\scriptsize $\mbox{{\bfseries r}}'$}}}}



\rput[c](23,27){{\scriptsize $\iso$}}
\rput[c](23,25){{\scriptsize\textsf{naturality of $r$} }}

{
\rput[c](48,34){\psset{unit=1mm,doubleline=true,arrowinset=0.6,arrowlength=0.5,arrowsize=0.5pt 2.1,nodesep=0pt,labelsep=-1pt}
\pcline{->}(3,3)(0,0) \naput{{\scriptsize $ 1_G.\phi^G$}}}}

{
\rput[c](58,20){\psset{unit=1mm,doubleline=true,arrowinset=0.6,arrowlength=0.5,arrowsize=0.5pt 2.1,nodesep=0pt,labelsep=0pt}
\pcline{->}(3,3)(0,0) \naput{{\scriptsize $\phi^G$}}}}

\endpspicture\]

$=$

\[\psset{unit=0.12cm,labelsep=2pt,nodesep=3pt}
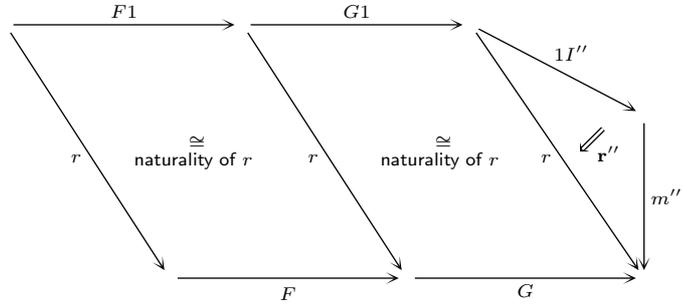
\pspicture(0,10)(80,44)

%

%
%
%

\rput(20,12){\rnode{a2}{$$}}  
\rput(46,12){\rnode{a4}{$$}}  


\rput(2,40){\rnode{d1}{$$}}  
\rput(28,40){\rnode{d3}{$$}}

\ncline{->}{d1}{d3} \naput{{\scriptsize $F1$}} 
\ncline{->}{d1}{a2} \nbput{{\scriptsize $r$}}

\rput(53,40){\rnode{d5}{$$}}  
\rput(72,30){\rnode{c6}{$$}}  
\rput(72,12){\rnode{a6}{$$}}  

\ncline{->}{d3}{d5} \naput{{\scriptsize $G1$}} 
\ncline{->}{a4}{a6} \nbput{{\scriptsize $G$}} 
\ncline{->}{d5}{c6} \naput{{\scriptsize $1I''$}} 
\ncline{->}{c6}{a6} \naput{{\scriptsize $m''$}}

\ncline{->}{d5}{a6} \nbput{{\scriptsize $r$}} 
\ncline{->}{d3}{a4} \nbput{{\scriptsize $r$}} 

\ncline{->}{a2}{a4} \nbput[labelsep=3pt]{{\scriptsize $F$}}


{
\rput[c](65,26){\psset{unit=1mm,doubleline=true,arrowinset=0.6,arrowlength=0.5,arrowsize=0.5pt 2.1,nodesep=0pt,labelsep=1pt}
\pcline{->}(3,3)(0,0) \naput{{\scriptsize $\mbox{{\bfseries r}}''$}}}}



\rput[c](23,27){{\scriptsize $\iso$}}
\rput[c](23,25){{\scriptsize\textsf{naturality of $r$} }}

\rput[c](50,27){{\scriptsize $\iso$}}
\rput[c](50,25){{\scriptsize\textsf{naturality of $r$} }}

%
%

\endpspicture\]

\end{itemize}

\end{definition}


\begin{mydefinition}

Let $F,G \: X \tra Y$ be weak $\cB$-functors such that $Fa=Ga$ for all objects $a \in X$.  A \demph{$\cB$-icon}
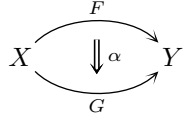
\[
\psset{unit=0.1cm,labelsep=2pt,nodesep=2pt}
\pspicture(20,20)

\rput(0,10){\rnode{a1}{$X$}}  
\rput(20,10){\rnode{a2}{$Y$}}  

\ncarc[arcangle=45]{->}{a1}{a2}\naput{{\scriptsize $F$}}
\ncarc[arcangle=-45]{->}{a1}{a2}\nbput{{\scriptsize $G$}}

\pcline[linewidth=0.6pt,doubleline=true,arrowinset=0.6,arrowlength=0.8,arrowsize=0.5pt 2.1]{->}(10,13)(10,7)  \naput{{\scriptsize $\alpha$}}


\endpspicture
\]
is given by, for all pairs of objects $a,b \in X$ a 2-cell

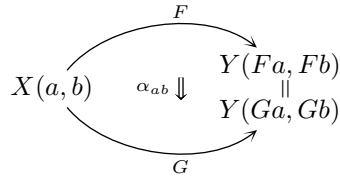
\[
\psset{unit=0.1cm,labelsep=2pt,nodesep=2pt}
\pspicture(0,-4)(30,24)

\rput(0,10){\rnode{a1}{$X(a,b)$}}  
\rput(30,13){\rnode{a2}{$Y(Fa,Fb)$}}  
\rput(31,10){\rotatebox{90}{$=$}}
\rput(30,7){\rnode{a3}{$Y(Ga,Gb)$}}  

\ncarc[arcangle=45]{->}{a1}{a2}\naput[npos=0.6]{{\scriptsize $F$}}
\ncarc[arcangle=-45]{->}{a1}{a3}\nbput[npos=0.6]{{\scriptsize $G$}}

{
\rput[c](17,8.5){\psset{unit=1mm,doubleline=true,arrowinset=0.6,arrowlength=0.5,arrowsize=0.5pt 2.1,nodesep=0pt,labelsep=2pt}
\pcline{->}(0,3)(0,0) \nbput{{\scriptsize $\alpha_{ab}$}}}}


\endpspicture
\]
satisfying the following axioms.

\begin{itemize}
 \item Composition:

\[
\psset{unit=0.1cm,labelsep=2pt,nodesep=2pt}
\pspicture(0,-15)(80,40)

%

\rput(0,30){\rnode{a1}{$(b,c)(a,b)$}}  
\rput(30,33){\rnode{a2}{\makebox[4em][l]{$(Fb,Fc)(Fa,Fb)$}}}  
\rput(38,30){\rotatebox{90}{$=$}}
\rput(30,27){\rnode{a3}{\makebox[4em][l]{$(Gb,Gc)(Ga,Gb)$}}}  

\rput(0,3){\rnode{b1}{$(a,c)$}}  
\rput(30,3){\rnode{b2}{\makebox[3em][l]{$(Ga,Gc)$}}}

\ncarc[arcangle=45]{->}{a1}{a2}\naput[npos=0.56]{{\scriptsize $FF$}}
\ncarc[arcangle=-45]{->}{a1}{a3}\nbput[npos=0.56]{{\scriptsize $GG$}}

\ncline{->}{a1}{b1} \nbput{{\scriptsize $m$}}
\ncline{->}{a3}{b2} \naput{{\scriptsize $m'$}}

\ncarc[arcangle=-45]{->}{b1}{b2}\nbput[npos=0.5]{{\scriptsize $G$}}

{
\rput[c](17,28.5){\psset{unit=1mm,doubleline=true,arrowinset=0.6,arrowlength=0.5,arrowsize=0.5pt 2.1,nodesep=0pt,labelsep=2pt}
\pcline{->}(0,3)(0,0) \nbput{{\scriptsize $\alpha\alpha$}}}}

{
\rput[c](17,8.5){\psset{unit=1mm,doubleline=true,arrowinset=0.6,arrowlength=0.5,arrowsize=0.5pt 2.1,nodesep=0pt,labelsep=2pt}
\pcline{->}(0,3)(0,0) \nbput{{\scriptsize $\phi^G$}}}}


\rput(60,20){$=$}

\endpspicture
\psset{unit=0.1cm,labelsep=2pt,nodesep=2pt}
\pspicture(0,-15)(50,40)

%
%

\rput(0,30){\rnode{a1}{$(b,c)(a,b)$}}  
\rput(30,30){\rnode{a2}{\makebox[4em][l]{$(Fb,Fc)(Fa,Fb)$}}}  

\rput(0,3){\rnode{b1}{$(a,c)$}}  
\rput(30,6){\rnode{b2}{\makebox[3em][l]{$(Fa,Fc)$}}}  

\rput(32,3){\rotatebox{90}{$=$}}
\rput(30,0){\rnode{b3}{\makebox[3em][l]{$(Ga,Gc)$}}}  

\ncarc[arcangle=45]{->}{a1}{a2}\naput[npos=0.56]{{\scriptsize $FF$}}

\ncline{->}{a1}{b1} \nbput{{\scriptsize $m$}}
\ncline{->}{a2}{b2} \naput{{\scriptsize $m'$}}

\ncarc[arcangle=45]{->}{b1}{b2}\naput[npos=0.56]{{\scriptsize $F$}}
\ncarc[arcangle=-45]{->}{b1}{b3}\nbput[npos=0.56]{{\scriptsize $G$}}

{
\rput[c](16,21){\psset{unit=1mm,doubleline=true,arrowinset=0.6,arrowlength=0.5,arrowsize=0.5pt 2.1,nodesep=0pt,labelsep=2pt}
\pcline{->}(0,3)(0,0) \nbput{{\scriptsize $\phi^F$}}}}

{
\rput[c](17,0.5){\psset{unit=1mm,doubleline=true,arrowinset=0.6,arrowlength=0.5,arrowsize=0.5pt 2.1,nodesep=0pt,labelsep=2pt}
\pcline{->}(0,3)(0,0) \nbput{{\scriptsize $\alpha$}}}}


\endpspicture
\]

\item Unit:

\[
\psset{unit=0.15cm,labelsep=2pt,nodesep=2pt}
\pspicture(45,20)


\rput(2,10){\rnode{a1}{$\1$}}
\rput(25,18){\rnode{a2}{$(a,a)$}}
\rput(25,2){\rnode{a3}{$(Fa,Fa)=(Ga,Ga)$}}

\ncline{->}{a1}{a2} \naput{{\scriptsize $I$}}
\ncline{->}{a1}{a3} \nbput{{\scriptsize $I'$}}

\ncarc[arcangle=45]{->}{a2}{a3}\naput[npos=0.5]{{\scriptsize $G$}}
\ncarc[arcangle=-45]{->}{a2}{a3}\nbput[npos=0.5]{{\scriptsize $F$}}

{
\rput[c](14,9){\psset{unit=1mm,doubleline=true,arrowinset=0.6,arrowlength=0.5,arrowsize=0.5pt 2.1,nodesep=0pt,labelsep=-2pt}
\pcline{->}(0,0)(3,3) \naput{{\scriptsize $\phi^F$}}}}

{
\rput[c](24,10){\psset{unit=1mm,doubleline=true,arrowinset=0.6,arrowlength=0.5,arrowsize=0.5pt 2.1,nodesep=0pt,labelsep=2pt}
\pcline{->}(0,0)(4,0) \nbput{{\scriptsize $\alpha$}}}}

\rput(40,10){$=$}

\endpspicture
\psset{unit=0.15cm,labelsep=2pt,nodesep=2pt}
\pspicture(30,20)


\rput(2,10){\rnode{a1}{$\1$}}
\rput(25,18){\rnode{a2}{$(a,a)$}}
\rput(25,2){\rnode{a3}{$(Fa,Fa)=(Ga,Ga)$}}

\ncline{->}{a1}{a2} \naput{{\scriptsize $I$}}
\ncline{->}{a1}{a3} \nbput{{\scriptsize $I'$}}

\ncarc[arcangle=45]{->}{a2}{a3}\naput[npos=0.5]{{\scriptsize $G$}}

{
\rput[c](17,9){\psset{unit=1mm,doubleline=true,arrowinset=0.6,arrowlength=0.5,arrowsize=0.5pt 2.1,nodesep=0pt,labelsep=-2pt}
\pcline{->}(0,0)(3,3) \naput{{\scriptsize $\phi^G$}}}}


\endpspicture\]

\end{itemize}

\end{mydefinition}


\begin{definition} \label{iconcomp}
We now define composition of \cl{B}-icons.  Vertical composition is inherited from $\cB$, that is, components are composed vertically as in $\cB$ and axioms then follow from those in $\cB$. 

Horizontal composition is similar, but the axioms require more care.   Given composable $\cB$-icons

\[
\psset{unit=0.1cm,labelsep=2pt,nodesep=2pt}
\pspicture(50,20)


\rput(0,10){\rnode{a1}{$X$}}  
\rput(20,10){\rnode{a2}{$Y$}}  
\rput(40,10){\rnode{a3}{$Z$}}

\ncarc[arcangle=45]{->}{a1}{a2}\naput{{\scriptsize $F$}}
\ncarc[arcangle=-45]{->}{a1}{a2}\nbput{{\scriptsize $G$}}

\ncarc[arcangle=45]{->}{a2}{a3}\naput{{\scriptsize $H$}}
\ncarc[arcangle=-45]{->}{a2}{a3}\nbput{{\scriptsize $K$}}

{
\rput[c](10,8.5){\psset{unit=1mm,doubleline=true,arrowinset=0.6,arrowlength=0.5,arrowsize=0.5pt 2.1,nodesep=0pt,labelsep=2pt}
\pcline{->}(0,3)(0,0) \naput{{\scriptsize $\alpha$}}}}

{
\rput[c](30,8.5){\psset{unit=1mm,doubleline=true,arrowinset=0.6,arrowlength=0.5,arrowsize=0.5pt 2.1,nodesep=0pt,labelsep=2pt}
\pcline{->}(0,3)(0,0) \naput{{\scriptsize $\beta$}}}}

\endpspicture
\]
the component $(\beta \ast \alpha)_{a,b}$ is given by the horizontal composite

\[
\psset{unit=0.18cm,labelsep=2pt,nodesep=2pt}
\pspicture(50,20)


\rput(0,10){\rnode{a1}{$X(a,b)$}}  
\rput(20,10){\rnode{a2}{$Y(Fa,Fb)$}}  
\rput(40,10){\rnode{a3}{\makebox[4em][l]{$Z(HFa,HFb)$}}}

\ncarc[arcangle=45]{->}{a1}{a2}\naput{{\scriptsize $F$}}
\ncarc[arcangle=-45]{->}{a1}{a2}\nbput{{\scriptsize $G$}}

\ncarc[arcangle=45]{->}{a2}{a3}\naput{{\scriptsize $H$}}
\ncarc[arcangle=-45]{->}{a2}{a3}\nbput{{\scriptsize $K$}}

{
\rput[c](8,9){\psset{unit=1mm,doubleline=true,arrowinset=0.6,arrowlength=0.5,arrowsize=0.5pt 2.1,nodesep=0pt,labelsep=2pt}
\pcline{->}(0,4)(0,0) \naput{{\scriptsize $\alpha_{a,b}$}}}}

{
\rput[c](28,9){\psset{unit=1mm,doubleline=true,arrowinset=0.6,arrowlength=0.5,arrowsize=0.5pt 2.1,nodesep=0pt,labelsep=2pt}
\pcline{->}(0,4)(0,0) \naput{{\scriptsize $\beta_{Fa,Fb}$}}}}

\endpspicture
\]
noting that this makes sense as $F$ and $G$ agree on objects, and likewise $H$ and $K$. We immediately check this satisfies the axioms for a $\cB$-icon.  The unit axiom follows immediately from the unit axioms for $\alpha$ and $\beta$.  The composition axiom is as follows.

\[
\psset{unit=0.1cm,labelsep=2pt,nodesep=3pt}
\pspicture(95,60)


\rput(0,44){\rnode{a1}{$(b,c)(a,b)$}}  
\rput(90,46){\rnode{a3}{\makebox[5em][l]{$(HFb, HFc)(HFa,HFb)$}}}  
\rput(45,20){\rnode{b2}{$(Gb,Gc)(Ga,Gb)$}}  

\rput(0,0){\rnode{c1}{$(a,c)$}}  
\rput(45,0){\rnode{c2}{$(Ga,Gc)$}}  
\rput(90,0){\rnode{c3}{$(KGa,KGc)$}}  


\rput(100,42){\rotatebox{90}{$=$}}
\rput(90,38){\rnode{a3b}{\makebox[5em][l]{$(KGb,KGc)(KGa,KGb)$}}}


\ncline{->}{c1}{c2} \nbput{{\scriptsize $G$}}
\ncline{->}{c2}{c3} \nbput{{\scriptsize $K$}}
\ncline{->}{a1}{c1} \nbput{{\scriptsize $m$}}
\ncline{->}{a3b}{c3} \naput{{\scriptsize $m''$}}
\ncline{->}{b2}{c2} \naput{{\scriptsize $m'$}}

\ncarc[arcangle=20]{->}{a1}{a3} \naput{{\scriptsize $HF.HF$}}
\ncarc[arcangle=-17]{->}{a1}{a3b} \naput{{\scriptsize $KG.KG$}}
\ncarc[arcangle=-15]{->}{a1}{b2} \nbput{{\scriptsize $G.G$}}
\ncarc[arcangle=-15]{->}{b2}{a3b} \nbput{{\scriptsize $K.K$}}

\rput[c](45,28.5){{\scriptsize $\iso$}}
\rput[c](45,26){{\scriptsize\textsf{interchange}}}

{
\rput[c](45,42){\psset{unit=1mm,doubleline=true,arrowinset=0.6,arrowlength=0.5,arrowsize=0.5pt 2.1,nodesep=0pt,labelsep=2pt}
\pcline{->}(0,4)(0,0) \naput{{\scriptsize $(\beta \ast \alpha)(\beta \ast \alpha) $}}}}

{
\rput[c](20,10){\psset{unit=1mm,doubleline=true,arrowinset=0.6,arrowlength=0.5,arrowsize=0.5pt 2.1,nodesep=0pt,labelsep=2pt}
\pcline{->}(0,4)(0,0) \naput{{\scriptsize $\phi^G$}}}}

{
\rput[c](70,10){\psset{unit=1mm,doubleline=true,arrowinset=0.6,arrowlength=0.5,arrowsize=0.5pt 2.1,nodesep=0pt,labelsep=2pt}
\pcline{->}(0,4)(0,0) \naput{{\scriptsize $\phi^K$}}}}


\endpspicture
\]

\[
\psset{unit=0.1cm,labelsep=2pt,nodesep=3pt}
\pspicture(95,63)


\rput(0,44){\rnode{a1}{$(b,c)(a,b)$}}  
\rput(90,46){\rnode{a3}{\makebox[5em][l]{$(HFb, HFc)(HFa,HFb)$}}}  
\rput(45,20){\rnode{b2}{$(Gb,Gc)(Ga,Gb)$}}  

\rput(0,0){\rnode{c1}{$(a,c)$}}  
\rput(45,0){\rnode{c2}{$(Ga,Gc)$}}  
\rput(90,0){\rnode{c3}{$(KGa, KGc)$}}  

\rput(100,42){\rotatebox{90}{$=$}}
\rput(90,38){\rnode{a3b}{\makebox[5em][l]{$(KGb,KGc)(KGa,KGb)$}}}


\rput(45,24){\rotatebox{90}{$=$}}
\rput(45,28){\rnode{bb}{$(Fb,Fc)(Fa,Fb)$}}

\ncline{->}{c1}{c2} \nbput{{\scriptsize $G$}}
\ncline{->}{c2}{c3} \nbput{{\scriptsize $K$}}
\ncline{->}{a1}{c1} \nbput{{\scriptsize $m$}}
\ncline{->}{a3b}{c3} \naput{{\scriptsize $m''$}}
\ncline{->}{b2}{c2} \naput{{\scriptsize $m'$}}

\ncarc[arcangle=20]{->}{a1}{a3} \naput{{\scriptsize $HF.HF$}}
\ncarc[arcangle=-15]{->}{a1}{b2} \nbput{{\scriptsize $G.G$}}
\ncarc[arcangle=-15]{->}{b2}{a3b} \nbput{{\scriptsize $K.K$}}

\ncarc[arcangle=15]{->}{a1}{bb} \naput{{\scriptsize $F.F$}}
\ncarc[arcangle=15]{->}{bb}{a3} \naput{{\scriptsize $H.H$}}

\rput[c](45,40){{\scriptsize $\iso$}}

{
\rput[c](20,32){\psset{unit=1mm,doubleline=true,arrowinset=0.6,arrowlength=0.5,arrowsize=0.5pt 2.1,nodesep=0pt,labelsep=2pt}
\pcline{->}(2,3)(0,0) \naput{{\scriptsize $\alpha.\alpha$}}}}

{
\rput[c](70,32){\psset{unit=1mm,doubleline=true,arrowinset=0.6,arrowlength=0.5,arrowsize=0.5pt 2.1,nodesep=0pt,labelsep=2pt}
\pcline{->}(-2,3)(0,0) \naput{{\scriptsize $\beta.\beta$}}}}

{
\rput[c](20,10){\psset{unit=1mm,doubleline=true,arrowinset=0.6,arrowlength=0.5,arrowsize=0.5pt 2.1,nodesep=0pt,labelsep=2pt}
\pcline{->}(0,4)(0,0) \naput{{\scriptsize $\phi^G$}}}}

{
\rput[c](70,10){\psset{unit=1mm,doubleline=true,arrowinset=0.6,arrowlength=0.5,arrowsize=0.5pt 2.1,nodesep=0pt,labelsep=2pt}
\pcline{->}(0,4)(0,0) \naput{{\scriptsize $\phi^K$}}}}

\rput(-15,20){$=$}


\endpspicture
\]

\[
\psset{unit=0.1cm,labelsep=2pt,nodesep=3pt}
\pspicture(0,-15)(95,63)


\rput(0,44){\rnode{a1}{$(b,c)(a,b)$}}  
\rput(90,46){\rnode{a3}{\makebox[5em][l]{$(HFb, HFc)(HFa,HFb)$}}}  
\rput(45,20){\rnode{b2}{$(Gb,Gc)(Ga,Gb)$}}  

\rput(0,-3){\rnode{c1}{$(a,c)$}}  
\rput(45,0){\rnode{c2}{$(Fa,Fc)$}}  
\rput(90,0){\rnode{c3}{$(HFa, HFc)$}}  


\rput(45,-3){\rotatebox{90}{$=$}}
\rput(45,-6){\rnode{c2b}{$(Ga,Gc)$}}  

\rput(90,-3){\rotatebox{90}{$=$}}
\rput(90,-6){\rnode{c3b}{$(KGa,KGc)$}}

\rput(45,24){\rotatebox{90}{$=$}}
\rput(45,28){\rnode{bb}{$(Fb,Fc)(Fa,Fb)$}}

\ncarc[arcangle=-20]{->}{c1}{c2b} \nbput{{\scriptsize $G$}}
\ncarc[arcangle=-20]{->}{c2b}{c3b} \nbput{{\scriptsize $K$}}

\ncarc[arcangle=20]{->}{c1}{c2} \naput{{\scriptsize $F$}}
\ncarc[arcangle=20]{->}{c2}{c3} \naput{{\scriptsize $H$}}

\ncline{->}{a1}{c1} \nbput{{\scriptsize $m$}}
\ncline{->}{a3}{c3} \naput{{\scriptsize $m''$}}
\ncline{->}{b2}{c2} \naput{{\scriptsize $m'$}}

\ncarc[arcangle=20]{->}{a1}{a3} \naput{{\scriptsize $HF.HF$}}

\ncarc[arcangle=15]{->}{a1}{bb} \naput{{\scriptsize $F.F$}}
\ncarc[arcangle=15]{->}{bb}{a3} \naput{{\scriptsize $H.H$}}

\rput[c](45,40){{\scriptsize $\iso$}}

{
\rput[c](22,-3){\psset{unit=1mm,doubleline=true,arrowinset=0.6,arrowlength=0.5,arrowsize=0.5pt 2.1,nodesep=0pt,labelsep=2pt}
\pcline{->}(0,3)(0,0) \naput{{\scriptsize $\alpha$}}}}

{
\rput[c](67,-5){\psset{unit=1mm,doubleline=true,arrowinset=0.6,arrowlength=0.5,arrowsize=0.5pt 2.1,nodesep=0pt,labelsep=2pt}
\pcline{->}(0,3)(0,0) \naput{{\scriptsize $\beta$}}}}

{
\rput[c](20,20){\psset{unit=1mm,doubleline=true,arrowinset=0.6,arrowlength=0.5,arrowsize=0.5pt 2.1,nodesep=0pt,labelsep=2pt}
\pcline{->}(0,4)(0,0) \naput{{\scriptsize $\phi^F$}}}}

{
\rput[c](70,20){\psset{unit=1mm,doubleline=true,arrowinset=0.6,arrowlength=0.5,arrowsize=0.5pt 2.1,nodesep=0pt,labelsep=2pt}
\pcline{->}(0,4)(0,0) \naput{{\scriptsize $\phi^H$}}}}


\rput(-15,20){$=$}

\endpspicture
\]

\noi where the first equality follows from functoriality axioms and the second from the first \cB-icon axiom applied on each side.

\end{definition}

\begin{proposition}\label{threepointseven}
Weak \cl{B}-categories, weak \cl{B}-functors and \cl{B}-icons organise themselves into a bicategory \cat{\cl{B}-Icon}.  Furthermore if the underlying bicategory of $\cl{B}$ is a 2-category, then \cat{\cl{B}-Icon} is also a 2-category.
\end{proposition}

\begin{prf}
Given a pair of weak $\cB$-categories $X$ and $Y$, the weak $\cB$-functors $X \tra Y$ and $\cB$-icons between those form a category---composition is vertical composition of $\cB$-icons and the identity $1_F: F \Tra F$ has all its components given by identity 2-cells in $\cB$.  The category axioms follow from those for the hom-categories in $\cB$.

We have defined the action of the composition functor on cells in Definitions~\ref{functorcomp} and \ref{iconcomp}.  For functoriality,  interchange of $\cB$-icons follows from interchange of 2-cells in $\cB$. Likewise the fact that horizontal composition of units is a unit follows from the corresponding fact in $\cB$.

We now define the identity weak $\cB$-functor. 

\begin{itemize}
 \item On objects it is the identity function.
\item On homs it is the identity 1-cell in \cB.
\item The constraint cells are given below. 

\end{itemize}

\[
\psset{unit=0.1cm,labelsep=2pt,nodesep=3pt}
\pspicture(40,30)


\rput(0,30){\rnode{a1}{$(b,c)(a,b)$}}  
\rput(40,30){\rnode{a2}{$(b,c)(a,b)$}}  
\rput(0,0){\rnode{b1}{$(a,c)$}}  
\rput(40,0){\rnode{b2}{$(a,c)$}}

\ncline{->}{a1}{a2} \naput{{\scriptsize $1.1$}}
\ncline{->}{b1}{b2} \nbput{{\scriptsize $1$}}
\ncline{->}{a1}{b1} \nbput{{\scriptsize $m$}}
\ncline{->}{a2}{b2} \naput{{\scriptsize $m$}}
\ncarc[arcangle=-15]{->}{a1}{b2} \naput{{\scriptsize $m$}}

\ncarc[arcangle=-35]{->}{a1}{a2} \nbput{{\scriptsize $1$}}

{
\rput[c](20,24){\psset{unit=1mm,doubleline=true,arrowinset=0.6,arrowlength=0.5,arrowsize=0.5pt 2.1,nodesep=0pt,labelsep=2pt}
\pcline{->}(0,3)(0,0) \naput{{\scriptsize $\iso$}}}}

{
\rput[c](10,6){\psset{unit=1mm,doubleline=true,arrowinset=0.6,arrowlength=0.5,arrowsize=0.5pt 2.1,nodesep=0pt,labelsep=2pt}
\pcline{->}(2,3)(0,0) \naput{{\scriptsize $\mbox{{\bfseries r}}^{-1}$}}}}

{
\rput[c](28,14){\psset{unit=1mm,doubleline=true,arrowinset=0.6,arrowlength=0.5,arrowsize=0.5pt 2.1,nodesep=0pt,labelsep=2pt}
\pcline{->}(2,3)(0,0) \naput{{\scriptsize {\bfseries l}}}}}

\endpspicture
\]

\[
\psset{unit=0.1cm,labelsep=2pt,nodesep=3pt,npos=0.4}
\pspicture(20,22)


\rput(0,20){\rnode{a1}{$\1$}}  
\rput(28,20){\rnode{a2}{$(a,a)$}}  
\rput(20,0){\rnode{a3}{$(a,a)$}}  

\ncline{->}{a1}{a2} \naput[npos=0.6]{{\scriptsize $I$}} 
\ncline{->}{a1}{a3} \nbput[npos=0.55]{{\scriptsize $I$}} 
\ncline{->}{a3}{a2} \nbput{{\scriptsize $1$}} 

{
\rput[c](14,11){\psset{unit=1mm,doubleline=true,arrowinset=0.6,arrowlength=0.5,arrowsize=0.5pt 2.1,nodesep=0pt,labelsep=1pt}
\pcline{->}(2,3)(0,0) \naput{{\scriptsize $\mbox{{\bfseries r}}^{-1}$}}}}

\endpspicture
\]
The axioms for a weak \cB-functor follow from coherence for the monoidal bicategory $\cB$.

We now deal with the associativity and unit constraints.  First note that in each case the constraint must be a $\cB$-icon, so has components which are 2-cells of $\cB$.  We simply define these to be the associativity and unit constraints of $\cB$, and it only remains to check that these components satisfy the $\cB$-icon axioms, which is straightforward as we now indicate.  For the first axiom we must check the following equality

\[
\psset{unit=0.1cm,labelsep=2pt,nodesep=3pt}
\pspicture(0,0)(95,60)

%

\rput(0,44){\rnode{a1}{$\ $}}  
\rput(90,44){\rnode{a3}{$\ $}}  
\rput(45,20){\rnode{b2}{$\ $}}  

\rput(22,15){\rnode{d4}{$\ $}}  
\rput(22,0){\rnode{c4}{$\ $}}

\rput(0,0){\rnode{c1}{$\ $}}  
\rput(45,0){\rnode{c2}{$\ $}}  
\rput(90,0){\rnode{c3}{$\ $}}  

\ncarc[arcangle=-10]{->}{a1}{d4} \nbput{{\scriptsize $F.F$}}
\ncarc[arcangle=-10]{->}{d4}{b2} \nbput{{\scriptsize $G.G$}}
\ncline{->}{d4}{c4} \nbput{{\scriptsize $m'$}}

\ncline{->}{c1}{c4} \nbput{{\scriptsize $F$}}
\ncline{->}{c4}{c2} \nbput{{\scriptsize $G$}}

\ncline{->}{c2}{c3} \nbput{{\scriptsize $H$}}
\ncline{->}{a1}{c1} \nbput{{\scriptsize $m$}}
\ncline{->}{a3}{c3} \naput{{\scriptsize $m'''$}}
\ncline{->}{b2}{c2} \naput{{\scriptsize $m''$}}

\ncarc[arcangle=17]{->}{a1}{a3} \naput{{\scriptsize $(HG)F.(HG)F$}}
\ncarc[arcangle=-17]{->}{a1}{a3} \naput{{\scriptsize $H(GF).H(GF)$}}
\ncarc[arcangle=-5]{->}{a1}{b2} \naput{{\scriptsize $GF.GF$}}
\ncarc[arcangle=-5]{->}{b2}{a3} \nbput{{\scriptsize $H.H$}}

\rput[c](25,24.5){{\scriptsize $\iso$}}
\rput[c](25,22){{\scriptsize\textsf{interchange}}}

\rput[c](45,30.5){{\scriptsize $\iso$}}
\rput[c](45,28){{\scriptsize\textsf{interchange}}}

{
\rput[c](45,43){\psset{unit=1mm,doubleline=true,arrowinset=0.6,arrowlength=0.5,arrowsize=0.5pt 2.1,nodesep=0pt,labelsep=2pt}
\pcline{->}(0,4)(0,0) \naput{{\scriptsize $a.a $}}}}

{
\rput[c](10,10){\psset{unit=1mm,doubleline=true,arrowinset=0.6,arrowlength=0.5,arrowsize=0.5pt 2.1,nodesep=0pt,labelsep=2pt}
\pcline{->}(0,4)(0,0) \naput{{\scriptsize $\phi^F$}}}}

{
\rput[c](33,5){\psset{unit=1mm,doubleline=true,arrowinset=0.6,arrowlength=0.5,arrowsize=0.5pt 2.1,nodesep=0pt,labelsep=2pt}
\pcline{->}(0,4)(0,0) \naput{{\scriptsize $\phi^G$}}}}

{
\rput[c](70,10){\psset{unit=1mm,doubleline=true,arrowinset=0.6,arrowlength=0.5,arrowsize=0.5pt 2.1,nodesep=0pt,labelsep=2pt}
\pcline{->}(0,4)(0,0) \naput{{\scriptsize $\phi^H$}}}}


\endpspicture
\]

\[
\psset{unit=0.1cm,labelsep=2pt,nodesep=3pt}
\pspicture(0,-20)(95,60)

%

\rput(0,44){\rnode{a1}{$\ $}}  
\rput(90,44){\rnode{a3}{$\ $}}  
\rput(45,20){\rnode{b2}{$\ $}}  

\rput(67,15){\rnode{d4}{$\ $}}  
\rput(67,0){\rnode{c4}{$\ $}}

\rput(0,0){\rnode{c1}{$\ $}}  
\rput(45,0){\rnode{c2}{$\ $}}  
\rput(90,0){\rnode{c3}{$\ $}}  

\ncarc[arcangle=-10]{->}{b2}{d4} \nbput{{\scriptsize $G.G$}}
\ncarc[arcangle=-10]{->}{d4}{a3} \nbput{{\scriptsize $H.H$}}
\ncline{->}{d4}{c4} \naput{{\scriptsize $m''$}}

\ncline{->}{c1}{c2} \nbput{{\scriptsize $F$}}
\ncline{->}{c2}{c4} \nbput{{\scriptsize $G$}}

\ncline{->}{c4}{c3} \nbput{{\scriptsize $H$}}
\ncline{->}{a1}{c1} \nbput{{\scriptsize $m$}}
\ncline{->}{a3}{c3} \naput{{\scriptsize $m'''$}}
\ncline{->}{b2}{c2} \nbput{{\scriptsize $m'$}}

\ncarc[arcangle=17]{->}{a1}{a3} \naput{{\scriptsize $(HG)F.(HG)F$}}
\ncarc[arcangle=-35]{->}{c1}{c3} \nbput{{\scriptsize $H(GF)$}}
\ncarc[arcangle=-5]{->}{a1}{b2} \naput{{\scriptsize $F.F$}}
\ncarc[arcangle=-5]{->}{b2}{a3} \naput{{\scriptsize $HG.HG$}}

\rput[c](66,24.5){{\scriptsize $\iso$}}
\rput[c](66,22){{\scriptsize\textsf{interchange}}}

\rput[c](45,40.5){{\scriptsize $\iso$}}
\rput[c](45,38){{\scriptsize\textsf{interchange}}}

{
\rput[c](45,-10){\psset{unit=1mm,doubleline=true,arrowinset=0.6,arrowlength=0.5,arrowsize=0.5pt 2.1,nodesep=0pt,labelsep=2pt}
\pcline{->}(0,4)(0,0) \naput{{\scriptsize $a$}}}}

{
\rput[c](22,10){\psset{unit=1mm,doubleline=true,arrowinset=0.6,arrowlength=0.5,arrowsize=0.5pt 2.1,nodesep=0pt,labelsep=2pt}
\pcline{->}(0,4)(0,0) \naput{{\scriptsize $\phi^F$}}}}

{
\rput[c](55,5){\psset{unit=1mm,doubleline=true,arrowinset=0.6,arrowlength=0.5,arrowsize=0.5pt 2.1,nodesep=0pt,labelsep=2pt}
\pcline{->}(0,4)(0,0) \naput{{\scriptsize $\phi^G$}}}}

{
\rput[c](78,10){\psset{unit=1mm,doubleline=true,arrowinset=0.6,arrowlength=0.5,arrowsize=0.5pt 2.1,nodesep=0pt,labelsep=2pt}
\pcline{->}(0,4)(0,0) \naput{{\scriptsize $\phi^H$}}}}

\rput(-15,20){$=$}


\endpspicture
\]
which follows from coherence for the monoidal bicategory $\cB$ and naturality of $a$, the associator for the 1-cells in $\cB$.  The unit axiom follows similarly.

Finally, to verify the bicategory axioms for \cat{\cl{B}-Icon} we only need to check them componentwise, so they all follow from the bicategory axioms for $\cB$.  Since the coherence constraints are built from those in $\cl{B}$, it follows that if the constraints in $\cl{B}$ are identities then the constraints in \cat{$\cB$-Icon} are also identities.  Thus if $\cB$ is a 2-category then \cat{$\cB$-Icon} is also a 2-category.  \end{prf}

\begin{theorem}  \label{smb}

If \cl{B} is a symmetric monoidal bicategory then so is \cat{\cl{B}-Icon}.  Thus in particular the icon construction can be iterated.
 
\end{theorem}

\begin{prf}
Suppose $\cB$ is a symmetric monoidal bicategory, with symmetry
\[R= R_{ab}: a \otimes b \tra b \otimes a.\]
We already know that \Bicon\ is a bicategory, and we now show that it is symmetric monoidal. 

We define a (weak) functor
\[\otimes \: \Bicon \times \Bicon \tra \Bicon.\]

\subsubsection*{Definition of $\otimes$ on weak $\cB$-categories.} 

Let $X,Y$ be weak $\cB$-categories.  Then $X \otimes Y$ has
\begin{itemize}
 \item objects pairs $\langle x,y \rangle \in \ob X \times \ob Y$, 

\item $(X \otimes Y) \big( \langle x_0, y_0 \rangle, \langle x_1, y_1 \rangle \big) = X(x_0, x_1).Y(y_0, y_1)$ where on the right the dot denotes the tensor product in $\cB$, 

\item composition: given $x_0, x_1, x_2 \in X_0$ and $y_0, y_1, y_2 \in Y_0$, composition is defined as below

\[
\psset{unit=0.1cm,labelsep=2pt,nodesep=2pt}
\pspicture(10,40)


\rput(5,3){\rnode{a1}{$X(x_0,x_2).Y(y_0, y_2)$,}}  %
\rput(5,20){\rnode{a2}{$X(x_1, x_2).X(x_0,x_1).Y(y_1,y_2).Y(y_0,y_1)$}}  %
\rput(5,37){\rnode{a3}{$X(x_1, x_2).Y(y_1,y_2).X(x_0,x_1).Y(y_0,y_1)$}}

\ncline{->}{a3}{a2} \naput{{\scriptsize $1.R.1$}}
\ncline{->}{a2}{a1} \naput{{\scriptsize $m.m$}}

\endpspicture\]

\item identities: for all $x \in X, y \in Y$, a 1-cell
\[ \1 \tra \1.\1 \tmap{I_x.I_y} X(x,x).Y(y,y),\]

\item associativity constraints: 

\vspace{5em}

\rotatebox{90}{
\begin{minipage}{50em}
\psset{unit=0.1cm,labelsep=2pt,nodesep=3pt,npos=0.4}
\pspicture(100,40)
\rput(0,40){\rnode{a1}{$(x_2,x_3)(y_2,y_3)(x_1,x_2)(y_1,y_2)(x_0,x_1)(y_0,y_1)$}}  
\rput(90,40){\rnode{a2}{$(x_2,x_3)(x_1,x_2)(y_2,y_3)(y_1,y_2)(x_0,x_1)(y_0,y_1)$}}  
\rput(170,40){\rnode{a3}{$(x_1,x_3)(y_1,y_3)(x_0,x_1)(y_0,y_1)$}}  
\rput(0,20){\rnode{b1}{$(x_2,x_3)(y_2,y_3)(x_1,x_2)(x_0,x_1)(y_1,y_2)(y_0,y_1)$}}  
\rput(90,20){\rnode{b2}{$(x_2,x_3)(x_1,x_2)(x_0,x_1)(y_2,y_3)(y_1,y_2)(y_0,y_1)$}}  
\rput(170,20){\rnode{b3}{$(x_1,x_3)(x_0,x_1)(y_1,y_3)(y_0,y_1)$}}  
\rput(0,0){\rnode{c1}{$(x_2,x_3)(y_2,y_3)(x_0,x_2)(y_0,y_2)$}}  
\rput(90,0){\rnode{c2}{$(x_2,x_3)(x_0,x_2)(y_2,y_3)(y_0,y_2)$}}  
\rput(170,0){\rnode{c3}{$(x_0,x_3)(y_0,y_3)$}}  
\ncline{->}{a1}{a2} \naput{{\scriptsize $1.R.1.1.1$}} 
\ncline{->}{a2}{a3} \naput{{\scriptsize $m.m.1.1$}} 
\ncline{->}{b1}{b2} \naput{{\scriptsize $1.R.1.1$}} 
\ncline{->}{b2}{b3} \naput{{\scriptsize $m.1.m.1$}} 
\ncline{->}{c1}{c2} \nbput{{\scriptsize $1.R.1$}} 
\ncline{->}{c2}{c3} \nbput{{\scriptsize $m.m$}} 
\ncline{->}{a1}{b1} \nbput{{\scriptsize $1.1.1.R.1$}} 
\ncline{->}{b1}{c1} \nbput{{\scriptsize $1.1.m.m$}} 
\ncline{->}{a2}{b2} \naput{{\scriptsize $1.1.R.1$}} 
\ncline{->}{b2}{c2} \naput{{\scriptsize $1.m.1.m$}} 
\ncline{->}{a3}{b3} \naput{{\scriptsize $1.R.1$}} 
\ncline{->}{b3}{c3} \naput{{\scriptsize $m.m$}} 
{
\rput[c](130,8){\psset{unit=1mm,doubleline=true,arrowinset=0.6,arrowlength=0.5,arrowsize=0.5pt 2.1,nodesep=0pt,labelsep=2pt}
\pcline{->}(0,4)(0,0) \naput{{\scriptsize $\mbox{\bfseries{a}}.\mbox{\bfseries{a}}$}}}}
\rput[c](45,12){{\scriptsize $\iso$}}
\rput[c](45,9){{\scriptsize\textsf{naturality of $r$}}}
\rput[c](130,32){{\scriptsize $\iso$}}
\rput[c](130,29){{\scriptsize\textsf{naturality of $r$}}}
\rput[c](45,32){{\scriptsize $\iso$}}
\rput[c](45,29){{\scriptsize\textsf{braid coherence}}}
\endpspicture
\end{minipage}}

\vspace{10em}

Note that in the above diagram, we have written $R$ for various different instances of braidings.  This is unambiguous given the sources and targets, but requires some care especially in the top left square.

\item unit constraints:

\[
\psset{unit=0.1cm,labelsep=2pt,nodesep=3pt}
\pspicture(80,40)


\rput(0,40){\rnode{a1}{$(x_0,x_1)(y_0,y_1).\1$}}  
\rput(45,40){\rnode{a2}{$(x_0,x_1)(y_0,y_1).\1.\1$}}  
\rput(100,40){\rnode{a3}{$(x_0,x_1)(y_0,y_1)(x_0,x_0)(y_0,y_0)$}}
  
\rput(45,20){\rnode{b2}{$(x_0,x_1).\1.(y_0,y_1).\1$}}  
\rput(100,20){\rnode{b3}{$(x_0,x_1)(x_0,x_0)(y_0,y_1)(y_0,y_0)$}}  
\rput(100,0){\rnode{c3}{$(x_0,x_1)(y_0,y_1)$}}

\ncline{->}{a1}{a2} \naput{{\scriptsize $\sim$}}
\ncline{->}{a2}{a3} \naput{{\scriptsize $1.1.I.I$}}
\ncline{->}{b2}{b3} \naput{{\scriptsize $1.I.1.I$}}
\ncline{->}{a2}{b2} \nbput{{\scriptsize $1.R.1$}}
\ncline{->}{a3}{b3} \naput{{\scriptsize $1.R.1$}}
\ncline{->}{b3}{c3} \naput{{\scriptsize $m.m$}}
\ncline{->}{b2}{c3} \nbput{{\scriptsize $\sim$}}

\ncarc[arcangleA=-40,arcangleB=-25]{->}{a1}{c3} \nbput{{\scriptsize $\sim$}}

\rput[c](70,32){{\scriptsize $\iso$}}
\rput[c](70,29){{\scriptsize\textsf{naturality of $R$}}}

\rput[c](22,32){{\scriptsize $\iso$}}
\rput[c](22,29.5){{\scriptsize\textsf{coherence}}}

{
\rput[c](84,11){\psset{unit=1mm,doubleline=true,arrowinset=0.6,arrowlength=0.5,arrowsize=0.5pt 2.1,nodesep=0pt,labelsep=2pt}
\pcline{->}(2,3)(0,0) \naput{{\scriptsize {\bfseries r.r}}}}}

\endpspicture
\]

\[
\psset{unit=0.1cm,labelsep=2pt,nodesep=3pt}
\pspicture(80,45)


\rput(0,40){\rnode{a1}{$\1.(x_0,x_1)(y_0,y_1)$}}  
\rput(45,40){\rnode{a2}{$\1.\1.(x_0,x_1)(y_0,y_1)$}}  
\rput(100,40){\rnode{a3}{$(x_1,x_1)(y_1,y_1)(x_0,x_1)(y_0,y_1)$}}
  
\rput(45,20){\rnode{b2}{$\1.(x_0,x_1).\1.(y_0,y_1)$}}  
\rput(100,20){\rnode{b3}{$(x_1,x_1)(x_0,x_1)(y_1,y_1)(y_0,y_1)$}}  
\rput(100,0){\rnode{c3}{$(x_0,x_1)(y_0,y_1)$}}

\ncline{->}{a1}{a2} \naput{{\scriptsize $\sim$}}
\ncline{->}{a2}{a3} \naput{{\scriptsize $I.I.1.1$}}
\ncline{->}{b2}{b3} \naput{{\scriptsize $I.1.I.1$}}
\ncline{->}{a2}{b2} \nbput{{\scriptsize $1.R.1$}}
\ncline{->}{a3}{b3} \naput{{\scriptsize $1.R.1$}}
\ncline{->}{b3}{c3} \naput{{\scriptsize $m.m$}}
\ncline{->}{b2}{c3} \nbput{{\scriptsize $\sim$}}

\ncarc[arcangleA=-40,arcangleB=-25]{->}{a1}{c3} \nbput{{\scriptsize $\sim$}}

\rput[c](70,32){{\scriptsize $\iso$}}
\rput[c](70,29){{\scriptsize\textsf{naturality of $R$}}}

\rput[c](22,32){{\scriptsize $\iso$}}
\rput[c](22,29.5){{\scriptsize\textsf{coherence}}}

{
\rput[c](84,11){\psset{unit=1mm,doubleline=true,arrowinset=0.6,arrowlength=0.5,arrowsize=0.5pt 2.1,nodesep=0pt,labelsep=2pt}
\pcline{->}(2,3)(0,0) \naput{{\scriptsize {\bfseries l.l}}}}}

\endpspicture
\]

\end{itemize}

We check that this data for $X \otimes Y$ satisfies the axioms for a weak $\cB$-category.  Ignoring all coherence in $\cB$, these are not difficult to check;  they follow from the axioms for $X$ and $Y$ tensored together.  The following diagram indicates how to prove the associativity axiom; we again use the cubical version of the axioms, as introduced in Definition~\ref{tractable}, and omit almost all the labels as they can be inferred.



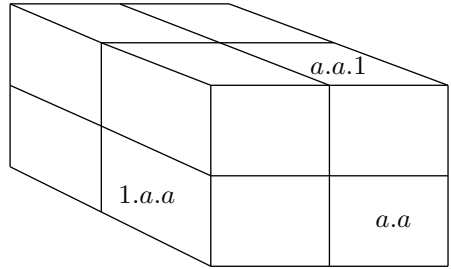
\[\psset{unit=0.12cm,labelsep=2pt,nodesep=0pt}
\pspicture(0,8)(60,40)

\rput(23,18){$1.a.a$}  
\rput(44,32.2){$a.a.1$}  
\rput(50,15){$a.a$}

\rput(30,10){\rnode{a}{$$}}  
\rput(43,10){\rnode{b}{$$}}  
\rput(56,10){\rnode{c}{$$}}  
\rput(30,20){\rnode{d}{$$}}  
\rput(43,20){\rnode{e}{$$}}  
\rput(56,20){\rnode{f}{$$}}  
\rput(30,30){\rnode{g}{$$}}  
\rput(43,30){\rnode{h}{$$}}  
\rput(56,30){\rnode{i}{$$}}  
\rput(18,16){\rnode{j}{$$}}  
\rput(18,25.5){\rnode{k}{$$}}  
\rput(18,34.7){\rnode{l}{$$}}  
\rput(31,34.7){\rnode{m}{$$}}  
\rput(43.5,34.7){\rnode{n}{$$}}  
\rput(8,21){\rnode{p}{$$}}  
\rput(8,30){\rnode{q}{$$}}  
\rput(8,39){\rnode{r}{$$}}  
\rput(20,39){\rnode{s}{$$}}  
\rput(32,39){\rnode{t}{$$}}

\ncline{a}{b} 
\ncline{b}{c} 
\ncline{d}{e} 
\ncline{e}{f} 
\ncline{g}{h} 
\ncline{h}{i} 

\ncline{a}{d} 
\ncline{d}{g} 
\ncline{b}{e} 
\ncline{e}{h} 
\ncline{c}{f} 
\ncline{f}{i}
 
\ncline{j}{k} 
\ncline{k}{l} 
\ncline{l}{m} 
\ncline{m}{n} 
\ncline{j}{a} 
\ncline{k}{d} 
\ncline{l}{g} 
\ncline{m}{h} 
\ncline{n}{i}
 
\ncline{p}{q} 
\ncline{q}{r} 
\ncline{r}{s} 
\ncline{s}{t} 
\ncline{p}{j} 
\ncline{q}{k} 
\ncline{r}{l} 
\ncline{s}{m} 
\ncline{t}{n}

\endpspicture
\]


\[\psset{unit=0.12cm,labelsep=2pt,nodesep=0pt}
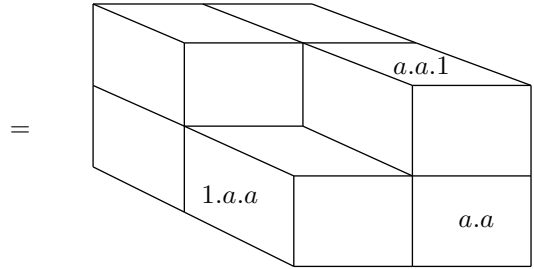
\pspicture(0,8)(60,40)

\rput(0,25){$=$}

\rput(23,18){$1.a.a$}  
\rput(44,32.2){$a.a.1$}  
\rput(50,15){$a.a$}

\rput(30,10){\rnode{a}{$$}}  
\rput(43,10){\rnode{b}{$$}}  
\rput(56,10){\rnode{c}{$$}}  
\rput(30,20){\rnode{d}{$$}}  
\rput(43,20){\rnode{e}{$$}}  
\rput(56,20){\rnode{f}{$$}}  
\rput(31,25.5){\rnode{g}{$$}}  
\rput(43,30){\rnode{h}{$$}}  
\rput(56,30){\rnode{i}{$$}}  
\rput(18,16){\rnode{j}{$$}}  
\rput(18,25.5){\rnode{k}{$$}}  
\rput(18,34.7){\rnode{l}{$$}}  
\rput(31,34.7){\rnode{m}{$$}}  
\rput(43.4,34.7){\rnode{n}{$$}}  
\rput(8,21){\rnode{p}{$$}}  
\rput(8,30){\rnode{q}{$$}}  
\rput(8,39){\rnode{r}{$$}}  
\rput(20,39){\rnode{s}{$$}}  
\rput(32,39){\rnode{t}{$$}}

\ncline{a}{b} 
\ncline{b}{c} 
\ncline{d}{e} 
\ncline{e}{f} 
\ncline{g}{e} 
\ncline{h}{i} 

\ncline{a}{d} 
\ncline{k}{g} 
\ncline{b}{e} 
\ncline{e}{h} 
\ncline{c}{f} 
\ncline{f}{i}
 
\ncline{j}{k} 
\ncline{k}{l} 
\ncline{l}{m} 
\ncline{m}{n} 
\ncline{j}{a} 
\ncline{k}{d} 
\ncline{m}{g} 
\ncline{m}{h} 
\ncline{n}{i}
 
\ncline{p}{q} 
\ncline{q}{r} 
\ncline{r}{s} 
\ncline{s}{t} 
\ncline{p}{j} 
\ncline{q}{k} 
\ncline{r}{l} 
\ncline{s}{m} 
\ncline{t}{n}

\endpspicture
\]


\[\psset{unit=0.12cm,labelsep=2pt,nodesep=0pt}
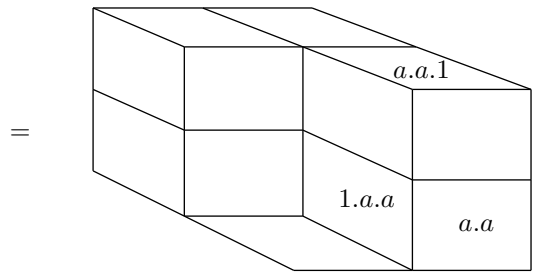
\pspicture(0,8)(60,40)

\rput(0,25){$=$}

\rput(38,18){$1.a.a$}  
\rput(44,32.2){$a.a.1$}  
\rput(50,15){$a.a$}

\rput(30,10){\rnode{a}{$$}}  
\rput(43,10){\rnode{b}{$$}}  
\rput(56,10){\rnode{c}{$$}}  
\rput(31,16){\rnode{d}{$$}}  
\rput(43,20){\rnode{e}{$$}}  
\rput(56,20){\rnode{f}{$$}}  
\rput(31,25.5){\rnode{g}{$$}}  
\rput(43,30){\rnode{h}{$$}}  
\rput(56,30){\rnode{i}{$$}}  
\rput(18,16){\rnode{j}{$$}}  
\rput(18,25.5){\rnode{k}{$$}}  
\rput(18,34.7){\rnode{l}{$$}}  
\rput(31,34.7){\rnode{m}{$$}}  
\rput(43.4,34.7){\rnode{n}{$$}}  
\rput(8,21){\rnode{p}{$$}}  
\rput(8,30){\rnode{q}{$$}}  
\rput(8,39){\rnode{r}{$$}}  
\rput(20,39){\rnode{s}{$$}}  
\rput(32,39){\rnode{t}{$$}}

\ncline{a}{b} 
\ncline{b}{c} 
\ncline{d}{b} 
\ncline{e}{f} 
\ncline{g}{e} 
\ncline{h}{i} 

\ncline{g}{d} 
\ncline{k}{g} 
\ncline{b}{e} 
\ncline{e}{h} 
\ncline{c}{f} 
\ncline{f}{i}
 
\ncline{j}{k} 
\ncline{k}{l} 
\ncline{l}{m} 
\ncline{m}{n} 
\ncline{j}{a} 
\ncline{j}{d} 
\ncline{m}{g} 
\ncline{m}{h} 
\ncline{n}{i}
 
\ncline{p}{q} 
\ncline{q}{r} 
\ncline{r}{s} 
\ncline{s}{t} 
\ncline{p}{j} 
\ncline{q}{k} 
\ncline{r}{l} 
\ncline{s}{m} 
\ncline{t}{n}

\endpspicture
\]


\[\psset{unit=0.12cm,labelsep=2pt,nodesep=0pt}
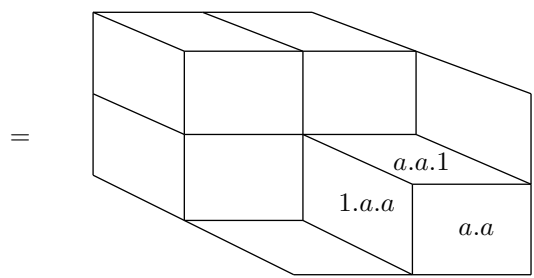
\pspicture(0,8)(60,45)

\rput(0,25){$=$}

\rput(38,18){$1.a.a$}  
\rput(44,22.5){$a.a.1$}  
\rput(50,15){$a.a$}

\rput(30,10){\rnode{a}{$$}}  
\rput(43,10){\rnode{b}{$$}}  
\rput(56,10){\rnode{c}{$$}}  
\rput(31,16){\rnode{d}{$$}}  
\rput(43,20){\rnode{e}{$$}}  
\rput(56,20){\rnode{f}{$$}}  
\rput(31,25.5){\rnode{g}{$$}}  
\rput(43.4,25.5){\rnode{h}{$$}}  
\rput(56,30){\rnode{i}{$$}}  
\rput(18,16){\rnode{j}{$$}}  
\rput(18,25.5){\rnode{k}{$$}}  
\rput(18,34.7){\rnode{l}{$$}}  
\rput(31,34.7){\rnode{m}{$$}}  
\rput(43.4,34.7){\rnode{n}{$$}}  
\rput(8,21){\rnode{p}{$$}}  
\rput(8,30){\rnode{q}{$$}}  
\rput(8,39){\rnode{r}{$$}}  
\rput(20,39){\rnode{s}{$$}}  
\rput(32,39){\rnode{t}{$$}}

\ncline{a}{b} 
\ncline{b}{c} 
\ncline{d}{b} 
\ncline{e}{f} 
\ncline{g}{e} 
\ncline{h}{f} 

\ncline{g}{d} 
\ncline{k}{g} 
\ncline{b}{e} 
\ncline{g}{h} 
\ncline{c}{f} 
\ncline{f}{i}
 
\ncline{j}{k} 
\ncline{k}{l} 
\ncline{l}{m} 
\ncline{m}{n} 
\ncline{j}{a} 
\ncline{j}{d} 
\ncline{m}{g} 
\ncline{n}{h} 
\ncline{n}{i}
 
\ncline{p}{q} 
\ncline{q}{r} 
\ncline{r}{s} 
\ncline{s}{t} 
\ncline{p}{j} 
\ncline{q}{k} 
\ncline{r}{l} 
\ncline{s}{m} 
\ncline{t}{n}

\endpspicture
\]


\[\psset{unit=0.12cm,labelsep=2pt,nodesep=0pt}
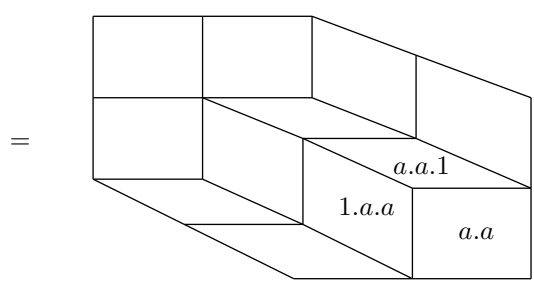
\pspicture(0,8)(60,40)

\rput(0,25){$=$}  

\rput(38,18){$1.a.a$}  
\rput(44,22.5){$a.a.1$}  
\rput(50,15){$a.a$}

\rput(30,10){\rnode{a}{$$}}  
\rput(43,10){\rnode{b}{$$}}  
\rput(56,10){\rnode{c}{$$}}  
\rput(31,16){\rnode{d}{$$}}  
\rput(43,20){\rnode{e}{$$}}  
\rput(56,20){\rnode{f}{$$}}  
\rput(31,25.5){\rnode{g}{$$}}  
\rput(43.4,25.5){\rnode{h}{$$}}  
\rput(56,30){\rnode{i}{$$}}  
\rput(18,16){\rnode{j}{$$}}  
\rput(20,21){\rnode{k}{$$}}  
\rput(20,30){\rnode{l}{$$}}  
\rput(32,30){\rnode{m}{$$}}  
\rput(43.4,34.7){\rnode{n}{$$}}  
\rput(8,21){\rnode{p}{$$}}  
\rput(8,30){\rnode{q}{$$}}  
\rput(8,39){\rnode{r}{$$}}  
\rput(20,39){\rnode{s}{$$}}  
\rput(32,39){\rnode{t}{$$}}

\ncline{a}{b} 
\ncline{b}{c} 
\ncline{d}{b} 
\ncline{e}{f} 
\ncline{g}{e} 
\ncline{h}{f} 

\ncline{g}{d} 
\ncline{k}{d} 
\ncline{b}{e} 
\ncline{g}{h} 
\ncline{c}{f} 
\ncline{f}{i}
 
\ncline{q}{l} 
\ncline{k}{l} 
\ncline{l}{m} 
\ncline{m}{t} 
\ncline{j}{a} 
\ncline{j}{d} 
\ncline{l}{g} 
\ncline{n}{h} 
\ncline{n}{i}
 
\ncline{p}{q} 
\ncline{q}{r} 
\ncline{r}{s} 
\ncline{s}{t} 
\ncline{p}{j} 
\ncline{p}{k} 
\ncline{s}{l} 
\ncline{h}{m} 
\ncline{t}{n}

\endpspicture
\]


\[\psset{unit=0.12cm,labelsep=2pt,nodesep=0pt}
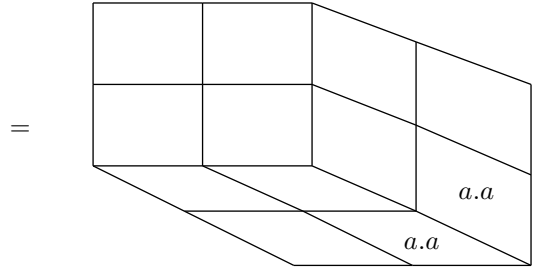
\pspicture(0,8)(60,40)

\rput(0,25){$=$}  

\rput(44,12.5){$a.a$}  
\rput(50,18){$a.a$}

\rput(30,10){\rnode{a}{$$}}  
\rput(43,10){\rnode{b}{$$}}  
\rput(56,10){\rnode{c}{$$}}  
\rput(31,16){\rnode{d}{$$}}  
\rput(43.4,16){\rnode{e}{$$}}  
\rput(56,20){\rnode{f}{$$}}  
\rput(32,21){\rnode{g}{$$}}  
\rput(43.4,25.5){\rnode{h}{$$}}  
\rput(56,30){\rnode{i}{$$}}  
\rput(18,16){\rnode{j}{$$}}  
\rput(20,21){\rnode{k}{$$}}  
\rput(20,30){\rnode{l}{$$}}  
\rput(32,30){\rnode{m}{$$}}  
\rput(43.4,34.7){\rnode{n}{$$}}  
\rput(8,21){\rnode{p}{$$}}  
\rput(8,30){\rnode{q}{$$}}  
\rput(8,39){\rnode{r}{$$}}  
\rput(20,39){\rnode{s}{$$}}  
\rput(32,39){\rnode{t}{$$}}

\ncline{a}{b} 
\ncline{b}{c} 
\ncline{d}{b} 
\ncline{e}{c} 
\ncline{d}{e} 
\ncline{h}{f} 

\ncline{g}{e} 
\ncline{k}{d} 
\ncline{h}{e} 
\ncline{g}{k} 
\ncline{c}{f} 
\ncline{f}{i}
 
\ncline{q}{l} 
\ncline{k}{l} 
\ncline{l}{m} 
\ncline{m}{t} 
\ncline{j}{a} 
\ncline{j}{d} 
\ncline{m}{g} 
\ncline{n}{h} 
\ncline{n}{i}
 
\ncline{p}{q} 
\ncline{q}{r} 
\ncline{r}{s} 
\ncline{s}{t} 
\ncline{p}{j} 
\ncline{p}{k} 
\ncline{s}{l} 
\ncline{h}{m} 
\ncline{t}{n}

\endpspicture
\]

This completes the definition of $X \otimes Y$ as a weak $\cB$-category.

\subsubsection*{Definition of $\otimes$ on weak $\cB$-functors.} 

Given weak $\cB$-functors
\[X \tmap{F} X'\]
\[Y \tmap{G} Y'\]
we define
\[X \otimes Y \tmap{F \otimes G} X' \otimes Y'\]
as follows:

\begin{itemize}
 \item on objects: $(F \otimes G) \: \langle  x,y \rangle \tmapsto \langle Fx, Gy \rangle$

\item on homs: $ X(x_0, y_0).Y(y_0, y_1) \tmap{F.G} X'(Fx_0,Fx_1).Y'(Gy_0,Gy_1)$ 

\item constraints 

\[\psset{unit=0.08cm,labelsep=0pt,nodesep=3pt}
\pspicture(90,45)


\rput(0,40){\rnode{a1}{\makebox[6em][r]{$(x_1,x_2)(y_1,y_2)(x_0,x_1)(y_0,y_1)$}}} 
\rput(70,40){\rnode{a2}{\makebox[8em][l]{$(Fx_1,Fx_2)(Gy_1,Gy_2)(Fx_0,Fx_1)(Gy_0,Gy_1)$}}} 

\rput(0,20){\rnode{b1}{\makebox[6em][r]{$(x_1,x_2)(x_0,x_1)(y_1,y_2)(y_0,y_1)$}}}   
\rput(70,20){\rnode{b2}{\makebox[8em][l]{$(Fx_1,Fx_2)(Fx_0,Fx_1)(Gy_1,Gy_2)(Gy_0,Gy_1)$}}}  

\rput(0,0){\rnode{c1}{\makebox[6em][r]{$(x_0,x_2)(y_0,y_2)$}}}   
\rput(70,0){\rnode{c2}{\makebox[8em][l]{$(Fx_0,Fx_2)(Gy_0,Gy_2)$}}}  

\psset{nodesep=3pt,labelsep=2pt,arrows=->}
\ncline{a1}{a2}\naput{{\scriptsize $FGFG$}} 
\ncline{b1}{b2}\naput{{\scriptsize $FFGG$}} 
\ncline{c1}{c2}\nbput{{\scriptsize $FG$}} 

\ncline{a1}{b1}\nbput{{\scriptsize $1R1$}} 
\ncline{b1}{c1}\nbput{{\scriptsize $mm$}} 

\ncline{a2}{b2}\naput{{\scriptsize $1R1$}} 
\ncline{b2}{c2}\naput{{\scriptsize $mm$}} 

\psset{labelsep=1.5pt}

{
\rput[c](33,8){\psset{unit=1mm,doubleline=true,arrowinset=0.6,arrowlength=0.5,arrowsize=0.5pt 2.1,nodesep=0pt,labelsep=2pt}
\pcline{->}(0,4)(0,0) \naput{{\scriptsize $\phi\phi$}}}}

\rput[c](33,33){{\scriptsize $\iso$}}
\rput[c](33,30){{\scriptsize\textsf{naturality of $R$} }}

\endpspicture\]

\[
\psset{unit=0.1cm,labelsep=2pt,nodesep=3pt,npos=0.4}
\pspicture(20,22)


\rput(0,20){\rnode{a1}{$\1$}}  
\rput(26,20){\rnode{a2}{$(x,x)(y,y)$}}  
\rput(26,0){\rnode{a3}{$(Fx,Fx)(Gy,Gy)$}}  

\ncline{->}{a1}{a2} \naput[npos=0.6]{{\scriptsize $II$}} 
\ncline{->}{a1}{a3} \nbput[npos=0.55]{{\scriptsize $I'I'$}} 
\ncline{->}{a2}{a3} \naput{{\scriptsize $FG$}} 

{
\rput[c](17,11){\psset{unit=1mm,doubleline=true,arrowinset=0.6,arrowlength=0.5,arrowsize=0.5pt 2.1,nodesep=0pt,labelsep=1pt}
\pcline{->}(0,0)(3,3) \naput{{\scriptsize $\phi\phi$}}}}

\endpspicture
\]

\end{itemize}

We now check that this data for $F\otimes G$ satisfies the axioms for a weak $\cB$-functor.  The calculations are similar to those for the weak $\cB$-category axioms for $X \otimes Y$ above, essentially by tensoring together the axioms for $F$ and $G$, and using coherence for $\cB$ \emph{qua} braided monoidal category.

\subsubsection*{Definition of $\otimes$ on $\cB$-icons.} 

Suppose we have weak $\cB$-functors
\[F,F'\: X \tra X'\]
\[G,G'\: Y \tra Y'\]
agreeing on objects, and $\cB$-icons
\[
\psset{unit=0.1cm,labelsep=2pt,nodesep=2pt}
\pspicture(20,20)

\rput(0,10){\rnode{a1}{$X$}}  
\rput(20,10){\rnode{a2}{$X'$}}  

\ncarc[arcangle=45]{->}{a1}{a2}\naput{{\scriptsize $F$}}
\ncarc[arcangle=-45]{->}{a1}{a2}\nbput{{\scriptsize $F'$}}

\pcline[linewidth=0.6pt,doubleline=true,arrowinset=0.6,arrowlength=0.8,arrowsize=0.5pt 2.1]{->}(10,13)(10,7)  \naput{{\scriptsize $\alpha$}}

\endpspicture
\]
\[
\psset{unit=0.1cm,labelsep=2pt,nodesep=2pt}
\pspicture(20,20)

\rput(0,10){\rnode{a1}{$Y$}}  
\rput(20,10){\rnode{a2}{$Y'$}}  

\ncarc[arcangle=45]{->}{a1}{a2}\naput{{\scriptsize $G$}}
\ncarc[arcangle=-45]{->}{a1}{a2}\nbput{{\scriptsize $G'$}}

\pcline[linewidth=0.6pt,doubleline=true,arrowinset=0.6,arrowlength=0.8,arrowsize=0.5pt 2.1]{->}(10,13)(10,7)  \naput{{\scriptsize $\beta$}}

\endpspicture
\]
we define the $\cB$-icon
\[
\psset{unit=0.1cm,labelsep=2pt,nodesep=2pt}
\pspicture(20,20)

\rput(-3,10){\rnode{a1}{$X \otimes Y$}}  
\rput(28,10){\rnode{a2}{$X' \otimes Y'$}}  

\ncarc[arcangle=45]{->}{a1}{a2}\naput{{\scriptsize $F \otimes G$}}
\ncarc[arcangle=-45]{->}{a1}{a2}\nbput{{\scriptsize $F' \otimes G'$}}

\pcline[linewidth=0.6pt,doubleline=true,arrowinset=0.6,arrowlength=0.8,arrowsize=0.5pt 2.1]{->}(10,13)(10,7)  \naput{{\scriptsize $\alpha \otimes \beta$}}

\endpspicture
\]
to have components

\[
\psset{unit=0.1cm,labelsep=2pt,nodesep=2pt}
\pspicture(0,-4)(30,24)

\rput(-10,10){\rnode{a1}{$(x_0,x_1)(y_0,y_1)$}}  
\rput(50,14){\rnode{a2}{$(Fx_0,Fx_1)(Gy_0,Gy_1)$}}  
\rput(51,10){\rotatebox{90}{$=$}}
\rput(50,6){\rnode{a3}{$(F'x_0,F'x_1)(G'y_0,G'y_1)$}}  

\ncarc[arcangle=45]{->}{a1}{a2}\naput[npos=0.55]{{\scriptsize $F.G$}}
\ncarc[arcangle=-45]{->}{a1}{a3}\nbput[npos=0.55]{{\scriptsize $F'.G'$}}

{
\rput[c](24,7){\psset{unit=1mm,doubleline=true,arrowinset=0.6,arrowlength=0.5,arrowsize=0.5pt 2.1,nodesep=0pt,labelsep=2pt}
\pcline{->}(0,5)(0,0) \nbput{{\scriptsize $\alpha_{x_0x_1}.\beta_{y_0y_1}$}}}}

\endpspicture
\]
and then the axioms for a $\cB$-icon follow easily from naturality of $R$ and the icon axioms for $\alpha$ and $\beta$ tensored together.

\bigskip

This completes the definition of the weak functor 
\[\otimes \: \Bicon \times \Bicon \tra \Bicon\]  
on cells.  Note that it follows immediately that this is locally functorial on hom-categories because this amounts to the functoriality of the tensor product of $\cB$.  We now define composition and unit constraints for the weak functoriality of $\otimes$.   The principle is as before.  We define putative $\cl{B}$-icons whose components are the corresponding constraints in \cl{B}; as axioms are then inherited from \cl{B},  we just have to check that these components really form $\cl{B}$-icons.  These calculations are routine, and are generally similar to those in Definition~\ref{iconcomp} and Proposition~\ref{threepointseven}, so we leave the details to the reader.

\subsubsection*{Definition of composition constraint for the functor $\otimes$.} 

Given weak \cl{B}-functors
\[X \tmap{F} X' \tmap{F'} X'' \]
\[Y \tmap{G} Y' \tmap{G'} Y'' \]
we need a constraint \cl{B}-icon
\[
\psset{unit=0.1cm,labelsep=2pt,nodesep=3pt}
\pspicture(30,20)

\rput(0,10){\rnode{a1}{$X \otimes Y$}}  
\rput(30,10){\rnode{a2}{$X'' \otimes Y''$.}}  

\ncarc[arcangle=45]{->}{a1}{a2}\naput{{\scriptsize $(F' \otimes G')\circ (F \otimes G)$}}
\ncarc[arcangle=-45]{->}{a1}{a2}\nbput{{\scriptsize $(F'\circ F) \otimes (G'\circ G)$}}

{
\rput[c](15,8.5){\psset{unit=1mm,doubleline=true,arrowinset=0.6,arrowlength=0.5,arrowsize=0.5pt 2.1,nodesep=0pt,labelsep=2pt}
\pcline{->}(0,3)(0,0) \naput{{\scriptsize $$}}}}

\endpspicture
\]
First observe that the source and target \cl{B}-functors agree on objects giving
\[(x,y) \tmapsto (F'Fx, G'Gy)\]
so it makes sense to seek such a \cl{B}-icon.  Now on homs the above \cl{B}-functors give us 1-cells in \cl{B}

\[
\psset{unit=0.1cm,labelsep=2pt,nodesep=3pt}
\pspicture(80,20)


\rput(0,12){\rnode{a1}{$X(x_0,x_1) $}}
\rput(37,12){\rnode{a2}{$X'(Fx_0, Fx_1)$}}
\rput(80,12){\rnode{a3}{$X''(F'Fx_0, F'Fx_1)$}}

\rput(0,0){\rnode{b1}{$Y(y_0,y_1)$}}
\rput(37,0){\rnode{b2}{$Y'(Gy_0, Gy_1)$}}
\rput(80,0){\rnode{b3}{$Y''(G'Gy_0, G'Gy_1)$}}

\ncline{->}{a1}{a2} \naput{{\scriptsize $F_{x_0,x_1}$}}
\ncline{->}{a2}{a3} \naput{{\scriptsize $F'_{Fx_0, Fx_1}$}}

\ncline{->}{b1}{b2} \naput{{\scriptsize $G_{y_0,y_1}$}}
\ncline{->}{b2}{b3} \naput{{\scriptsize $G'_{Gy_0, Gy_1}$}}

\endpspicture\]
so our constraint \cl{B}-icon needs components that are 2-cells in \cl{B} of the form

\[
\psset{unit=0.1cm,labelsep=3pt,nodesep=3pt}
\pspicture(0,-10)(80,40)

\rput(0,10){\rnode{a1}{$X(x_0,x_1). Y(y_0,y_1)$}}  
\rput(70,10){\rnode{a2}{\makebox[8em][l]{$X''(F'Fx_0,F'Fx_1).  Y''(G'Gy_0,G'Gy_1)$}}}

\ncarc[arcangle=40]{->}{a1}{a2}\naput{{\scriptsize $(F'_{Fx_0,Fx_1} . G'_{Gy_0,Gy_1})\circ(F_{x_0,x_1} . G_{y_0,y_1})$}}
\ncarc[arcangle=-40]{->}{a1}{a2}\nbput{{\scriptsize $(F'_{Fx_0,Fx_1}\circ F_{x_0,x_1}) . (G'_{Gy_0,Gy_1}\circ G_{y_0,y_1})$}}

{
\rput[c](35,8.5){\psset{unit=1mm,doubleline=true,arrowinset=0.6,arrowlength=0.5,arrowsize=0.5pt 2.1,nodesep=0pt,labelsep=0pt}
\pcline{->}(0,4)(0,0) \naput{{\scriptsize $$}}}}

\endpspicture
\]
so we can simply use the constraint for the functoriality of the tensor product in \cl{B} (recall that in all these diagrams the tensor product in $\cB$ is written as a dot).     We must check this satisfies the axioms for a \cl{B}-icon; as before this follows from coherence in \cl{B}.


\subsubsection*{Definition of unit constraint for the functor $\otimes$.} 

Recall that for all objects $X$ in \cl{B} we have an identity weak \cl{B}-functor $1_X$.  We need for all $X,Y$ a coherence constraint
\[ 1_{X \otimes Y} \tra 1_X \otimes 1_Y.\]
As before, observe that the source and target functors agree on objects, both acting as the identity.  So our constraint \cl{B}-icon needs components that are 2-cells in \cl{B} of the form
\[
\psset{unit=0.1cm,labelsep=3pt,nodesep=3pt}
\pspicture(0,-10)(50,30)

\rput(0,10){\rnode{a1}{$X(x_0,x_1). Y(y_0,y_1)$}}  
\rput(50,10){\rnode{a2}{$X(x_0,x_1). Y(y_0,y_1).$}}

\ncarc[arcangle=40]{->}{a1}{a2}\naput{{\scriptsize $1_{X(x_0,x_1) . Y(y_0,y_1)}$}}
\ncarc[arcangle=-40]{->}{a1}{a2}\nbput{{\scriptsize $1_{X(x_0,x_1)} . 1_{Y(y_0,y_1)}$}}

{
\rput[c](25,8.5){\psset{unit=1mm,doubleline=true,arrowinset=0.6,arrowlength=0.5,arrowsize=0.5pt 2.1,nodesep=0pt,labelsep=0pt}
\pcline{->}(0,4)(0,0) \naput{{\scriptsize $$}}}}

\endpspicture
\]
Again we use the constraint for the tensor product in \cl{B}.   We must check this satisfies the axioms for a \cl{B}-icon, and again this follows from coherence in \cl{B}.

\subsubsection*{Functoriality axioms for $\otimes$.}

It suffices to check the axioms componentwise.  Note that all the constraints have components given by constraint cells for the functoriality of the tensor product of $\cB$ as a functor
\[\cB \times \cB \lra \cB.\]
Hence by coherence for functors of bicategories, the axioms hold.  

Thus we have defined a weak functor
\[\otimes \: \Bicon \times \Bicon \tra \Bicon.\]

\subsubsection*{Unit object for the monoidal structure on \Bicon.}

The unit object $I$ in \Bicon\ has
\begin{itemize}

\item one object
\item the hom object is the unit in \cl{B}
\item composition is given by $l$.

\end{itemize}

All the rest of the structure is given by coherence.  Note that we could instead use $r$ for composition, which would give an equivalent object of \Bicon.

\subsubsection*{Constraints $a,l,r$ for the monoidal structure on \Bicon.}

We must define an adjoint equivalence $(a, \eta, \epsilon)$
\[
\psset{unit=0.1cm,labelsep=2pt,nodesep=3pt}
\pspicture(30,20)

\rput(0,10){\rnode{a1}{$(X \otimes Y) \otimes Z$}} 
\rput(30,10){\rnode{a2}{$X \otimes (Y \otimes Z)$}}  

\ncarc[arcangle=45]{->}{a1}{a2}\naput{{\scriptsize $a$}}
\ncarc[arcangle=45]{->}{a2}{a1}\naput{{\scriptsize $a^*$}}

\endpspicture
\]
where $a$ and $a^*$ are weak \cl{B}-functors, and $\eta$ and $\epsilon$ must be \cl{B}-icons.

\begin{itemize}
 \item On objects $a$ and $a^*$ are given by associativity in \Set.
\item On homs $a$ and $a^*$ are given by the adjoint equivalences for associativity of the tensor product in \cl{B}.
\end{itemize}

We show that this makes $a$ into a weak \cl{B}-functor.  First we must define constraints for $a$.  The composition constraint is given by the composite below. Note that here $X,Y,Z$ are shorthand for hom-objects of the weak $\cB$-categories, with 0-cells suppressed.  Also note that we have suppressed the associators necessary to make the given braidings applicable, as these are taken care of by coherence for the braiding.

\[\psset{unit=0.08cm,labelsep=0pt,nodesep=3pt}
\pspicture(0,-5)(90,65)


\rput(0,60){\rnode{z1}{{$\big((XY)Z\big)\big((XY)Z\big)$}}} 
\rput(70,60){\rnode{z2}{{$\big(X(YZ)\big)\big(X(YZ)\big)$}}} 

\rput(0,40){\rnode{a1}{{$\big((XY)(XY)\big)(ZZ)$}}} 
\rput(70,40){\rnode{a2}{{$(XX)\big((YZ)(YZ)\big)$}}} 

\rput(0,20){\rnode{b1}{{$\big((XX)(YY)\big)(ZZ)$}}}   
\rput(70,20){\rnode{b2}{{$(XX)\big((YY)(ZZ)\big)$}}}  

\rput(0,0){\rnode{c1}{{$(XY)Z$}}}   
\rput(70,0){\rnode{c2}{{$X(YZ)$}}}  

\psset{nodesep=3pt,labelsep=2pt,arrows=->}
\ncline{z1}{z2}\naput{{\scriptsize $aa$}} 
\ncline{b1}{b2}\naput{{\scriptsize $a$}} 
\ncline{c1}{c2}\nbput{{\scriptsize $a$}} 

\ncline{z1}{a1}\nbput{{\scriptsize $1R1$}} 
\ncline{a1}{b1}\nbput{{\scriptsize $1R111$}} 
\ncline{b1}{c1}\nbput{{\scriptsize $mmm$}} 

\ncline{z2}{a2}\naput{{\scriptsize $1R1$}} 
\ncline{a2}{b2}\naput{{\scriptsize $111R1$}} 
\ncline{b2}{c2}\naput{{\scriptsize $mmm$}} 

\psset{labelsep=1.5pt}


\rput[c](35,11){{\scriptsize $\iso$}}
\rput[c](35,8){{\scriptsize\textsf{naturality of $a$} }}

\rput[c](35,48){{\scriptsize $\iso$}}
\rput[c](35,45){{\scriptsize\textsf{braid coherence}}}

\endpspicture\]

The unit constraint is defined similarly using coherence cells for the monoidal structure of $\cB$ and naturality of associator for the tensor product of $\cB$. The axioms for a weak \cl{B}-functor follow from coherence in \cl{B}.  The result for $a^*$ follows similarly.

Next we define $\eta$ and $\epsilon$, the unit and counit for the adjoint equivalence.  As usual, these are $\cB$-icons with components given by the corresponding units and counits in $\cB$; all axioms are given by coherence in \cl{B}.

The constraints $l$ and $r$ follow in exactly the same way.

\subsubsection*{Constraints $\pi, \mu, \lambda, \rho$ for the monoidal structure on $\cB$-Icon.}

We now need the following constraints for a monoidal bicategory: $\pi, \mu, \lambda, \rho$.  In each case the components are the corresponding constraints from \cl{B}, and the \cl{B}-icon axioms follow from coherence in the usual manner. 

\bigskip

Finally the monoidal bicategory axioms come from the monoidal bicategory axioms for \cl{B} together with the monoidal category axioms for \Set.  This completes the proof that \Bicon\ is a monoidal bicategory.  

\subsubsection*{Definition of symmetric structure.} 

We now construct a braiding $R\: X \otimes Y \tra Y \otimes X$ for \Bicon.

\begin{itemize}
 \item On objects $R$ is the symmetry in \Set.
\item On homs $R$ is the symmetry in \cl{B}.
\end{itemize}

It has the structure of a weak functor with constraints coming from the syllepsis in \cl{B} and some coherence cells as follows.   

\[\psset{unit=0.1cm,labelsep=0pt,nodesep=3pt}
\pspicture(-10,0)(30,42)


\rput(-5,40){\rnode{a1}{$XYXY$}} 
\rput(25,40){\rnode{a2}{$YXYX$}} 

\rput(-5,20){\rnode{b1}{$XXYY$}}   
\rput(25,20){\rnode{b2}{$YYXX$}}  

\rput(-5,0){\rnode{c1}{$XY$}}   
\rput(25,0){\rnode{c2}{$YX$}}  

\psset{nodesep=3pt,labelsep=2pt,arrows=->}
\ncline{a1}{a2}\naput{{\scriptsize $RR$}} 
\ncline{b1}{b2}\nbput{{\scriptsize $R$}} 
\ncline{c1}{c2}\nbput{{\scriptsize $R$}} 

\ncline{a1}{b1}\nbput{{\scriptsize $1R1$}} 
\ncline{b1}{c1}\nbput{{\scriptsize $mm$}} 

\ncline{a2}{b2}\naput{{\scriptsize $1R1$}} 
\ncline{b2}{c2}\naput{{\scriptsize $mm$}} 

\psset{labelsep=1.5pt}

\rput[c](10,32){{\scriptsize $\iso$}}
\rput[c](10,29){{\scriptsize\textsf{symmetry} }}

\rput[c](10,12){{\scriptsize $\iso$}}
\rput[c](10,9){{\scriptsize\textsf{naturality of $R$} }}

\endpspicture\]
Note that the top square is given as follows.  The 1-cell source and target have the following underlying braids respectively.


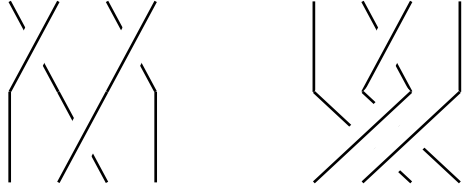
\[\psset{unit=1cm,arrowsize=1.5pt 2,arrowlength=1.4,linewidth=1pt,labelsep=5pt,nodesep=0pt}
 \psset{unit=0.08cm}
\pspicture(50,30)


\rput(0,30){\rnode{a1}{$$}}
\rput(8,30){\rnode{a2}{$$}}
\rput(16,30){\rnode{a3}{$$}}
\rput(24,30){\rnode{a4}{$$}}

\rput(0,15){\rnode{b1}{$$}}
\rput(8,15){\rnode{b2}{$$}}
\rput(16,15){\rnode{b3}{$$}}
\rput(24,15){\rnode{b4}{$$}}

\rput(0,0){\rnode{c1}{$$}}
\rput(8,0){\rnode{c2}{$$}}
\rput(16,0){\rnode{c3}{$$}}
\rput(24,0){\rnode{c4}{$$}}

\psset{nodesep=0pt,angleB=90,angleA=270,border=5.8pt}

\ncline{a1}{b2} 
\ncline{b2}{c3}

\ncline{a3}{b4} 
\ncline{b4}{c4}

\ncline{a2}{b1} 
\ncline{b1}{c1}

\ncline{a4}{b3} 
\ncline{b3}{c2}

\endpspicture
\pspicture(24,30)


\rput(0,30){\rnode{a1}{$$}}
\rput(8,30){\rnode{a2}{$$}}
\rput(16,30){\rnode{a3}{$$}}
\rput(24,30){\rnode{a4}{$$}}

\rput(0,15){\rnode{b1}{$$}}
\rput(8,15){\rnode{b2}{$$}}
\rput(16,15){\rnode{b3}{$$}}
\rput(24,15){\rnode{b4}{$$}}

\rput(0,0){\rnode{c1}{$$}}
\rput(8,0){\rnode{c2}{$$}}
\rput(16,0){\rnode{c3}{$$}}
\rput(24,0){\rnode{c4}{$$}}

\psset{nodesep=0pt,angleB=90,angleA=270,border=5.8pt}

\ncline{a2}{b3}
\ncline{a3}{b2}
\ncline{a1}{b1} 
\ncline{a4}{b4}

\ncline{b1}{c3}
\ncline{b2}{c4}
\ncline{b3}{c1}
\ncline{b4}{c2}

\endpspicture\]
\noi Comparing these braids, we see that there is an isomorphism between them using the syllepsis, and the isomorphism is unique by coherence for symmetric monoidal bicategories.  

The unit constraint has the following form
\[
\psset{unit=0.1cm,labelsep=2pt,nodesep=3pt,npos=0.4}
\pspicture(20,22)


\rput(0,20){\rnode{a1}{$\1$}}  
\rput(26,20){\rnode{a2}{$XY$}}  
\rput(26,0){\rnode{a3}{$YX$}}  

\ncline{->}{a1}{a2} \naput[npos=0.6]{{\scriptsize $I_XI_Y$}} 
\ncline{->}{a1}{a3} \nbput[npos=0.55]{{\scriptsize $I_YI_X$}} 
\ncline{->}{a2}{a3} \naput{{\scriptsize $R$}} 

{
\rput[c](19,13){{\scriptsize $\iso$}}}

\endpspicture
\]
and is given by a unique isomorphism from the braided structure.  Using the cubical notation as in Definition~\ref{tractable} the first axiom for a weak $\cl{B}$-functor is

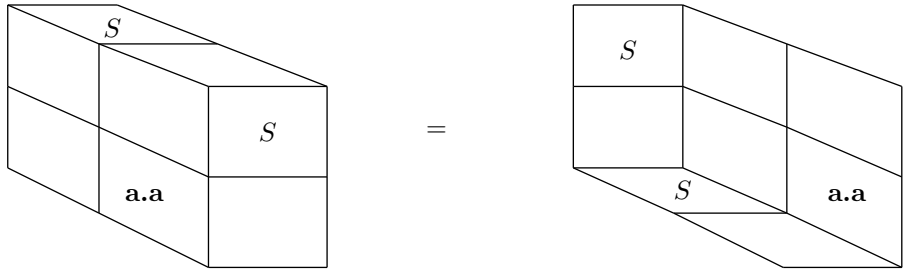
\[\psset{unit=0.12cm,labelsep=2pt,nodesep=0pt}
\pspicture(0,8)(60,40)

\rput(23,18){\mbox{\bf{a.a}}}  
\rput(36.5,25){$S$}  
\rput(19.5,36.5){$S$}

\rput(30,10){\rnode{a}{$$}}  
\rput(43,10){\rnode{b}{$$}}  
\rput(30,20){\rnode{d}{$$}}  
\rput(43,20){\rnode{e}{$$}}  
\rput(30,30){\rnode{g}{$$}}  
\rput(43,30){\rnode{h}{$$}}  
\rput(18,16){\rnode{j}{$$}}  
\rput(18,25.5){\rnode{k}{$$}}  
\rput(18,34.7){\rnode{l}{$$}}  
\rput(31,34.7){\rnode{m}{$$}}  
\rput(8,21){\rnode{p}{$$}}  
\rput(8,30){\rnode{q}{$$}}  
\rput(8,39){\rnode{r}{$$}}  
\rput(20,39){\rnode{s}{$$}}  

\ncline{a}{b} 
\ncline{d}{e} 
\ncline{g}{h} 

\ncline{a}{d} 
\ncline{d}{g} 
\ncline{b}{e} 
\ncline{e}{h} 
 
\ncline{j}{k} 
\ncline{k}{l} 
\ncline{l}{m} 
\ncline{j}{a} 
\ncline{k}{d} 
\ncline{l}{g} 
\ncline{m}{h} 
 
\ncline{p}{q} 
\ncline{q}{r} 
\ncline{r}{s} 
\ncline{p}{j} 
\ncline{q}{k} 
\ncline{r}{l} 
\ncline{s}{m} 

\endpspicture
%
%
%
%
\pspicture(10,8)(60,40)

\rput(5,25){$=$}  

\rput(32,18.5){$S$}  
\rput(26,34){$S$}  
\rput(50,18){\mbox{\bf{a.a}}}

\rput(43,10){\rnode{b}{$$}}  
\rput(56,10){\rnode{c}{$$}}  
\rput(31,16){\rnode{d}{$$}}  
\rput(43.4,16){\rnode{e}{$$}}  
\rput(56,20){\rnode{f}{$$}}  
\rput(32,21){\rnode{g}{$$}}  
\rput(43.4,25.5){\rnode{h}{$$}}  
\rput(56,30){\rnode{i}{$$}}  
\rput(20,21){\rnode{k}{$$}}  
\rput(20,30){\rnode{l}{$$}}  
\rput(32,30){\rnode{m}{$$}}  
\rput(43.4,34.7){\rnode{n}{$$}}  
\rput(20,39){\rnode{s}{$$}}  
\rput(32,39){\rnode{t}{$$}}

\ncline{b}{c} 
\ncline{d}{b} 
\ncline{e}{c} 
\ncline{d}{e} 
\ncline{h}{f} 

\ncline{g}{e} 
\ncline{k}{d} 
\ncline{h}{e} 
\ncline{g}{k} 
\ncline{c}{f} 
\ncline{f}{i}
 
\ncline{k}{l} 
\ncline{l}{m} 
\ncline{m}{t} 
\ncline{m}{g} 
\ncline{n}{h} 
\ncline{n}{i}
 
\ncline{s}{t} 
\ncline{s}{l} 
\ncline{h}{m} 
\ncline{t}{n}

\endpspicture
\]
where $S$ indicates a unique coherence 2-cell involving the syllepsis of \cB, and as before \mbox{\bf{a}} indicates the associativity constraint for a weak \cB-category.  The equality then follows from a calculation manipulating cubes in a similar fashion to earlier, using coherence for symmetric monoidal bicategories.   The other axiom follows similarly.

Now we need a unit and counit.  These come from the adjoint equivalence for $R$ in the braided structure on \cB.  We just need to check that $\eta$ and $\epsilon$ are \cB-icons; this follows from coherence again, this time for symmetric monoidal bicategories.  The modifications $R_{(-|-,-)}$ and $R_{(-,-|-)}$ are given on components by the corresponding cells in \cB.  We only need to check the icon axioms, and these follow from coherence for symmetric monoidal bicategories and the corresponding modification axioms for the braiding in \cB.  The braiding axioms follow from those in \cB.

Finally we construct a syllepsis: the components come from the syllepsis in \cB, so the axioms (and the symmetry axiom) are immediate.  We just need to check the icon axioms, and these follow from coherence for symmetric monoidal bicategories.

\bigskip

This completes the proof that \Bicon\ is a symmetric monoidal bicategory as required.   \end{prf}

\section{Higher dimensions}\label{hd}

Now that we know \Bicon\ is a symmetric monoidal bicategory, we can iterate the icon construction to produce the higher-dimensional structures that motivated this work.  Recall that putting $\cB=\Cat$ gives $\cat{Icon}$, a convenient 2-category of bicategories.  The idea is that iterating this gives us convenient 2-categories of 3-categories, 4-categories, and so on.  In fact these are not strict $n$-categories, nor are they semi-strict in the usual sense---while they can be thought of as ``$n$-categories with extra strictness conditions'', in fact the coherence constraints have different shapes.  For some purposes this might not be a critical difference, but it is crucial for our study of totalities of degenerate structures.  We are finally able to give a satisfying account of doubly degenerate tricategories, as promised in \cite{cg3}; we refer the reader to that paper for a discussion of the pitfalls of more naive approaches.

This suggests a general scheme for studying totalities of $k$-degenerate $(n+k)$-categories, and comparing them with entries in the ``Periodic Table'' of $n$-categories.

\subsection{Iteration}

The main example we have in mind starts the iteration with \Cat, but we make the following general definition of iterated icons.

\begin{mydefinition}
We define $\cat{\cl{B}-nIcon}$ for each $n \geq 0$ as follows.
\begin{itemize}
 \item $\cat{\cl{B}-0Icon} = \cl{B}$.
 \item $\cat{\cl{B}-nIcon} = \cat{(\cl{B}-(n-1)Icon)-Icon}$.
\end{itemize}

\end{mydefinition}


%
%
%

Recall that if $\cB$ is a symmetric monoidal \emph{2-category} then \Bicon\ is a symmetric monoidal 2-category.  Thus starting with $\cB=\Cat$ we get symmetric monoidal 2-categories of $n$-categories.  These are a subtly semi-strict version of $n$-category in which low-dimensional composition is strict.  We will now exhibit this explicitly for  $\cl{B}=\cat{Icon}$ with monoidal structure given by cartesian product.  This gives the 2-category \cat{Icon-Icon}.

A 0-cell $X$ is a ``weak \cat{Icon}-category'' (see Definition~\ref{weakbcat}).  We now unpack this definition to reveal what sort of 3-dimensional structure is produced.

\begin{itemize}
 \item a set $X_0$ of objects,
\item for all $a,b \in X_0$ a bicategory $X(a,b)$,

\item composition: for all $a,b,c \in X_0$ a weak functor
\[\begin{array}{ccccc}
m_{abc} &:& X(b,c) \times X(a,b) &\lra & X(a,c)\\
&& (g,f) & \mapsto & g \otimes f\\
&& (\beta,\alpha) & \mapsto & \beta \otimes \alpha,
\end{array}\]

\item identities: for all $a \in X_0$ a weak functor
\[\begin{array}{ccccc}
I_a &:& \1 & \lra & X(a,a)\\
&& \ast & \mapsto & I_a
\end{array}\]
where here $\1$ is the terminal bicategory,

\item associativity constraints: for all $a,b,c,d \in X_0$ an invertible icon

 \[
 \xy
    (-30,28)*+{ \big(X(c,d) \times  X(b,c)\big) \times  X(a,b)}="tl";
    (-30,15)*+{ X(c,d) \times  \big(X(b,c) \times  X(a,b)\big)}="ml";
        {\ar_{\cong} "tl";"ml"};
(-30,-10)*+{X(c,d) \times  X(a,c)}="tr";
    (30,28)*+{X(b,d)\times  X(a,b)}="bl";
    (30,-10)*+{X(a,d)}="br";
        {\ar_-{1 \times m_{a,b,c}} "ml";"tr"};
        {\ar^-{m_{b,c,d} \times 1} "tl";"bl"};
        {\ar^-{m_{a,b,d}} "bl";"br"};
        {\ar_-{m_{a,c,d}} "tr";"br"};
    {\ar@{=>}^{\;a_{a,b,c,d}} (8,12);(-6,2)};
 \endxy
 \]
and

 \item unit constraints: for all $a$, $b$ invertible icons
  \[\xy
  (20,20)*+{X(a,b) \times X(a,a)}="tl";
  (-20,20)*+{X(a,b) \times \1 }="tr";
  (20,0)*+{ X(a,b)}="bl";
    {\ar^{m_{a,a,b}} "tl";"bl"};
    {\ar^{1\times I_a} "tr";"tl"};
    {\ar_{\cong} "tr";"bl"};
 {\ar@{=>}_{r_{a,b}} (15,15);(9,9)};
  \endxy
 \qquad \quad
  \xy
  (-20,20)*+{X(b,b) \times X(a,b)}="tl";
  (20,20)*+{ \1 \times X(a,b)}="tr";
  (-20,0)*+{ X(a,b)}="bl";
    {\ar_{m_{a,b,b}} "tl";"bl"};
    {\ar_{I_b \times 1} "tr";"tl"};
    {\ar^{\cong} "tr";"bl"};
{\ar@{=>}^{\ell_{a,b}} (-15,15);(-9,9)};
  \endxy\]

\end{itemize}
satisfying the axioms given in Definition~\ref{weakbcat}.

The key now is that since $a,\ell$ and $r$ are icons, some additional coherence is implicit, as follows.  Each of these constraint icons has a source and target functor that must agree on objects.  So we must have, for all $a \map{f} b \map{g} c \map{h} d$
\[\begin{array}{rcl}
   (h \otimes g) \otimes f &=& h \otimes (g \otimes f) \\
f \otimes I_a &=& f \\
I_b \otimes f &=& f
  \end{array}\]
that is, composition of 1-cells in the a weak \cat{Icon}-category is strict.  

Then the 2-cell components of the constraint icons give the following data:

\begin{itemize}
 \item  for all composable triples of 2-cells as shown
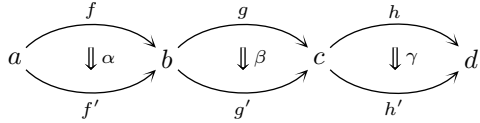
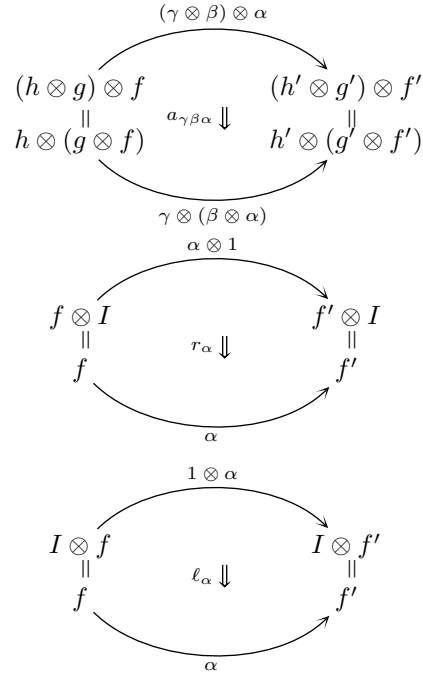
\[
\psset{unit=0.1cm,labelsep=2pt,nodesep=2pt}
\pspicture(50,20)


\rput(0,10){\rnode{a1}{$a$}}  
\rput(20,10){\rnode{a2}{$b$}}  
\rput(40,10){\rnode{a3}{$c$}}
\rput(60,10){\rnode{a4}{$d$}}

\ncarc[arcangle=45]{->}{a1}{a2}\naput{{\scriptsize $f$}}
\ncarc[arcangle=-45]{->}{a1}{a2}\nbput{{\scriptsize $f'$}}

\ncarc[arcangle=45]{->}{a2}{a3}\naput{{\scriptsize $g$}}
\ncarc[arcangle=-45]{->}{a2}{a3}\nbput{{\scriptsize $g'$}}

\ncarc[arcangle=45]{->}{a3}{a4}\naput{{\scriptsize $h$}}
\ncarc[arcangle=-45]{->}{a3}{a4}\nbput{{\scriptsize $h'$}}

{
\rput[c](10,8.5){\psset{unit=1mm,doubleline=true,arrowinset=0.6,arrowlength=0.5,arrowsize=0.5pt 2.1,nodesep=0pt,labelsep=2pt}
\pcline{->}(0,3)(0,0) \naput{{\scriptsize $\alpha$}}}}

{
\rput[c](30,8.5){\psset{unit=1mm,doubleline=true,arrowinset=0.6,arrowlength=0.5,arrowsize=0.5pt 2.1,nodesep=0pt,labelsep=2pt}
\pcline{->}(0,3)(0,0) \naput{{\scriptsize $\beta$}}}}

{
\rput[c](50,8.5){\psset{unit=1mm,doubleline=true,arrowinset=0.6,arrowlength=0.5,arrowsize=0.5pt 2.1,nodesep=0pt,labelsep=2pt}
\pcline{->}(0,3)(0,0) \naput{{\scriptsize $\gamma$}}}}

\endpspicture
\]
invertible 3-cells
\[
\psset{unit=0.1cm,labelsep=2pt,nodesep=2pt}
\pspicture(0,-4)(30,24)


\rput(0,14){\rnode{a1}{$(h \otimes g) \otimes f$}} 
\rput(1,10){\rotatebox{90}{$=$}}
\rput(0,7){\rnode{a2}{$h \otimes (g \otimes f)$}}  

\rput(35,14){\rnode{b1}{$(h' \otimes g') \otimes f'$}}  
\rput(36,10){\rotatebox{90}{$=$}}
\rput(35,7){\rnode{b2}{$h' \otimes (g' \otimes f')$}}  

\ncarc[arcangle=45]{->}{a1}{b1}\naput[npos=0.5]{{\scriptsize $(\gamma \otimes \beta) \otimes \alpha$}}
\ncarc[arcangle=-45]{->}{a2}{b2}\nbput[npos=0.5]{{\scriptsize $\gamma \otimes (\beta \otimes \alpha)$}}

{
\rput[c](19,8.5){\psset{unit=1mm,doubleline=true,arrowinset=0.6,arrowlength=0.5,arrowsize=0.5pt 2.1,nodesep=0pt,labelsep=2pt}
\pcline{->}(0,3)(0,0) \nbput{{\scriptsize $a_{\gamma\beta\alpha}$}}}}

\endpspicture
\]
\[
\psset{unit=0.1cm,labelsep=2pt,nodesep=2pt}
\pspicture(0,-4)(30,24)


\rput(0,14){\rnode{a1}{$f \otimes I$}} 
\rput(1,11){\rotatebox{90}{$=$}}
\rput(0,7){\rnode{a2}{$f$}}  

\rput(35,14){\rnode{b1}{$f' \otimes I$}}  
\rput(36,11){\rotatebox{90}{$=$}}
\rput(35,7){\rnode{b2}{$f'$}}  

\ncarc[arcangle=45]{->}{a1}{b1}\naput[npos=0.5]{{\scriptsize $\alpha \otimes 1$}}
\ncarc[arcangle=-45]{->}{a2}{b2}\nbput[npos=0.5]{{\scriptsize $\alpha$}}

{
\rput[c](19,8.5){\psset{unit=1mm,doubleline=true,arrowinset=0.6,arrowlength=0.5,arrowsize=0.5pt 2.1,nodesep=0pt,labelsep=2pt}
\pcline{->}(0,3)(0,0) \nbput{{\scriptsize $r_\alpha$}}}}

\endpspicture
\]
\[
\psset{unit=0.1cm,labelsep=2pt,nodesep=2pt}
\pspicture(0,-4)(30,24)


\rput(0,14){\rnode{a1}{$I \otimes f$}} 
\rput(1,11){\rotatebox{90}{$=$}}
\rput(0,7){\rnode{a2}{$f$}}  

\rput(35,14){\rnode{b1}{$I \otimes f'$}}  
\rput(36,11){\rotatebox{90}{$=$}}
\rput(35,7){\rnode{b2}{$f'$}}  

\ncarc[arcangle=45]{->}{a1}{b1}\naput[npos=0.5]{{\scriptsize $1 \otimes \alpha$}}
\ncarc[arcangle=-45]{->}{a2}{b2}\nbput[npos=0.5]{{\scriptsize $\alpha$}}

{
\rput[c](19,8.5){\psset{unit=1mm,doubleline=true,arrowinset=0.6,arrowlength=0.5,arrowsize=0.5pt 2.1,nodesep=0pt,labelsep=2pt}
\pcline{->}(0,3)(0,0) \nbput{{\scriptsize $\ell_\alpha$}}}}

\endpspicture
\]
natural in all arguments.

%
%
%

\end{itemize}

It is straightforward to characterise the 1- and 2-cells in \cat{Icon-Icon} in a similar fashion.

\begin{remarks}  It is tempting to think that a weak \cat{Icon}-category is a special, slightly strict kind of tricategory.  While this is a reasonable way to think about it informally, it is not strictly true.  This corresponds to and arises from the fact that we say icons are ``identity component oplax natural transformations''  \emph{informally} although this is also not strictly true.  The crucial point here is that the identity components are omitted from the data, so the 2-cell data have different shapes, yielding different, stricter, composition.  Thus in a weak \cat{Icon}-category, this phenomenon applies to the constraint icons.  

Informally, we may regard a weak \cat{Icon}-category as a ``tricategory in which 1-cell composition is strict'' as in \cite{shul2}, where the notion is called ``iconic tricategory''.  However, technically this would give a different shape of constraint cell.  For example, an associator 3-cell would look like

\[\psset{unit=0.11cm,labelsep=0pt,nodesep=3pt}
\pspicture(40,22)


\rput(0,20){\rnode{a1}{$(h \otimes g) \otimes f$}} 
\rput(40,20){\rnode{a2}{$(h' \otimes g') \otimes f'$}} 

\rput(0,0){\rnode{b1}{$h \otimes (g \otimes f)$}}   
\rput(40,0){\rnode{b2}{$h' \otimes (g' \otimes f')$}}  

\psset{nodesep=3pt,labelsep=2pt,arrows=->}
\ncline{a1}{a2}\naput{{\scriptsize $(\gamma \otimes \beta) \otimes \alpha$}} 
\ncline{b1}{b2}\nbput{{\scriptsize $\gamma \otimes (\beta \otimes \alpha)$}} 
\ncline{a1}{b1}\nbput{{\scriptsize $I$}} 
\ncline{a2}{b2}\naput{{\scriptsize $I$}} 

\psset{labelsep=1.5pt}
\pnode(23,13){a3}
\pnode(17,7){b3}
\ncline[doubleline=true,arrowinset=0.6,arrowlength=0.8,arrowsize=0.5pt 2.1]{a3}{b3} \nbput[npos=0.4]{{\scriptsize $a_{\gamma\beta\alpha}$}}

\endpspicture\]
instead of the globular shape above, and then all the axioms involving it would be much more complicated to take into account the unit 1-cells at the sides.  The correct statement should be that there is a bicategory of iconic tricategories and appropriate higher cells between them, which is biequivalent to \cat{Icon-Icon}.  Shulman \cite{shul2} alludes to this but does not attempt the technicalities of producing the correct higher cells and their composition.  

A weak \cat{Icon}-category is not only stricter than a general tricategory but also more streamlined.  This is crucially what enables us to cut down the dimensions of the totality, cleanly describe the correct higher cells, and produce the correct degenerate structures corresponding to those in the periodic table.  

Note, however, that if we consider icons between 2-categories and 2-functors (so everything is strict), then these issues go away: icons in this case really are identity component oplax natural transformations.  Similarly \cat{Gray}-categories  really are a strict form of weak \cat{Icon}-category---those for which 
\begin{itemize}
 \item the homs are 2-categories,
\item the composition and unit functors are cubical, and
\item the icons $a,\ell, r$ are identities. 
\end{itemize}
\end{remarks}

%
%

%
%

Recall that in the usual tricategory of bicategories, 1-cells compose strictly.  In fact, a naturally occurring example of a (large) weak \cat{Icon}-category is given by the totality of bicategories, as follows.

\begin{theorem}
There is a large weak \cat{Icon}-category with
\begin{itemize}
 \item 0-cells all (small) bicategories, and
\item given bicategories $X, Y$, $\cat{Bicat}(X,Y)$ is the bicategory of weak functors, weak transformations and modifications from $X$ to $Y$. 
\end{itemize}

\end{theorem}

\begin{prf}
We must give weak functors for units and composition, and icons for associativity and unit constraints.  It is not difficult to modify the proofs given in \cite[Section 5.1]{gur3} by omitting the relevant identity components.
\end{prf}

%
%
%
%
%
%
%
%

\begin{example}
 The category of bicategories and weak functors is monoidal under the cartesian product.  A category enriched in this category can be viewed as a weak \cat{Icon}-category in which $a, \ell$ and $r$ are identity icons.  Thus not only is 1-cell composition strict, but all composition along bounding 0-cells is strict.  An example of such a structure is the tricategory of topological spaces as constructed in \cite{gur4}. 
\end{example}

\subsection{Degeneracy}

%
%
%
%
%

We now show that \cat{Icon-Icon} is the correct framework for studying doubly degenerate 3-dimensional structures.  That is, restricting \cat{Icon-Icon} to the doubly degenerate 0-cells gives us the 2-category of braided monoidal categories.

In fact, a careful study must proceed via the intermediate structure of ``2-tuply monoidal categories''.  Essentially, these are categories with two monoidal structures and interchange up to isomorphism; precisely, they are pseudomonoids in the monoidal 2-category of monoidal categories and weak functors.  

Joyal and Street state that 2-tuply monoidal categories are the same as braided monoidal categories \cite[Proposition 5.3 and Example 5.4]{js1} but do not give full details about the totalities of such structures.  

As we saw in \cite{cg3}, passing between these two types of structure is quite delicate in a way that has a serious impact on doubly degenerate case. For this reason we separate our comparison into two distinct steps.

\begin{enumerate}
 \item Compare doubly degenerate weak \cat{Icon}-categories with 2-tuply monoidal categories.
\item Compare 2-tuply monoidal categories with braided monoidal categories.
\end{enumerate}

The present work makes the first step immediate, resolving the difficulties that arise from more naive approaches.  For the second we refer to \cite{js1}.

\begin{theorem}
Let \cl{J} be the full sub-2-category of \cat{Icon-Icon} whose 0-cells are the doubly degenerate weak \cat{Icon}-categories.  That is, the weak \cat{Icon}-categories with only one object, and whose single hom-bicategory has only one 0-cell.  Let \cat{2MonCat} be the 2-category of pseudomonoids in \cat{MonCat}.
Then there is a 2-equivalence $\cl{J} \catequiv \cat{2MonCat}$.
\end{theorem}

\begin{prf}
We characterise \cl{J} explicitly and will see that the only difference between \cl{J} and \cat{2MonCat} is in a choice of unique 0- and 1-cell.

An object of \cl{J} consists of
\begin{itemize}
 \item a monoidal category $M$ with unit $I$, say,
\item a monoidal functor 
\[\begin{array}{ccc}
   \1 & \tmap{U} & M\\
\ast & \tmapsto & u,
  \end{array}\]
\item a monoidal functor $\boxtimes \: M \times M \tra M$,
\item monoidal natural isomorphisms
\[\begin{array}{ccc}
   (x \boxtimes y) \boxtimes z & \tmap{a_{x,y,z}} & x \boxtimes (y \boxtimes z) \\
x \boxtimes u & \tmap{r_x} & x \\
u \boxtimes x & \tmap{\ell_x} & x
  \end{array}\]
\end{itemize}
and axioms making $(M,u,r,\ell,a)$ into a monoidal category.  This is exactly the structure of a 2-tuply monoidal category.

A 1-cell $F\: M \tra M'$ of \cl{J} consists of
\begin{itemize}
\item a monoidal functor $F \: M \tra M'$ (monoidal with respect to $\otimes$),
\item a monoidal isomorphism
\[\psset{unit=0.08cm,labelsep=0pt,nodesep=3pt}
\pspicture(40,22)


\rput(0,20){\rnode{a1}{$M \times M$}} 
\rput(40,20){\rnode{a2}{$M' \times M',$}} 

\rput(0,0){\rnode{b1}{$M$}}   
\rput(40,0){\rnode{b2}{$M'$}}  

\psset{nodesep=3pt,labelsep=2pt,arrows=->}
\ncline{a1}{a2}\naput{{\scriptsize $F \times F$}} 
\ncline{b1}{b2}\nbput{{\scriptsize $F$}} 
\ncline{a1}{b1}\nbput{{\scriptsize $\boxtimes$}} 
\ncline{a2}{b2}\naput{{\scriptsize $\boxtimes'$}} 

{
\rput[c](17,8){\psset{unit=1mm,doubleline=true,arrowinset=0.6,arrowlength=0.5,arrowsize=0.5pt 2.1,nodesep=0pt,labelsep=1pt}
\pcline{<-}(0,0)(3,3) \naput{{\scriptsize $\phi_2$}}}}

\endpspicture\]

\item a monoidal isomorphism
\[
\psset{unit=0.1cm,labelsep=2pt,nodesep=3pt,npos=0.4}
\pspicture(20,22)


\rput(0,20){\rnode{a1}{$\1$}}  
\rput(26,20){\rnode{a2}{$M$}}  
\rput(26,0){\rnode{a3}{$M'$}}  

\ncline{->}{a1}{a2} \naput[npos=0.6]{{\scriptsize $U$}} 
\ncline{->}{a1}{a3} \nbput[npos=0.55]{{\scriptsize $U'$}} 
\ncline{->}{a2}{a3} \naput{{\scriptsize $F$}} 

{
\rput[c](17,11){\psset{unit=1mm,doubleline=true,arrowinset=0.6,arrowlength=0.5,arrowsize=0.5pt 2.1,nodesep=0pt,labelsep=1pt}
\pcline{<-}(0,0)(3,3) \naput{{\scriptsize $\phi_0$}}}}

\endpspicture
\]
and axioms making $(F,\phi_2, \phi_0)$ into a monoidal functor with respect to $\boxtimes$ and $\boxtimes'$.

\end{itemize}

A 2-cell
\[
\psset{unit=0.1cm,labelsep=2pt,nodesep=3pt}
\pspicture(20,20)

\rput(0,10){\rnode{a1}{$M$}}  
\rput(20,10){\rnode{a2}{$M'$}}  

\ncarc[arcangle=45]{->}{a1}{a2}\naput{{\scriptsize $F$}}
\ncarc[arcangle=-45]{->}{a1}{a2}\nbput{{\scriptsize $F'$}}

{
\rput[c](10,8.5){\psset{unit=1mm,doubleline=true,arrowinset=0.6,arrowlength=0.5,arrowsize=0.5pt 2.1,nodesep=0pt,labelsep=2pt}
\pcline{->}(0,3)(0,0) \naput{{\scriptsize $\alpha$}}}}

\endpspicture
\]
of $\cl{J}$ consists of a monoidal transformation $\alpha$ with respect to $\otimes$ (because icons between degenerate bicategories are precisely monoidal transformations between the corresponding monoidal categories, together), and axioms making it monoidal with respect to $\boxtimes$.

Thus once we have forgotten the unique 0- and 1-cells from the doubly degenerate objects of $\cl{J}$, we have exactly the 0-, 1- and 2-cells of \cat{2MonCat}.
\end{prf}

Note that we could perform the icon construction again to produce the 2-category $\cat{Cat-3Icon} = \cat{Icon-Icon-Icon}$.  Triply degenerate 0-cells of this are 3-tuply monoidal categories; Joyal and Street \cite{js1} observe that these correspond to symmetric monoidal categories.

\begin{theorem}
Let \cl{K} be the full sub-2-category of \cat{Icon-Icon-Icon} whose 0-cells are the triply degenerate ones.   Let \cat{3MonCat} be the 2-category of pseudomonoids in \cat{2MonCat}.  Then there is a 2-equivalence $\cl{K} \catequiv \cat{3MonCat}$.
\end{theorem}


We now propose a general construction that should provide the correct totality of $k$-degenerate $(n+k)$-categories for every $n$ and $k$.  The idea is that \emph{a priori} the $k$-degenerate $(n+k)$-categories form an $(n+k+1)$-category, but when we regard them as $n$-categories with extra structure, they should form an $(n+1)$-category; naive methods of making this dimension change do not produce the correct structure.  However, iterating the icon construction does produce the correct structure.

%

The generalised icon construction described in this paper gives the case $n=1$. The case for $n >1 $ involves the $(n+1)$-category of weak $n$-categories, a structure that is not very well understood at present.  Nevertheless, we conjecture the following generalisations of the above constructions.

\begin{conjecture}
Given a symmetric monoidal $(n+1)$-category $\cl{E}$, we can define a symmetric monoidal $(n+1)$-category \cat{\cl{E}-Icon} of ``categories weakly enriched in \cl{E}'', weak functors, and icon-like higher morphisms.  
\end{conjecture}

Note also that the weakness of the enrichment increases as $n$ increases.  
However, the higher morphisms must be iconic only for the lowest dimension of component, indexed by 0-cells.  For example, a transformation between tricategories is in general given by
\begin{itemize}
 \item for every 0-cell in the source, a 1-cell in the target,
\item for every 1-cell in the source, a 2-cell in the target, and
\item for every 2-cell in the source, a 3-cell in the target.  
\end{itemize}
Thus the icon version omits the 1-cell data, and alters the shape of the 2- and 3-cell data.  Likewise, a modification between tricategories is in general given by
\begin{itemize}
 \item for every 0-cell in the source, a 2-cell in the target, and
\item for every 1-cell in the source, a 3-cell in the target,  
\end{itemize}
so an ``icon-like'' modification would omit the 2-cell data, and alter the shape of the 3-cell data accordingly.  For higher cells with more dimensions of icon-like behaviour, we need to iterate the construction.  We propose the following definition, although the formalism to make it precise does not yet exist.

\begin{conjecture}
Fix $n \geq 1$. We define a symmetric monoidal $(n+1)$-category $\cE_{n,k}$ for each $k \geq 0$ as follows.
\begin{itemize}
 \item $\cE_{n,0} = \cat{$n$-Cat}$, and

\item $\cE_{n,k} = \cat{$\cE_{n,k-1}$-Icon}$.
\end{itemize}
\end{conjecture}

Thus objects of $\cE_{n,k}$ are a bit like ``slightly strict'' $(n+k)$-categories.

\begin{remark}
 This indexing is convenient for considering the $k$-degenerate version, but another point of view is to fix the total dimension and ask what lower dimensional totalities of $m$-category we have, lower than the canonical $(m+1)$-category thereof.  The idea is that we should now get a $k$-category of $m$-categories for every $k \leq m$ as well, but for different values of $k$ we have different starting points in the above iteration.  Thus in each case we have a different notion of slightly strict $m$-category.  This is different from the approach of \cite{gg1} in which low-dimensional totalities of tricategories are studied but the tricategories and functors remain fully weak; constructing a tricategory of such is then very complex.    
\end{remark}

Finally we make the following conjectures for $k$-degenerate $(n+k)$-categories ones.

\begin{conjecture}
 Given a symmetric monoidal $(n+1)$-category $\cE$ there is a symmetric monoidal $(n+1)$-category $\cat{Mon}(\cE)$ of weak monoids in $\cE$.  Iterating this construction, we write
\begin{itemize}
 \item $\cat{0Mon}(\cE) = \cE$, and
\item $\cat{(k+1)Mon}(\cE) = \cat{Mon}(\cat{kMon}(\cE))$.
\end{itemize}

\end{conjecture}

\begin{conjecture}
 Write $\cJ_{n,k}$ for the full sub $(n+1)$-category of $\cE_{n,k}$ whose objects are the $k$-degenerate ones.  Then there is an equivalence of $(n+1)$-categories
\[\cJ_{n,k} \catequiv \cat{kMon}(\cat{$n$-Cat}).\]
\end{conjecture}

It remains to compare the $k$-tuply monoidal $n$-categories with the appropriate entry in the Periodic Table.



%


\ed